\documentclass{article}
\usepackage{amsmath,amsfonts,amsthm,amssymb,amscd,latexsym, bm,bbm}
\usepackage{hyperref}
\usepackage[usenames,dvipsnames]{color}
\usepackage{authblk}
\hypersetup{colorlinks,
citecolor=Blue,
linkcolor=Blue,
urlcolor=Blue}

\usepackage[T2A]{fontenc}
\usepackage{graphicx}
\usepackage{mathrsfs}

{\theoremstyle {definition} \newtheorem {defi} {Definition} [section] }
{\theoremstyle {plain}  \newtheorem {theo} [defi] {Theorem}}
{\theoremstyle {plain}  \newtheorem {cor} [defi]{Corollary}}
{\theoremstyle {plain}  \newtheorem {lem} [defi]{Lemma}}
{\theoremstyle {plain}  \newtheorem {prop} [defi]{Proposition}}
{\theoremstyle {plain} \newtheorem {rem}[defi] {Remark}}
{\theoremstyle {definition} \newtheorem {example}[defi] {Example}}

\numberwithin{equation}{section}
\newtheorem*{defi*}{Definition}
\newtheorem*{problem*}{Problem}

\newtheorem*{rem*}{Remark}
\newtheorem*{note*}{Note}

\makeatletter \@addtoreset{equation}{section} \makeatother

\newcommand{\mC}{\mathbb{C}}
\newcommand{\mR}{\mathbb{R}}

\newcommand{\mZ}{\mathbb{Z}}
\newcommand{\mN}{\mathbb{N}}
\newcommand{\mB}{\mathbb{B}}
\newcommand{\mI}{\mathbb{I}}

\newcommand{\mH}{\mathbb{H}}

\newcommand{\CC}{{\cal C}}
\newcommand{\DD}{{\cal D}}

\newcommand{\FF}{{\cal F}}
\newcommand{\GG}{{\cal G}}

\newcommand{\NN}{{\cal N}}

\newcommand{\SSS}{{\cal S}}

\newcommand{\KK}{{\cal K}}
\newcommand{\XX}{{\cal X}}

\newcommand{\eps}{\varepsilon}
\newcommand{\ph}{\varphi}

\newcommand{\be}{\beta}

\newcommand{\tht}{\theta}
\newcommand{\al}{\alpha}
\newcommand{\ga}{\gamma}

\newcommand{\de}{\delta}
\newcommand{\De}{\Delta}

\newcommand{\supp}{\operatorname{supp}}
\newcommand{\tr}{\operatorname{tr}}

\newcommand{\Ree}{\operatorname{Re}}

\newcommand{\Cov}{\operatorname{Cov}}
\newcommand{\Var}{\operatorname{Var}}

\newcommand{\ov}{\overline}
\newcommand{\wid}{\widetilde}

\newcommand{\ssk}{\smallskip}
\newcommand{\msk}{\medskip}
\newcommand{\bsk}{\bigskip}

\newcommand{\MO}{\mbox{\bf E}\,}
\newcommand{\MON}{\mbox{\bf E}_N}

\newcommand{\PR}{\mbox{\bf P}\,}

\newcommand{\ds}{\displaystyle{}}
\newcommand{\ra}{\rightarrow}
\newcommand{\ran}{\rangle}
\newcommand{\lan}{\langle}
\newcommand{\volna}{\thicksim}
\newcommand{\sm}{\setminus}
\newcommand{\raw}{\rightharpoonup}
\newcommand{\os}{\overset}

\newcommand{\fr}{\frac}

\newcommand{\qmb}{\quad\mbox}
\newcommand{\qu}{\quad}
\newcommand{\qnd}{\quad\mbox{and}\quad}

\newcommand{\sli}{\sum\limits}
\newcommand{\ili}{\int\limits}
\newcommand{\ilif}{\ili_{-\infty}^\infty}
\newcommand{\lbl}{\label}

\newcommand{\ass}{\quad\mbox{as}\quad}

\newcommand{\rprop}{Proposition \nolinebreak}
\newcommand{\rtheo}{Theorem \nolinebreak}
\newcommand{\rlem}{Lemma \nolinebreak}
\newcommand{\rcor}{Corollary \nolinebreak}

\newcommand{\rsec}{Section \nolinebreak}

\newcommand{\bee}{\begin{equation}}
\newcommand{\eee}{\end{equation}}
\newcommand{\btt}{\begin{theo}}
\newcommand{\ett}{\end{theo}}
\newcommand{\bl}{\begin{lem}}
\newcommand{\el}{\end{lem}}
\newcommand{\bpp}{\begin{prop}}
\newcommand{\epp}{\end{prop}}
\newcommand{\bcc}{\begin{cor}}
\newcommand{\ecc}{\end{cor}}
\newcommand{\bdd}{\begin{def}}
\newcommand{\edd}{\end{def}}
\newcommand{\brr}{\begin{rem}}
\newcommand{\err}{\end{rem}}

\newcommand{\non}{\nonumber}
\newcommand{\sck}{\substack}
\newcommand{\Zp}{\mathbb{Z}^d_+}
\newcommand{\sq}{\sqrt}
\newcommand{\Conf}{\mbox{Conf}}
\newcommand{\CoRr}{\Conf(\mR^m)}

\title{A functional limit theorem for the sine-process}
\author[1,2,3,4,5]{Alexander I. Bufetov\thanks{bufetov@mi.ras.ru}}
\author[2,3]{Andrey V. Dymov\thanks{dymov@mi.ras.ru}}
\affil[1]{Aix-Marseille Universit\'e, CNRS, Centrale Marseille, I2M, UMR 7373, 39 rue F. Joliot Curie, Marseille, FRANCE}
\affil[2]{Steklov Mathematical Institute of RAS, Moscow}
\affil[3]{National Research University Higher School of Economics, Moscow}
\affil[4]{Institute for Information Transmission Problems, Moscow}
\affil[5]{The Chebyshev Laboratory, Saint-Petersburg State University, Saint-Petersburg, RUSSIA}
\date{}

\begin{document}

\large\maketitle

\begin{abstract}
The main result of this paper is a functional limit theorem for  the
 sine-process.
 In particular, we study the limit distribution, in the space of
 trajectories, for the
 number of particles in a growing interval. The sine-process has the
 Kolmogorov property  and satisfies the Central Limit Theorem, but our
 functional limit theorem is very different from the Donsker Invariance
 Principle.
 We show that the time integral of our process can be approximated by
 the sum of a linear Gaussian process and independent Gaussian fluctuations whose
 covariance matrix is computed explicitly.
We interpret these results in terms of the Gaussian Free Field convergence for the random matrix models. 
 The proof relies on a general form of the multidimensional Central
 Limit Theorem  under the sine-process for linear statistics of two types:
those having growing variance
and those with bounded variance corresponding to observables
of Sobolev regularity
 $1/2$.
\end{abstract}

\tableofcontents

\section{Introduction}
\lbl{sec:intro}

\subsection{Formulation of the main result}
\lbl{sec:ii}

In this paper we
study the asymptotic
behaviour of trajectories
of determinantal random point processes;
for basic definitions and background concerning
determinantal point processes, see Section \ref{sec:detproc} below.
Mostly we deal with the sine-process given by the kernel
\bee\lbl{iKsine}
\KK_{\rm{sine}}(x,y)=\fr{\sin(x-y)}{\pi(x-y)}
\qmb{if}\qu
x\neq y
\qmb{and}\qu
\KK_{\rm{sine}}(x,x)\equiv \fr{1}{\pi}, \qu x,y\in\mR.
\eee
The sine-process is a strongly  chaotic stationary process: it satisfies  the Kolmogorov property \cite{Ly}, \cite{BQS}, \cite{OO},   having, therefore, Lebesgue spectrum and positive entropy,   and  enjoys an analogue of
the Gibbs property \cite{Buf14,Buf14a}, namely, the quasi-invariance under the group of diffeomorphisms with compact support.
At the same time, the sine-process is  rigid  in the sense of Ghosh and Peres \cite{G,GP}: the number of particles in a bounded interval is almost surely determined by the configuration in its exterior. The reason for the rigidity is the slow growth of the variance for the sine-process: for instance,  the number of particles $\#_{[0,N]}$
in the interval $[0,N]$ satisfies
\bee\lbl{iVar}
\Var\#_{[0,N]} =\fr{1}{\pi^2} \ln N + O(1),
\eee
see e.g. Exercise 4.2.40 from \cite{AGZ}. This slow growth of the variance
can be  seen from the form $\rho(\theta)=|\theta|$ of the spectral density  for the sine-process, cf. \cite{So00}, \cite{BDQ}.
  Costin and Lebowitz \cite{CL} showed
that the sine-process
satisfies the Central Limit Theorem: 
 the random variable
\bee\lbl{ixi}
\hat\xi^N:=\fr{\#_{[0,N]}-\MO \#_{[0,N]}}{\sqrt{\Var\#_{[0,N]}}}
\eee
converges in distribution to the normal law:
\bee\lbl{iCLT}
\DD(\hat\xi^N)\raw \NN(0,1) \ass N\ra\infty.
\eee
Here $\MO$, $\Var$ and $\DD$
stand for the expectation, variance and distribution under the sine-process.

The Central Limit Theorem was subsequently proven
for arbitrary determinantal processes governed by self-adjoint kernels and arbitrary additive statistics with growing variance 
(see \cite{So00,SoAB,So00b,So01, HKPV}),
in particular, for the Airy and Bessel processes (\cite{SoAB}),
for which the variance  of the number of particles has logarithmic growth and rigidity  holds \cite{Buf16}.

Many classical dynamical systems satisfying the Central Limit Theorem
also satisfy the Donsker Invariance Principle, which, informally speaking, states that trajectories of the system can be approximated by the Brownian motion, cf. \cite{Sinai89}.
The  main result of this paper is a functional limit theorem for the sine-process.
The limit dynamics is completely different  from  Brownian motion.
As far as we know,
this is the first example of such behaviour
in the theory of dynamical systems.
More specifically, we investigate asymptotic behaviour,  as $N\ra\infty$,
of the piecewise continuous random process
\bee\lbl{ixit}
\xi^N_t=\fr{\#_{[0,tN]}-\MO \#_{[0,tN]}}{\pi^{-1}\sqrt{\ln N}}, \qu 0\leq t\leq 1,
\eee
under the sine-process.
Trajectories of the process $\xi^N_t$
become extremely irregular when $N$ grows
(see \rprop\ref{lem:fd}),
so that the sequence of distributions of trajectories $\DD(\xi^N)$ does not have a limit
in any separable metric space.
That is why, instead of the process $\xi^N_t$ itself
we study its time integral
$\ili_0^t \xi^N_s\, ds$
in the space of continuous functions $C([0,1],\mR)$.
We fix $0<\tau\leq 1$ and set
\bee\lbl{ietaz}
\eta^N:=\fr{1}{\tau}\ili_0^{\tau} \xi^N_s \, ds
\qnd
z^N_t:=\pi^{-1}\sqrt{\ln N} \Big(\ili_0^t \xi^N_s \, ds - t\eta^N\Big),
\eee
so that
\bee\lbl{intxi}
\ili_0^t \xi^N_s \, ds = t\eta^N + \fr{z^N_t}{\pi^{-1}\sq{\ln N}}.
\eee
The parameter $\tau$ is fixed throughout the paper, and we skip it in the notation.
Recall that the Gaussian Free Field in the plain $\mR^2$ is a  generalized centred Gaussian process in $\mR^2$ given by the covariance function
\bee\lbl{GFFcov}
\GG(t,s)=-\fr1{2\pi}\ln|t-s|,
\quad t,s\in\mR.
\eee
Let
\bee\lbl{GFF}
W(t,s):=\pi^{-1}\Big(
	\ili_0^t\ili_0^s
	-\fr{s}{\tau}\ili_0^t\ili_0^\tau - \fr{t}{\tau}\ili_0^\tau\ili_0^s
	+\fr{ts}{\tau^2}\ili_0^\tau\ili_0^\tau \Big)\, \big(\GG(u,v)-\GG(u,0)-\GG(0,v)\big)\,dudv.
\eee
Computing \eqref{GFF} explicitly, we find the following expression. Set
$\ds{\tht(t):=\fr{t^2\ln|t|}{4\pi^2}}$ and $\tht(0):=0$. Then
$$
W(t,s)=w(t,s)+w(s,t),
$$
where
$$
w(t,s)
            :=\fr{1}{2}\tht(t-s)
            - \Big(1-\fr{s}{\tau}\Big)\tht(t)
            -\fr{s}{\tau}\tht(t-\tau)
            +\fr{t}{\tau}\Big(1-\fr{s}{\tau}\Big)\tht(\tau).
$$
%
%
%
Our first main result is
\btt\lbl{th:FCLTS}
For any $0<\tau\leq 1$,
under the sine-process
we have
the weak convergence of measures
\bee\lbl{reqconv}
\DD (\eta^N,z^N)\raw \DD(\eta,z) \ass
N\ra\infty \qmb{in} \qu \mR\times C([0,1],\mR),
\eee
where $\eta$ and $z$ are independent, $\eta\volna \NN(0,1/2)$
and $z_t$ is a centred continuous Gaussian random process with the covariances
\bee\lbl{covarz}
\MO z_t z_s
=W(t,s),
\quad
0\leq t,s\leq 1.
\eee
\ett
Proof of Theorem \ref{th:FCLTS} is given in \rsec \ref{sec:FCLTS}.
Informally, \rtheo \ref{th:FCLTS} states that,
up to terms of the size $o\big((\ln N)^{-1/2}\big)$,
the process
$\ili_0^t \xi^N_s\, ds$
can be decomposed to a linear random process
$\eta^N t$
and small Gaussian fluctuations
$z_t/\pi^{-1}\sqrt{\ln N}$,
see figure (\ref{pic:xi}).
Here $z_t$ is a continuous centred Gaussian process
whose covariances (\ref{covarz}), which we compute explicitly, 
are governed by the Gaussian Free Field in the plain; 
for explanation of this phenomenon see Section \ref{sec:GFF}.
About the linear process $\eta^N t$
we know that asymptotically it is governed by
the process $\eta t$,
where $\eta\volna\NN(0,1/2)$ is independent from $z$.
For the rate of convergence
of $\eta^N$ to $\eta$
we have
\begin{figure}[t]
\centering
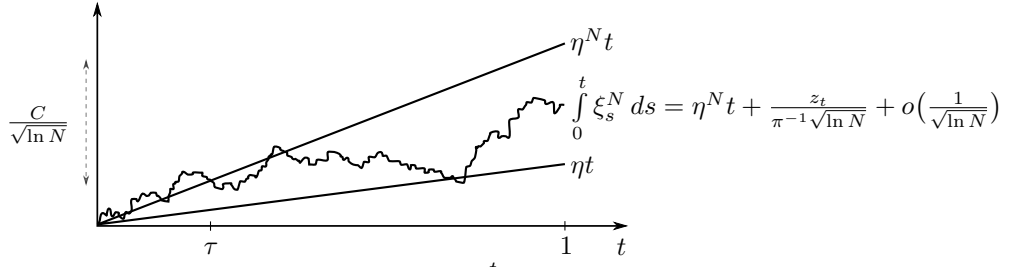
\caption{Up to terms of the size $o\big((\ln N)^{-1/2}\big)$, the integral  $\ili_0^t\xi^N_s \, ds$ decomposes to the sum of the linear in time process $\eta^N t$ and Gaussian fluctuations ${\ds\fr{z_t}{\pi^{-1}\sq{\ln N}}}$.
Deviation of the process $\eta^N t$ from the process $\eta t$ is of the size ${\ds\fr{C}{\sq{\ln N}}.}$ }
\lbl{pic:xi}
\end{figure}
\bpp\lbl{lem:i2}
Cumulants $(A^N_k)$ and $(A_k)$
of the random variables $\eta^N$ and $\eta$
satisfy $A^N_1=A_1=0$,
$|A^N_2-A_2|\leq C_2(\ln N)^{-1}$ and
\bee\lbl{cumtheo}
|A^N_k-A_k|\leq \fr{C_k}{(\ln N)^{k/2-1}} \qmb{for all} \qu k\geq 3,
\eee
with some constants $C_k$.
\epp
For a short reminding about cumulants see the beginning of \rsec\ref{sec:cumulants}.
Proof of \rprop \ref{lem:i2} is given in \rsec \ref{sec:FCLTS}.
Informally, \rprop \ref{lem:i2} states that
the deviation of the process
$\eta^Nt$, specifying the linear growth of the integral (\ref{intxi}),
from the process $\eta t$
is of the size $(\ln N)^{-1/2}$.
So that, it coincides
with the size
of the term $z_t/\pi^{-1}\sqrt{\ln N}$,
specifying the nonlinear fluctuations.


\rtheo\ref{th:FCLTS} has the following statistical interpretation.
\footnote{We are deeply grateful to Leonid Petrov for this remark.}
In order to predict behaviour of the process
$\ili_0^t \xi^N_s \, ds$ on the whole time interval $0\leq t\leq 1$
it suffices to know its realization at
arbitrarily small positive time $\tau$.
Indeed, then we determine $\eta^N$ by the formula (\ref{ietaz})
and approximate the integral $\ili_0^t \xi^N_s \, ds$
by the sum $\ds{\eta^N t+\fr{z_t}{\pi^{-1}\sqrt{\ln N}}}$,
where $z_t$ is the Gaussian process from \rtheo \ref{th:FCLTS}.

The main order asymptotic
$\DD(\eta^N)\overset{N\ra\infty}\raw\DD(\eta)$
from \rtheo\ref{th:FCLTS}
only uses the logarithmic growth of the variance
and holds for a general determinantal process
with logarithmically growing variance
(in particular, similar convergence takes place under the Airy and Bessel processes).
We show this in \rsec\ref{sec:FCLTw}.
To prove the asymptotic
$\DD(z^N)\raw\DD(z)$ however
we crucially use the form of the sine-kernel (\ref{iKsine}).
More specifically,
this asymptotic relies on a multidimensional Central Limit \rtheo \ref{th:CLTweak}
discussed in \rsec \ref{sec:iCLT}.
To establish the latter we analyse the corresponding cumulants
using a combinatorial identity (\ref{jSochG}) which is due to Soshnikov \cite{So00b}
and is specific for the sine-process.
While we expect the result to hold for the discrete sine-process,
additional arguments are needed.
It would be interesting to establish the convergence analogous to
$\DD(z^N)\raw\DD(z)$ for a general determinantal process with logarithmically growing variance
by using some different method, e.g. that of contour integrals developed in \cite{BF} (for discussion of this paper see Section~\ref{sec:GFF}).

The rest of \rsec\ref{sec:intro} is organized as follows.
In the next subsection we
describe motivation behind \rtheo\ref{th:FCLTS}.
In Section~\ref{sec:GFF} we explain why the Gaussian Free Field correlation function appears in \eqref{GFF}.
In \rsec \ref{sec:iergint} we state \rtheo\ref{th:ergintFCLTS}
which is our second main result.
There we show that a large class of observables ~---
ergodic integrals corresponding to a shift operator on the space of configurations,
has exactly the same asymptotic behaviour under the sine-process
as that described by \rtheo\ref{th:FCLTS}.
In \rsec\ref{sec:iCLT} we discuss the multidimensional Central Limit Theorem
mentioned above,
which is the main ingredient of the proofs of
Theorems \nolinebreak \ref{th:FCLTS} and \ref{th:ergintFCLTS}.
In \rsec\ref{sec:sketch} we outline the proof of \rtheo  \ref{th:FCLTS}.

\subsection{Finite dimensional distributions and motivation behind Theorem \ref{th:FCLTS}}

We first look at the finite-dimensional distributions of the process $\xi^N_t$.
Let $\eta_t$, $0\leq t \leq 1$, be a family of independent identically distributed Gaussian random variables satisfying $\eta_t\volna\NN(0,1/2)$.  
\bpp\lbl{lem:fd}
For any $0\leq t_1<\ldots<t_d\leq 1$, $d\geq 1$, we have
\bee\lbl{fdim_r}
\DD(\xi^N_{t_1},\ldots,\xi^N_{t_d}) \raw \DD(\eta_0-\eta_{t_1},\ldots,\eta_0-\eta_{t_d}) \ass N\ra\infty.
\eee
\epp
Proposition \ref{lem:fd} generalizes convergence (\ref{iCLT}) to many dimensions and is established in \rsec \ref{sec:FCLTS}.
Without a detailed proof a similar result was stated by Soshnikov, see \cite{So00}, p. 962 and \cite{SoAB}, p. 499.
Note that the covariance matrix $(b_{ij})_{1\leq i,j\leq N}$ of the limiting vector $(\eta_0-\eta_{t_1},\ldots,\eta_0-\eta_{t_d})$ has the form
\bee\lbl{iCov}
b_{ij}= 1/2+ \de_{ij}/2,
\eee
where $\de_{ij}$ is the Kronecker symbol; here we assume $t_1\neq 0$ since otherwise $b_{1j}=0$ due to the identity $\eta_0-\eta_{t_1}=0$.
The covariance matrix $(b_{ij})$ is independent
from the choice of times $t_1,\ldots, t_d$.
In particular, this means that if the limit as $N\ra\infty$
of the process $\xi^N_t$
exists in some sense,
then it can not be a continuous process,
so nothing as Brownian motion can appear.
Note that the proof of the convergence \eqref{fdim_r}
crucially uses the logarithmic growth of the variance (\ref{iVar}).
\brr\lbl{rem:fB} In the sense of finite-dimensional distributions,
asymptotic for large $N$ behaviour of the process $\xi^N_t$
is close to behaviour of the fractional Brownian motion $B^H_t$
with small parameter $H\ll 1$.
Indeed, in Lemma 4.1 of \cite{BMNZ} it is pointed out that
$\DD(B^H_{t_1},\ldots,B^H_{t_d})\raw \DD(\De\eta)$
as $H\ra 0$,
where  $\De\eta$ is the Gaussian vector from \rprop\ref{lem:fd}.
\err
In view of the convergence \eqref{fdim_r}, we expect that the limiting behaviour of the process $\xi^N_t$
is governed by the increments $\eta_0-\eta_t$ of the process $\eta_t$.
However, the process $\eta_t$ does not exist in a classical sense
(more precisely, it cannot be defined over a separable metric space).
That is why, in order to regularize the limiting dynamics,
instead of the process $\xi^N_t$
we study its time integral $\ili_0^t \xi^N_s \,ds$.
We expect that when $N\ra\infty$ the latter is governed by the difference
$
\ili_0^t \eta_0 - \eta_s \,ds= \eta_0 t - \ili_0^t \eta_s \,ds
$
where the integral $\ili_0^t \eta_s \,ds$ should be defined appropriately.
Here one can draw an analogy with the white noise,
which is not defined in the classical sense but its time integral gives the Brownian motion.
However, this heuristic idea leads us to the following rigorous result.
\bpp\lbl{ith:FCLTw}
For any function $\phi\in L^1[0,1]$ we have
\bee\lbl{FCLTw}
\DD\big(\ili_0^1 \phi(s)\xi^N_s\, ds\big) \raw \DD\big(\eta_0\ili_0^1 \phi(s)\,ds\big)  \ass N\ra\infty,
\eee
where $\eta_0\volna \NN(0,1/2)$.
\epp
\brr\lbl{rem:zerocovdelta}
Here and below by the normal law with zero expectation and variance
we understand the Dirac delta-measure at zero $\delta_0$.
In particular, if $\ili_0^1 \phi(s)\,ds=0$ then $\DD\big(\ili_0^1 \phi(s)\xi^N_s\, ds\big)\raw\delta_0$.
\err
Choosing
$\phi=\mI_{[0,t]}$
we find the leading term of the asymptotic
for the process $\ili_0^t \xi^N_s \,ds$, claimed in \rtheo \ref{th:FCLTS}:
\bee\lbl{mainorderintro}
\DD\Big(\ili_0^t \xi^N_s \,ds\Big) \raw \DD(\eta t) \ass N\ra\infty,
\eee
where we set $\eta:=\eta_0$.
Thus, we do not observe the integral $\ili_0^t \eta_s \,ds$.
The reason is that the process $\eta_t$ is completely uncorrelated in time
and has a bounded variance
(in difference with the white noise whose variance is the delta-function).
So that, $\eta_t$ oscillates fast with not very large amplitude
and averages out under the integration over the interval $[0,t]$.
Note that convergence (\ref{FCLTw}) takes place
even for a very rough observable $\phi$:
only integrability of $\phi$ is assumed.

Proposition \ref{ith:FCLTw} is a particular case of Proposition \ref{th:FCLTw},
in which we establish a stronger result
for an important class of determinantal point processes
including those with logarithmically growing variance,
as the sine, Airy and Bessel processes;
see Section \nolinebreak \ref{sec:FCLTw}.

Proposition \ref{ith:FCLTw} gives some information about
the asymptotic behaviour
of the process $\xi_t^N$.
But we lose a lot:
we do not observe any influence of the process $\eta_t$
which we find at the level of finite dimensional distributions.
Our next goal is to <<catch>> the process $\eta_t$.
The informal identity $\ili_0^t \eta_s \,ds=0$
resembles the law of large numbers.
To observe
the influence of $\eta_t$
we try to look at
the Central Limit Theorem scaling.
Since we expect that, informally,
$$
\ili_0^t \xi_s^N \, ds -\eta_0 t \ra -\ili_0^t \eta_s \, ds
\ass N\ra\infty,
$$
we need to find a sequence $\al_N\ra\infty$ as $N\ra\infty$,
such that the random process
\bee\lbl{iz}
z^N_t=\al_N\Big(\ili_0^t\xi^N_s \,ds -\eta_0 t\Big)
\eee
converges to a non-trivial limit.
However, joint distribution of the process $\xi^N_t$ and the random variable $\eta_0$
is undefined.
To overcome this difficulty we note that,
due to (\ref{mainorderintro}),
$\DD(\eta^N)\os{N\ra\infty}\raw \DD(\eta_0)$
where
$\eta^N$
is defined in (\ref{ietaz}),
and
replace  in  the definition (\ref{iz}) of the process $z^N_t$ the random variable $\eta_0$ by $\eta^N$ .
Then, setting
 $\al_N = \pi^{-1}\sqrt{\ln N}$
we arrive at
\rtheo \ref{th:FCLTS}.
\brr
It could seem that influence of
the process $\eta_t$
could be discovered by
consideration of some {\it nonlinear}
functional of the process $\xi^N_t$
such as,
for example, the integral
$
\ili_0^1 \phi(t)(\xi^N_t)^m\, dt
$
for integer $m\geq 2$,
where $\phi\in L^1[0,1]$.
However, this is not the case.
Indeed, we expect that
$(\xi_t^N)^m \approx (\eta_0-\eta_t)^m=\sli_{k=0}^m (-1)^kC_m^k\eta_t^k\eta_0^{m-k}$,
if $N$ is large.
Since terms $\eta_t^k$ and $\eta_s^k$  are independent for $t\neq s$,
the situation here is similar to that of Proposition~\ref{ith:FCLTw}:
the integral $\ili_0^1 \phi(t)\eta_t^k \, dt$
averages the terms $\eta_t^k$,
so feels only their means $\MO\eta_t^k$.
More precisely, one can prove that
\bee\lbl{FCLTwnl}
\DD\big(\ili_0^1 \phi(t)(\xi^N_t)^m\, dt\big)
\os{N\ra\infty}\raw
\DD\big(\ili_0^1 \phi(s)\, ds \sli_{k=0}^m (-1)^kC_m^k\eta_0^{m-k} \MO\eta_t^k\big) .
\eee
Comparing with the right-hand side of (\ref{FCLTw}),
the r.h.s. of (\ref{FCLTwnl})
depends on the moments
$\MO\eta_t^k$,
so that now we feel the ``noise'' $\eta_t$ but in a trivial way.
Indeed, all the randomness is still due to $\eta_0$,
although modified by the moments of $\eta_t$.
\err

\subsection{Connection with the Gaussian Free Field}
\lbl{sec:GFF}
In this section we explain in heuristic way the appearance of the Gaussian Free Field in the formula \eqref{GFF}.
Let us recall that the Gaussian Free Field in the plane $\mR^2$ is a  generalized centred Gaussian process in $\mR^2$ whose covariance is given by the Green function $\GG$ of the Laplace operator in the plane, $\Delta_t\GG(t,s)=\delta(t-s)$, where $\delta$ is the Dirac delta-function and $t,s\in\mR$. Explicitly, the function $\GG$ takes the form \eqref{GFFcov}.
The Gaussian Free Field in the upper half-plain $\mH:=\{z\in\mC:\,\Im z>0\}$ is defined similarly: its correlation function $\wid\GG$ is given by the Green function of the Laplace operator in $\mH$ with the Dirichlet boundary conditions. Explicitly,
$$
\wid\GG(z,w)=-\fr1{2\pi}\ln\Big|\fr{z-w}{z-\bar w}\Big|, 
\quad z,w\in\mC.
$$ 
The Gaussian Free Field is known to describe the asymptotic behaviour of the so-called height function in many models that have the random matrix-type behaviour. 
In particular, this is the case for a class of random surfaces \cite{dT,K,Du,BF,Ku,P} and random matrices \cite{Bor,BG}.
Below we explain this in more details on the examples provided by \cite{BF,Bor}.

Let $A^N$ be an $N\times N$ random matrix from the Gaussian Unitary Ensemble.
That is, $A^N=(a^N_{ij})_{1\leq i,j\leq N}$ is Hermitian, 
where $a_{ll}^N\volna\NN(0,1)$, $\Re a_{ij}^N,\Im a_{ij}^N\volna \NN(0,1/2)$ for $i\neq j$
and $a_{ll}, \Re a_{ij}, \Im a_{ij}$ are independent for all $l$ and $i>j$. 
Introduce the height function $h^N_t$ by the formula
\bee\non
h^N_t:=\{\mbox{the number of eigenvalues of $A^N$ that are} \geq \sqrt{N}t\}.
\eee
In \cite {Bor} it is shown that the limiting behaviour of the random process $h^N_t$ is governed by the Gaussian Free Field in the upper half-plane. In particular, it is proven that, roughly speaking, 
\bee\lbl{cov_h_r}
\Cov(h^N_t,h^N_s)\to const\, \wid\GG(y_t,y_s) \ass N\to\infty\quad\mbox{for}\quad t\neq s,
\eee 
where $y_t,y_s\in\mH$ are defined appropriately.

In \cite{BF} a family of stochastic growth models in $2+1$ dimensions is considered, that belongs to the anisotropic KPZ class. For an appropriately defined height function $h^N_t$, where $t$ is a parameter from $\mR^3$, the Central Limit Theorem is proven,
\bee\lbl{h_CLT}
\DD\Big(\fr{h_t^N-\MO h_t^N}{\sqrt{\ln N}}\Big)\raw \DD(\eta_t) 
\ass N\to \infty
\qmb{for any }t,
\eee 
where $\eta_t\volna\NN(0,\sigma^2)$, for some $\sigma^2>0$ independent from $t$. In particular, we have $\Var h_t^N\volna \sigma^2\ln N$. 
On the other hand, as in \eqref{cov_h_r}, it is shown that for $t\neq s$ the limiting covariance $\Cov(h^N_t,h^N_s)$
of the height function is governed by the Gaussian Free Field in the upper half-plain, without the normalization factor $(\ln N)^{-1}$ (actually, in \cite{BF} a similar result is also proven for higher order moments). 

Since the asymptotic behaviour of the spectrum of a large class of Wigner matrices 
(in particular, of the Gaussian Unitary Ensemble) is governed by the sine-process,
the appearance of the Gaussian Free Field in \rtheo\ref{th:FCLTS} is not surprising.
Indeed, the random process $\#_{[0,Nt]}$ we are interested in can be represented as
\bee
\lbl{prochigh_r}
\#_{[0,Nt]}=h_0^N-h_t^N,
\eee
where we set $h_s^N:=\#_{[sN,\infty)}$. 
However, the height function $h_s^N$ in this case is defined only informally, since we have $h_s^N=\infty$ almost surely.
Then, heuristically,
the process $\xi_t^N$ defined in \eqref{ixit}  takes the form
$$
\xi^N_t=\fr{h_0^N-h_t^N-\MO(h_0^N-h_t^N)}{\pi^{-1}\sqrt{\ln N}}.
$$
Drawing an analogy with \eqref{cov_h_r} and \eqref{h_CLT}, one could expect the convergence 
$$
\DD(\xi^N_t)\raw\DD(\eta_0-\eta_t) 
\qmb{for any }t,
$$
where $\eta_0,\eta_t$ are independent identically distributed centred Gaussian random variables.
Due to Proposition~\ref{lem:fd}, we see that, indeed, this convergence takes place. 
Thus, representation \eqref{prochigh_r}, also being informal,
gives a correct intuition for the asymptotic behaviour of the process $\#_{[0,Nt]}$.
Let us now informally express the covariance $\Cov(z_t^N,z_s^N)$ of the process $z_t^N$, 
defined in \eqref{ietaz}, through the covariance of the height function $h_t^N$.
We have
\begin{align}\lbl{covzz_r}
&\Cov(z_t^N,z_s^N)
=
\fr{\ln N}{\pi^2}\Cov\Big(\ili_0^t\xi_u^N\,du-\frac{t}{\tau}\ili_0^\tau\xi^N_u\,du,
\ili_0^s\xi_v^N\,dv-\frac{s}{\tau}\ili_0^\tau\xi^N_v\,dv
\Big)\\\non
&=\Big(
\ili_0^t\ili_0^s\, dudv
-\fr{t}{\tau}\ili_0^\tau\ili_0^s\, dudv
-\fr{s}{\tau}\ili_0^t\ili_0^\tau \, dudv
+\fr{ts}{\tau^2}\ili_0^\tau\ili_0^\tau \, dudv
\Big)\Cov(\#_{[0,uN]},\#_{[0,vN]}).
\end{align}
Due to \eqref{prochigh_r},
\bee\lbl{cornp_r}
\Cov(\#_{[0,uN]},\#_{[0,vN]})=
\Var h_0^N+\Cov(h_u^N,h_v^N)-\Cov(h_0^N,h_u^N)-\Cov(h_0^N,h_v^N).
\eee
Then
\begin{align}\non
\Cov(z_t^N,z_s^N)
=\Big(&\ili_0^t\ili_0^s\, dudv
-\fr{t}{\tau}\ili_0^\tau\ili_0^s\, dudv
-\fr{s}{\tau}\ili_0^t\ili_0^\tau \, dudv
+\fr{ts}{\tau^2}\ili_0^\tau\ili_0^\tau \, dudv
\Big) \\\non
&\big(\Cov(h_u^N,h_v^N)-\Cov(h_0^N,h_u^N)-\Cov(h_0^N,h_v^N)\big),
\end{align}
since the term $\Var h_0^N$ is cancelled by the integration.
In view of \eqref{cov_h_r}, the formula above perfectly agrees with  \eqref{GFF},\eqref{covarz}. The only difference is that in our case the covariance $\Cov(h_t^N,h_s^N)$
is expected to be governed by the Gaussian Free Field in the plane (as., e.g., in \cite{dT}) rather than in the half-plain as in \eqref{cov_h_r}. Namely, we expect the convergence
\bee\lbl{h_th_stoGFF}
\Cov(h_t^N,h_s^N)\to\pi^{-1}\GG(t,s) \ass N\to\infty, \qmb{for }t\neq s.
\eee
Note that the only term in \eqref{cornp_r} which is expected to grow with $N$ 
is $\Var h_0^N$. 
In particular, the reason for which the covariance 
$\Cov(z_t^N,z_s^N)$ does not grow is that this term  is cancelled by the integration in \eqref{covzz_r}.
This gives an indication to the fact that this is the fixed end $0$ of the interval $[0,Nt]$ who is responsible for the growth of the covariance \eqref{cornp_r}.

To make sure that the heuristic argument above is correct, we next consider the random process 
$\#_{[aNt, bNt]}$, where for definiteness we assume $t>0$ and $a<b$, $a,b\in\mR\setminus\{0\}$.
On an intuitive level, we have
$$
\#_{[aNt, bNt]}=h^N_{at}-h^N_{bt}.
$$ 
Then
$$
\Cov(\#_{[aNt, bNt]},\#_{[aNs, bNs]})=
\Cov(h_{at}^N,h_{as}^N)+\Cov(h_{bt}^N,h_{bs}^N)-
\Cov(h_{at}^N,h_{bs}^N)-\Cov(h_{bt}^N,h_{as}^N).
$$
Consequently, if the ends of the intervals $[at,bt]$ and $[as,bs]$ do not coincide,  we expect the covariance 
$\Cov(\#_{[aNt, bNt]},\#_{[aNs, bNs]})$
to be bounded in $N$ and to be governed by the Gaussian Free Field in the plain 
accordingly to the convergence 
\eqref{h_th_stoGFF}. 
This turns out to be true, namely in Section~\ref{sec:ergintFCLTS} we prove
\bpp\lbl{lem:covnotfixed}
Assume that $t,s>0$ satisfy $t\neq s$, $at\neq bs$ and $bt\neq as$, so that the ends of the intervals 
$[at,bt]$ and $[as,bs]$ do not coincide. Then, as $N\to\infty$, we have
$$
\Cov(\#_{[aNt, bNt]},\#_{[aNs, bNs]})\to
\pi^{-1}\big(\GG(at,as) + \GG(bt,bs)-\GG(at,bs)-\GG(bt,as)\big).
$$
\epp
In particular, one could study directly the integral 
$I_t^N:=\ili_0^t \#_{[aNs, bNs]} - \MO\#_{[aNs, bNs]}\,ds$, 
without the normalization $(\ln N)^{-1/2}$ and  
without extracting the linear term $\eta^N t$ 
as in \eqref{intxi}, which appear because of the fixed left end of the interval $[0,t]$. 
Indeed, due to Proposition~\ref{lem:covnotfixed}, 
the covariances $\Cov(I_t^N,I_s^N)$ converge as $N\to\infty$ for any $t,s>0$
and are governed by the Gaussian Free Field.  

\subsection{Functional limit theorem for ergodic integrals}
\lbl{sec:iergint}

In this section we explain that ergodic integrals corresponding to a shift operator
acting on the space of configurations
possess the same asymptotic behaviour as the number of particles
$\#_{[0,N]}$.
Denote by $\Conf(\mR)$ the space of locally finite configurations on $\mR$,
$$
\Conf(\mR)=\big\{\XX\subset\mR\big| \XX\mbox{ does not have limit points in }\mR\big\}.
$$
Let $ T^u$, $u\geq 0$, be a shift operator acting on $\Conf(\mR)$ as
$$
T^u:\Conf(\mR)\mapsto \Conf(\mR),\qu
T^u(\XX)=\XX-u.
$$
Consider the dynamical system
\bee\lbl{dynsyst}
\big(\Conf(\mR), (T^u)_{u\geq 0},\PR\big),
\eee
where $\PR$ is the probability measure on $\Conf(\mR)$,
given by the sine-process.
Take a bounded measurable function $\ph:\mR\mapsto\mR$ with compact support.
The {\it linear statistics} $\SSS_\ph$
corresponding to the function $\ph$
is introduced by the formula
\bee\lbl{lstatdef}
\SSS_\ph:\, \Conf(\mR)\mapsto\mR,
\qu
\SSS_\ph(\XX):=\sli_{x\in\XX}\ph(x).
\eee
In particular, if
$\ph=\mI_{[a,b]},$
we have
$\SSS_\ph=\#_{[a,b]}$.
Assume that the function $\ph$
satisfies the normalization requirement
\bee\lbl{avnonzero}
\ili_{-\infty}^\infty \ph(u)\,du =1.
\eee
Consider the ergodic integral
$$
\ili_0^{tN} \SSS_\ph\circ T^u \,du,
\qu
0\leq t\leq 1,
$$
where
\bee\lbl{ergint1}
\SSS_\ph\circ T^u (\XX)=\sli_{x\in T^u(\XX)}\ph(x)
=\sli_{x\in \XX}\ph(x-u).
\eee
Let
$$
\ph_t^N:=\ili_0^{tN} \ph(\cdot -u) \, du.
$$
Then, exchanging the integral with the sum, we see that the ergodic integral coincides with the linear statistics
$\SSS_{\ph_t^N}$,
\bee\lbl{ergint}
\ili_0^{tN} \SSS_\ph\circ T^u \,du = \SSS_{\ph_t^N}.
\eee
Consider the random process
$$
\xi^N_{\ph,t}:=
\fr{\SSS_{\ph_{t}^N} -\MO\SSS_{\ph_{t}^N}}
{\pi^{-1}\sqrt{\ln N}}.
$$
In the next theorem, which is our second main result,
we show that
under the sine-process
the process $\xi^N_{\ph,t}$ possesses exactly the same asymptotic
behaviour as the process $\xi^N_t$ given by (\ref{ixit}).
Fix $0<\tau\leq 1$ and set
$$\eta^N=\fr{1}{\tau}\ili_0^{\tau} \xi_{\ph,s}^N\, ds.$$
Choose the random process $z_t^N$ in such a way that
\bee\non
\ili_0^t  \xi_{\ph,s}^N \, ds = t\eta^N + \fr{z^N_t}{\pi^{-1}\sq{\ln N}}.
\eee
\btt\lbl{th:ergintFCLTS} Under the sine-process we have
\begin{enumerate}
\item
For any $t>0$,
$$
\Var\SSS_{\ph_{t}^N} = \pi^{-2}\ln N + O(\sqrt{\ln N})
\ass N\ra\infty.
$$

\item
For any $0< t_1<\ldots<t_d\leq 1$, $d\geq 1$,
distribution of the random vector
$$\xi^N_{\ph}:=(\xi^N_{\ph,t_1}, \ldots, \xi^N_{\ph,t_d})$$
satisfies
$\DD(\xi^N_{\ph})\raw\DD(\De\eta)$ as $N\ra\infty$,
where $\De\eta$ the Gaussian random vector
from \rprop\ref{lem:fd}.

\item
The distribution $\DD(\eta^N,z^N)$ satisfies
$
\DD(\eta^N,z^N)\raw\DD(\eta,z)
$
as $N\ra\infty$
in $\mR\times C([0,1],\mR)$,
where the random variable $\eta$ and the random process $z_t$
are as in \rtheo\ref{th:FCLTS}.

\item Cumulants $(A_k^N)$ and $(A_k)$ of the random variables
$\eta^N$ and $\eta$ satisfy $A_1^N=A_1=0$,
$|A^N_2-A_2|\leq C_2(\ln N)^{-1/2}$,
and (\ref{cumtheo}) for $k\geq 3$ and
some constants $C_k$.
\end{enumerate}
\ett
Theorem \ref{th:ergintFCLTS} is proven in \rsec\ref{sec:ergintFCLTS}.
To see the connection between the processes $\xi^N_t$ and $\xi^N_{\ph,t}$
observe that the function $\ph^N_t$ has the form
as shown on figure \nolinebreak \ref{pic:phi}:
it has a flat part of the length $\volna N$
where $\ph^N_t(x)=\ili_{-\infty}^\infty \ph(x)\,dx=1$,
and <<tails>> with the length of order one.
So that, $\ph^N_t$  <<almost coincides>>
with a shifted indicator function $\mI_{[0,N]}$,
if $N$ is large.
But the linear statistics $\SSS_{\mI_{[0,N]}}$
is equal exactly to the number of particles $\#_{[0,N]}$.
\begin{figure}[t]
\centering
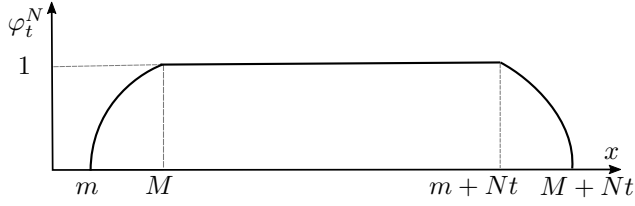
\caption{Function $\ph_t^N.$ Here $m:=\inf\supp\ph$ and $M:=\sup\supp\ph$. }
\lbl{pic:phi}
\end{figure}

\subsection{Central Limit Theorem for linear statistics}
\lbl{sec:iCLT}

Proofs of Theorems \ref{th:FCLTS} and \ref{th:ergintFCLTS}
follow the same pattern
and rely on the multidimensional Central Limit Theorem
\nolinebreak \ref{th:CLTj},
which we state below in a simpler form.
Recall that the linear statistic $\SSS_\ph$ of a function $\ph$
is defined in (\ref{lstatdef}).
\btt\lbl{th:CLTweak}
Let $f_1,\ldots,f_p,g_1,\ldots,g_q:\mR\mapsto\mR$, $p,q\geq 0$,
be measurable bounded functions
with compact supports.
Set $f_i^N:=f_i(\cdot/N)$, $g_j^N:=g_j(\cdot/N)$ and
consider the corresponding linear statistics
$$
\SSS_{f_1^N},\ldots, \SSS_{f_p^N}, \SSS_{g_1^N},\ldots, \SSS_{g_q^N}
$$
as random variables under the sine-process.
Assume that
\begin{enumerate}
\item
There exists a sequence
$V_N\ra\infty$ as $N\ra\infty$
and numbers $b_{ij}^f$ satisfying $b_{ii}^f>0$,
such that for any $i,j$
\bee\lbl{iCovN}
\fr{\Cov(\SSS_{f_i^N},\SSS_{f_j^N})}{V_N} \ra b^f_{ij} \ass N\ra\infty.
\eee

\item
The functions $g_i$
belong to the Sobolev space $H^{1/2}(\mR)$.
\end{enumerate}
Let $(\xi^N_f,\xi^N_g)$ be the random vector with components
$$
\xi^N_{f_i}:=\fr{\SSS_{f_i^N} -\MO\SSS_{f_i^N}}{\sqrt{V_N}}
\qmb{and}\qu
\xi^N_{g_j}:=\SSS_{g_j^N} -\MO\SSS_{g_j^N}.
$$
Then we have the weak convergence
$\DD(\xi^N_f,\xi^N_g) \raw \DD(\xi_f,\xi_g)$,
where $(\xi_f,\xi_g)$ is a centred Gaussian random vector
with the covariance matrix
$\begin{pmatrix}
(b^f_{ij})&0 \\
0&(b^g_{kl})
\end{pmatrix}$
and $b^g_{kl} = \lan g_k,g_l\ran_{1/2} $,
where the pairing $\lan\cdot,\cdot\ran_{1/2}$
is given by (\ref{pairing}).
\ett
Note that under the assumption $g_i\in H^{1/2}$
the variances $\Var\SSS_{g^N_i}$ do not grow at all,
so that assumption (\ref{iCovN})
can not be satisfied for the functions $g_i$.
Conversely, the inclusion $f_i\in H^{1/2}$
can not take place once (\ref{iCovN}) holds.

The difference between Theorems \ref{th:CLTweak} and \ref{th:CLTj}
is that in the latter we admit more general
dependence of the functions $f^N_i, g^N_j$
on $N$ than in \rtheo \ref{th:CLTweak}.
This is needed for the proof of \rtheo\ref{th:ergintFCLTS}.

The marginal convergence
$\DD(\xi^N_f)\overset{N\ra\infty}\raw\DD(\xi_f)$
does not use the special structure of the sine-kernel and takes place
under a large class of determinantal point processes,
once (\ref{iCovN}) holds.
We prove this in \rtheo\ref{th:CLT}
and use in \rsec\ref{sec:FCLTw},
where we establish the main order asymptotic from \rtheo\ref{th:FCLTS}
for a general determinantal process with logarithmically growing variance.
To establish the convergence
$\DD(\xi^N_g)\overset{N\ra\infty}\raw\DD(\xi_g)$,
however,
we crucially use the form of the sine-kernel.
Indeed, proof of our Central Limit Theorem is based on analysis of
cumulants $(A_k^N)_{k\in\mZ_+^{p+q}}$ of the random vector $(\xi^N_f,\xi^N_g)$.
In particular, we
show that $A_k^N\overset{N\ra\infty}\ra 0$
once $|k|>2$.
For the cumulants corresponding to the component $\xi^N_f$
the latter convergence follows from general estimates obtained in \rsec\ref{sec:cumulants}
and decay of the normalization factor $V_N^{-1}$.
For the component $\xi^N_g$ such normalization is lacking
and the analysis is more delicate.
We rely  on the combinatorial identity (\ref{jSochG}) obtained by Soshnikov in \cite{So00b},
while application of the latter requires the relation (\ref{Soshid})
which is specific for the sine-process.

The main novelty of  \rtheo \ref{th:CLTweak}
is that we study asymptotic behaviour of the \emph{joint}
linear statistics $(\SSS_{f_i^N},\SSS_{g_j^N})$,
so that we work simultaneously on two different scales,
corresponding to the growing and bounded variance.
Indeed, the marginal convergence $\DD(\xi^N_{f_i})\ra \DD(\xi_{f_i})$
in the generality of \rtheo\ref{th:CLT}
generalizes convergences obtained by
Costin and Lebowitz \cite{CL} and Soshnikov \cite{So00,SoAB,So01},
see \rsec\ref{sec:CLT} for the discussion.
The convergence
$\DD(\xi^N_{g_i})\ra \DD(\xi_{g_i})$
was proven by Spohn \cite{Sp} and Soshnikov
\cite{So00b,So01}.
For further developments see also works \cite{JL,L15,L15a,BD16,BD17},
where certain Central Limit Theorems were established  for linear
  statistics with bounded variance, related to the marginal convergence
$\DD(\xi^N_{g_i})\ra \DD(\xi_{g_i})$.
More precisely, in \cite{JL,L15} and \cite{BD16} the Central Limit Theorems were proven for linear statistics of various orthogonal polynomial ensembles on mesoscopic scales.
In \cite{L15a} and \cite{BD17} those were obtained for linear statistics of certain biorthogonal ensembles.

\subsection{Outline of the proofs of Theorems \ref{th:FCLTS} and \ref{th:ergintFCLTS} }
\lbl{sec:sketch}

First we discuss \rtheo \ref{th:FCLTS}. We note that, due to (\ref{ixit}),
$$
\eta^N=\fr{\tau^{-1}\ili_0^{\tau} \big(\#_{[0,sN]} -\MO \#_{[0,sN]}\big) \,ds}{\pi^{-1}\sq{\ln N}}
=\fr{\SSS_{f^N}-\MO \SSS_{f^N}}{\pi^{-1}\sqrt{ \ln N}},
$$
where $\SSS_{f^N}$ is the linear statistics corresponding to the function
$f^N(x)=f(x/N)$ with
$f(x)=\fr{1}{\tau}\ili_0^{\tau} \mI_{[0,s]}(x)\,ds$.
Similarly,
$$
z_t^N
=\ili_0^t \big(\#_{[0,sN]} -\MO \#_{[0,sN]}\big) \,ds
- \fr{t}{\tau} \ili_0^{\tau} \big(\#_{[0,sN]} - \MO \#_{[0,sN]}\big) \, ds
=\SSS_{g^N_t}-\MO \SSS_{g^N_t},
$$
where $g^N_t(x)=g_t(x/N)$ and
the functions $g_t$, $0\leq t\leq 1$,
are given by
\bee\lbl{igtN}
g_t(x):=\ili_0^t \mI_{[0,s]}(x)\, ds - \fr{t}{\tau}\ili_0^{\tau} \mI_{[0,s]}(x)\, ds.
\eee
It is easy to see that the functions $f$ and $g_t$ have compact support, are piecewise linear,
and the functions $g_t$ are continuous (see (\ref{gtx})-(\ref{gtx1}) for the explicit form of $g_t$).
In particular, $g_t\in H^1(\mR)$ for all $0\leq t\leq 1$.

Next we show that
$\ds{\Var\SSS_{f^N}\volna \fr{1}{2\pi^{2}} \ln N}$
and that the pairing
$\lan g_t,g_s\ran_{1/2}$
equals to the right-hand side of (\ref{covarz}).
Thus,
for any
$0\leq t_1<\ldots< t_d\leq 1$
the functions
$(f, g_{t_1},\ldots, g_{t_d})$
satisfy assumptions of \rtheo\ref{th:CLTweak},
with $V_N=\pi^{-2} \ln N$, $b^f_{11}=1/2$
and
$
b^g_{ij}=\mbox{ r.h.s. of (\ref{covarz})}.
$
The latter implies the convergence
\bee\lbl{iconvfinite}
\DD(\eta^N,z^N_{t_1},\ldots,z^N_{t_d})
\raw
\DD(\eta,z_{t_1},\ldots,z_{t_d})
\ass N\ra\infty,
\eee
where the random variable $\eta$ and the random process $z_t$
are as in \rtheo \ref{th:FCLTS}.
Then, using a compactness argument in a standard way,
we show that convergence (\ref{iconvfinite})
implies assertion of the theorem.

Proof of \rtheo \ref{th:ergintFCLTS} uses similar argument. 
Its main difference from the proof of \rtheo\ref{th:FCLTS} is that the functions $f^N$ and $g_t^N$ depend on $N$ in a more complicated way. 
That is why instead of \rtheo\ref{th:CLTweak} 
we use more general 
\rtheo\ref{th:CLTj}.

\subsection{Organization of the paper}

In \rsec\ref{sec:preliminaries} we first introduce notation
which will be used throughout the paper.
Then we recall some basic definitions concerning determinantal point processes
and establish some simple facts needed in the sequel.
In \rsec \ref{sec:cumulants}
we compute and estimate
cumulants of linear statistics first under a general determinantal process
and then specify our attention on the sine-process.
Results obtained there
are used in \rsec \ref{sec:CLTgen},
where
we establish the Central Limit Theorems
\nolinebreak
\ref{th:CLT} and  \ref{th:CLTj},
which are discussed
\rsec \ref{sec:iCLT}.
\rsec \ref{sec:FCLT} is devoted to the proofs of our main results:
Theorems \nolinebreak \ref{th:FCLTS}, \ref{th:ergintFCLTS} and 
Propositions \ref{lem:i2}, \ref{lem:fd}, \ref{lem:covnotfixed}.
In \rsec \ref{sec:FCLTw} we prove an analogue of \rprop \ref{ith:FCLTw} for an important class of determinantal processes,
including those with logarithmically growing variance (in particular, the Airy and Bessel processes).

\section{Preliminaries}
\lbl{sec:preliminaries}

\subsection{Notation}
\lbl{sec:Notation}
\begin{enumerate}
\item By $C,C_1,\ldots$ we denote various positive constants. By $C(a), \ldots$ we denote constants depending on a parameter $a$.
Unless otherwise stated, the constants never depend on $N$.

\item For $d\geq 1$ we set $\Zp:=\{\mZ^d\ni k=(k_1,\ldots,k_d)\neq 0:\, k_j\geq 0\; \forall 1\leq j\leq d\}$.

\item\lbl{not:!} For $k\in\mZ^d_+$ and $z\in\mC^d$ we denote
$|k|:=k_1+\cdots+k_d$, $k!:=k_1!\cdots k_d!$ and $z^k:=z_1^{k_1}\cdots z_d^{k_d}$.

\item \lbl{not:F} Our convention for the Fourier transform is as follows:
$
\hat h(t)=\FF(h)=\ili_{-\infty}^\infty h(x)e^{-itx}\, dx.
$
For the inverse Fourier transform we write
$
\FF^{-1}(\hat h)(x)=(2\pi)^{-1}\ili_{-\infty}^\infty \hat h(t)e^{itx}\, dt.
$
\item We denote by $\|\cdot\|$ the usual operator norm, by $\|\cdot\|_{HS}$
~--- the Hilbert-Schmidt norm
and by $\|\cdot\|_\infty$ and $\|\cdot\|_{L^m}$, $m\geq 1$,
~--- the Lebesgue $L^\infty$ and  $L^m$-norms.
By $H^n(\mR)$, $n>0$, we denote the Sobolev space of order  $n$
and for functions $f,g\in H^n(\mR)$ we set
\bee\lbl{pairing}
\|f\|^2_{n} := \fr{1}{4\pi^2}\ili_{-\infty}^\infty |u|^{2n} |\hat f(u)|^2 \,du,
\qu
\lan f,g\ran_{n} := \fr{1}{4\pi^2} \ili_{-\infty}^\infty |u|^{2n}
\hat f(u) \ov{\hat g(u)}\,du,
\eee
and $\|f\|_{H^n}^2:=\|f\|^2_{L^2} + 2\pi\|f\|^2_n$.

\item By $\CoRr$ we denote the space of locally finite configurations of particles in
$\mR^m$, $m\geq 1$,
\bee\lbl{defConf}
\CoRr:=\big\{\XX\subset\mR^{m}\big| \XX\mbox{ does not have limit points in }\mR^m\big\}.
\eee

\item
Let $\XX\in\CoRr$.
By $\#_B(\XX):=\#\{B\cap\XX\}$ we denote the number of particles
from the configuration $\XX$
intersected with the set $B$.

\item For a bounded compactly supported function
$h:\mR^m\mapsto\mR$, by $\SSS_h$
we denote the corresponding linear statistics,
$$
\SSS_h:\CoRr\mapsto \mR,
\qu
\SSS_h(\XX)= \sli_{x\in\XX} h(x).
$$
\item By $\mI_B$ we denote the indicator function of a set $B\in\mR^m$.
\end{enumerate}

\subsection{Determinantal point processes}
\lbl{sec:detproc}

In this section we recall some basic definitions and facts
concerning determinantal processes.
Determinantal (or fermion) random point processes
form a special class of random point processes,
which was introduced by Macchi in seventies (see \cite{Ma75,Ma77, DVJ}).
They play an important role in the random matrix theory,
statistical and quantum mechanics, probability, representation and number theory.
For detailed background see \cite{So00,ST,STa},
see also Chapter 4.2 in \cite{AGZ}.

Consider on the space of locally finite configurations  $\CoRr$,
defined in (\ref{defConf}),  a $\sigma$-algebra $\FF$
generated by cylinder sets
$$
\CC^n_B=\{\XX\in\CoRr:\, \#_B(\XX)=n\},
$$
where $n$ and $B$ run over natural numbers and bounded Borel subsets of $\mR^m$ correspondingly.
The triple $(\CoRr,\FF,\PR)$,
where $\PR$ is a probability measure on $(\CoRr,\FF)$,
is called a {\it  random point process}.

Assume that there exists a family of locally integrable nonnegative functions
$
\rho_n:(\mR^m)^n\mapsto\mR, \; n\geq 1,
$
such that for any $n\geq 1$ and any mutually disjoint Borel subsets $B_1,\ldots,B_n$ of $\mR^m$
we have
$$
\MO \#_{B_1}\cdots \#_{B_n}
= \int_{B_1\times\ldots\times B_n}\rho_n(x_1,\ldots,x_n)\, dx_1 \cdots dx_n.
$$
The functions $\rho_n$ are called {\it correlation functions}.
Under natural assumptions the family  $(\rho_n)_{n\geq 1}$  determines the probability $\PR$ uniquely, see e.g. \cite{So00}.

Consider a non-negative integral operator
$K: L^2(\mR^m,dx)\mapsto L^2(\mR^m,dx)$
with a Hermitian kernel $\KK:\mR^m\times\mR^m\mapsto \mC$,
\bee\lbl{K}
K f(x)=\int_{\mR^m} \KK(x,y) f(y)\,dy, \qu K\geq 0.
\eee
Assume that $K$ is locally trace class,
i.e. for any bounded Borel set $B\subset\mR^m$
the operator
$\mI_B K \mI_B$
is trace class.
Lemmas \nolinebreak 1 and 2 from \cite{So00} imply
that it is possible to choose the kernel $\KK$ in such a way
that for any bounded Borel sets $B_1,\ldots,B_n$, $n\geq 1,$ we have
\bee\lbl{trint}
\tr \mI_{B_n}K\mI_{B_1}K\mI_{B_2}\ldots K\mI_{B_n}
= \ili_{B_1\times\ldots\times B_n} \KK(x_1,x_2)\KK(x_2,x_3)\cdots\KK(x_n,x_1)\, dx_1\ldots dx_n.
\eee
In particular, for $n=1$ we have
$\ds{\tr \mI_{B_1}K\mI_{B_1}=\ili_{B_1}\KK(x,x)\, dx.}$
Assume that (\ref{trint}) is satisfied.

\begin{defi} A random point process is called {\it determinantal}
if it has the correlation functions of the form
$$
\rho_n(x_1,\ldots,x_n)\equiv \det
\begin{pmatrix}
\KK(x_1,x_1)&\ldots& \KK(x_1,x_n) \\
\vdots& &\vdots \\
\KK(x_n,x_1) & \ldots & \KK(x_n,x_n)
\end{pmatrix} \qmb{for all }n\geq 1.
$$
\end{defi}
Determinantal processes possess the following property,
which can be viewed as their equivalent definition.
Take any bounded measurable function $h:\mR^m\mapsto \mR$
with a compact support $D:=\supp h$.
Consider the corresponding linear statistics $\SSS_h$.
Then the generating function $\MO z^{\SSS_h}$, $z\in\mC$, takes the form
\bee\lbl{Xar}
\MO z^{\SSS_h}= \det\big(1+ (z^{h}-1)K \mI_D\big),
\eee
where the expectation is taken under the determinantal process and $\det$ denotes the Fredholm determinant.
The latter is well-defined since the operator $K$ is locally trace class.

\subsection{Elementary inequalities for the trace}
\lbl{sec:basic}

We will need the following elementary inequalities.
Consider a determinantal point process on $\mR^m$, $m\geq 1$, given by a Hermitian kernel $\KK$.
Take a bounded measurable function $h: \mR^m\mapsto\mR$ with  compact support.
Set
$$
D:=\supp h
\qnd\qu
K_D:=\mI_D K \mI_D,
$$
where the integral operator $K$ is defined in (\ref{K}).
\bpp\lbl{K-K^2}
We have
\bee\non
h(K_D-K_D^2)h\geq 0.
\eee
\epp
{\it Proof.}
It is well-known that
$0 \leq K\leq \mbox{Id}$.
Consequently, $0 \leq K_D\leq \mbox{Id}$
and $K_D-K_D^2=K_D(\mbox{Id}-K_D)\geq 0$.
Then,
denoting by $\lan\cdot,\cdot\ran_{L^2}$ the scalar product in $L^2(\mR^m,dx)$,
for any function $f\in L^2(\mR^m,dx)$ we obtain
$$
\lan h(K_D-K_D^2)h f,f\ran_{L^2} = \lan(K_D-K_D^2)h f,h f\ran_{L^2}\geq 0.
$$
\qed

\ssk
It is well-known that
\begin{align}\non
\Var \SSS_h&=\ili_D h^2(x)\KK(x,x)\,dx
- \ili_D\ili_D h(x)h(y)|\KK(x,y)|^2\, dx dy, \\ \lbl{Var}
&=\tr h^2 K_D - \tr (h K_D)^2.
\end{align}
The traces above are well-defined since the operator $K$ is locally trace class,
so that the operators $h^2 K_D$ and $(h K_D)^2$ are trace class.
Denote by $[\cdot,\cdot]$ the commutator,
$[A,B]=AB-BA$.
We have
\begin{align}\non
\|[K_D,h]\|_{HS}^2&=2\ili_D\ili_D h^2(x)|\KK(x,y)|^2\, dx dy
- 2\ili_D\ili_D h(x)h(y)|\KK(x,y)|^2\, dxdy \\ \lbl{HS}
&= 2\big(\tr h^2 K_D^2 - \tr(h K_D)^2 \big).
\end{align}
\bpp
We have
\bee\lbl{trK-K^2}
0\leq \tr h^2(K_D-K_D^2) \leq \Var\SSS_{h}
\qmb{and}\qu
\|[K_D,h]\|_{HS}^2\leq 2\Var \SSS_h.
\eee
\epp
{\it Proof.}
\rprop\ref{K-K^2} together with cyclicity of the trace implies $\tr h^2(K_D-K_D^2)\geq 0$.
Next, subtracting (\ref{HS}) divided by two from (\ref{Var}), we get
$$
\Var \SSS_h - \fr12\|[h,K_D]\|_{HS}^2 = \tr h^2(K_D-K_D^2).
$$
Since $\|[h,K_D]\|_{HS}^2, \tr h^2(K_D-K_D^2)\geq 0$, we obtain (\ref{trK-K^2}).
\qed

\bpp\lbl{lem:[]}
For any linear operators $G_1,\ldots,G_n, F$, $n\geq 1$, we have
$$
[G_1\cdots G_n,F]=\sli_{l=1}^n G_1\cdots G_{l-1}[G_l,F]G_{l+1}\cdots G_n.
$$
\epp
{\it Proof.} By induction.

\section{Cumulants of linear statistics}
\lbl{sec:cumulants}

\subsection{Cumulants and traces}

In this section we
compute cumulants of linear statistics
viewed as random variables
under a determinantal process and obtain some estimates for them.
Despite that in the present paper we mainly work with the sine-process,
first in this section we consider a general determinantal process.
This is needed for the proof of Theorem \ref{th:CLT}.

Recall  that the numbers $(J_k)_{k\in\Zp}$ are called cumulants of a random vector
$\nu=(\nu_1,\ldots,\nu_d) \in\mR^d$ if for any
sufficiently small $y\in\mR^d$ we have
$$
\ln \MO e^{iy\cdot\nu}=\sli_{k\in\Zp} J_k\frac{(iy)^k}{k!},
$$
where $\cdot$ denotes the standard scalar product in $\mR^d$ while
$(iy)^k$ and $k!$ are defined in
item \nolinebreak \ref{not:!} of Section \ref{sec:Notation}.
A cumulant $J_k$ can be expressed
through the moments $(m_l)_{|l|\leq |k|}$ of the random vector $\nu$
and the other way round.
If $(e_1,\ldots,e_d)$ is the standard basis of $\mZ^d$ then
\bee\lbl{cumgeneral}
J_{e_i}=\MO \nu_i \qmb{and}\qu  J_{e_i+e_j}=\Cov(\nu_i,\nu_j) \qmb{for any } 1\leq i,j\leq d.
\eee
The vector $\nu$ is Gaussian iff $J_k=0$ for all $|k|\geq 3$.
For more information see e.g. \cite{Shi}, Section 2.12.

Let $h_1,\ldots, h_d:\mR^m\mapsto\mR$, $d\geq 1$, be
bounded Borel measurable functions with compact supports
and
$h:=(h_1,\ldots,h_d)$.
Consider the vector of linear statistics
\bee\lbl{Sh}
\SSS_h:=(\SSS_{h_1},\ldots,\SSS_{h_d})
\eee
as a random vector under a determinantal process given by a Hermitian kernel $\KK$.
Denote
\bee\lbl{DDD}
D:=\cup_{i=1}^d\supp h_i \qmb{and}\qu K_D=\mI_D K\mI_D,
\eee
where the locally trace class operator $K$ is given by (\ref{K}).
The proofs of the following \rlem\ref{lem:cum} and \rprop\ref{lem:gener} are routine (cf. formulas (1.14) and (2.7) from \cite{So00b}) and we include them for completeness.
\bl\label{lem:cum}
For any $k\in\Zp$ satisfying $|k|\geq 2$ the cumulant $B_k$ of the random vector (\ref{Sh}) has the form
\bee
\lbl{B^N_k}
B_k=k!\sli_{j=1}^{|k|}\fr{(-1)^{j+1}}{j}\sli_{\sck{a^1,\ldots,a^j\in\Zp: \\ a^1+\cdots+a^j=k}}
\fr{\tr h^{a^1} K \cdots h^{a^{j-1}} K h^{a^j} K_D}{a^1!\cdots a^j!}.
\eee
\el
{\it Proof.}
Due to (\ref{Xar}) with $z:=e^i$ and $h:=h\cdot y$, we have
$$
\ln \MO e^{i\SSS_h\cdot y}
=\ln\det\big(1+ (e^{ih\cdot y}-1) K_D \big).
$$
Then, Lemma XIII.17.6 from \cite{RS} implies that for a sufficiently small $y\in\mR^d$ we have
\begin{align}\non
\ln \MO e^{i\SSS_h\cdot y}
&=\ln\exp\Big(\sli_{j=1}^\infty\fr{(-1)^{j+1}}{j} \tr\big( (e^{ih\cdot y}-1) K_D\big)^j \Big)
\\ \lbl{cum1}
&= \sli_{j=1}^\infty\fr{(-1)^{j+1}}{j}\sli_{l^1,\ldots, l^j=1}^\infty
\fr{\tr (iy\cdot h)^{l^1} K\cdots(iy\cdot h)^{l^j} K_D}{l^1!\cdots l^j!}.
\end{align}
Note that
\bee\label{l^n}
(iy\cdot h)^{l^n}=\sum_{m_1,\ldots,m_{l^n}=1}^d iy_{m_1}h_{m_1}\cdots iy_{m_{l^n}}h_{m_{l^n}}.
\eee
We have
$$
iy_{m_1}h_{m_1}\cdots iy_{m_{l^n}}h_{m_{l^n}}
=(iy_1h_1)^{a^n_1}\cdots (iy_dh_d)^{a^n_d}=(iy)^{a^n} h^{a^n},
$$
where
\bee\label{a^n}
\qu a^n:=(a^n_1,\ldots,a^n_d)\in\Zp \qmb{and}\qu a^n_r:=\#\{q, 1\leq q\leq l^n: m_q=r\}.
\eee
Next we replace in (\ref{l^n}) the summation over $m_1,\ldots,m_{l^n}$ by that over $a^n\in\mZ^d_+$. To this end, we note that
$|a^n|=l^n$ and
 for a given vector $a^n$ the number of vectors $(m_1,\ldots,m_{l^n})$ satisfying (\ref{a^n}) is equal to $l^n!/a^n!$. Then
$$
(iy\cdot h)^{l^n}=\sli_{a^n\in\Zp: |a^n|=l^n}\fr{l^n!}{a^n!}(iy)^{a^n}h^{a^n}.
$$
Now (\ref{cum1}) implies
\begin{align}\non
\ln \MO e^{i\SSS_h\cdot y}
&=\sli_{j=1}^\infty\fr{(-1)^{j+1}}{j}\sli_{l^1,\ldots, l^j=1}^\infty
\sli_{\sck{a^1,\ldots, a^j\in\Zp: \\ |a^1|=l^1,\ldots,|a^j|=l^j}}
\fr{\tr (iy)^{a^1}h^{a^1}K\cdots (iy)^{a^j}h^{a^j}K_D}{l^1!\cdots l^j!}
\fr{l^1!\cdots l^j!}{a^1!\cdots a^j!} \\\non
&=\sli_{k\in\Zp}(iy)^k\sli_{j=1}^{|k|}\fr{(-1)^{j+1}}{j}
\sli_{\sck{a^1,\ldots, a^j\in\Zp: \\ a^1+\cdots+a^j=k}}
\fr{\tr h^{a^1}K\cdots h^{a^j}K_D}{a^1!\cdots a^j!},
\end{align}
where in the last equality the second sum is taken only over $j\leq |k|$ since for $j>|k|$ the relation $a^1+\cdots+a^j=k$ with $a^1,\ldots,
a^j\in\Zp$ is impossible.
\qed

\bpp\lbl{lem:gener}
For any $k\in\Zp$ satisfying $|k|\geq 2$ we have
\bee
\lbl{gener}
\sli_{j=1}^{|k|}\fr{(-1)^{j+1}}{j}\sli_{\sck{a^1,\ldots,a^j\in\Zp: \\ a^1+\cdots+a^j=k}}
\fr{1}{a^1!\cdots a^j!} =0.
\eee
\epp
{\it Proof.}
Denote the left-hand side of (\ref{gener}) by $T_k$.
Represent the function
\bee\lbl{combg}
g(x):=x_1+\ldots+x_d
\qmb{where}\qu
x=(x_1,\ldots,x_d),
\eee
in the form
$g(x)=\ln \big(1+(e^{x_1+\cdots+x_d}-1 )\big)$.
Developing the logarithm and exponents to the series,
we see that
$g(x)=\sli_{k\in\Zp} T_k x^k$. Indeed,
\begin{align}\non
g(x)&=      \sli_{j=1}^{\infty}\fr{(-1)^{j+1}}{j}
        \Big( \sli_{n_1=0}^\infty \fr{x_1^{n_1}}{n_1!} \cdots \sli_{n_d=0}^\infty \fr{x_d^{n_d}}{n_d!} - 1 \Big)^j
    =   \sli_{j=1}^{\infty}\fr{(-1)^{j+1}}{j}
        \Big( \sli_{n\in\Zp} \fr{x^n}{n!} \Big)^j
        \\ \non
&=   \sli_{j=1}^{\infty}\fr{(-1)^{j+1}}{j}
        \sli_{a^1,\ldots,a^j\in\Zp}\fr{x^{a^1}\cdots x^{a^j}}{a^1!\cdots a^j!}
    =   \sli_{k\in\Zp} T_k x^k.
\end{align}
Thus, due to (\ref{combg}), we have
$T_k=0$ for $|k|\geq 2$.
\qed

\ssk
Lemma \ref{lem:cum} together with Proposition \ref{lem:gener} immediately implies
\begin{cor}
\lbl{cor:cum}
For any $k\in\Zp$ satisfying $|k|\geq 2$, the cumulants $B_k$
of the random vector (\ref{Sh}) can be represented in the form
\bee
\lbl{B^N_k'}
B_k=k!\sli_{j=1}^{|k|}\fr{(-1)^{j+1}}{j}\sli_{\sck{a^1,\ldots,a^j\in\Zp: \\ a^1+\cdots+a^j=k}}
\fr{\tr h^{a^1} K\cdots h^{a^j} K_D-\tr h^k K_D}{a^1!\cdots a^j!}.
\eee
\end{cor}

In the next lemma we estimate the right-hand side of (\ref{B^N_k'}).
\bl\lbl{lem:cumest}
Let $k\in\Zp$, $|k|\geq 2$, and vectors $a^1,\ldots, a^j\in\Zp$, $j\geq 1$, satisfy $a^1+\cdots+a^j=k$. Then
\bee\lbl{cumest}
|\tr h^{a^1}K\cdots h^{a^j}K_D - \tr h^k K_D|
\leq C(|k|,d,j) \max\limits_{1\leq i\leq d} \|h_i\|^{|k|-2}_\infty \sli_{l=1}^d \Var\SSS_{h_l}.
\eee
\el
Proof of \rlem \ref{lem:cumest} follows a scheme similar to that used in the proof of \rlem 3.2 from \cite{BD15}. However, in \cite{BD15} only the
case when $K$ is a projection was considered
and the operators $h_lK$, $Kh_l$ were assumed to be of the trace class.
We do not impose these restrictions.

{\it Proof.}
{\it Step 1.} We argue by induction.
If $j=1$ then the left-hand side of (\ref{cumest}) is equal to zero. Consider the case $j=2$.
Using cyclicity of the trace, by a direct computation we get
$$
\tr h^{a^1}K h^{a^2}K_D
=\tr h^{a^1}K_D h^{a^2}K_D
= \frac12\tr[h^{a^1},K_D][h^{a^2},K_D] + \tr{h^kK_D^2}.
$$
Then
\bee\lbl{|tr-tr|}
|\tr h^{a^1}Kh^{a^2}K_D-\tr h^kK_D|\leq \fr12 \|[h^{a^1},K_D]\|_{HS} \|[h^{a^2},K_D]\|_{HS}
+ |\tr (h^k K_D^2 - h^kK_D)|.
\eee
We estimate the terms of the right-hand side above separately.
Set
\bee\lbl{gaaa}
\ga(x):=\sqrt{h_1^2(x)+\cdots+h_d^2(x)}.
\eee
Using the convention $\fr{0}{0}=:0$, we obtain
\bee\lbl{trKK2}
|\tr (h^k K_D - h^kK_D^2)|=\Big|\tr \fr{h^k}{\ga^2}\ga(K_D - K_D^2)\ga\Big|
\leq \Big\|\fr{h^k}{\ga^2}\Big\|_\infty \tr \ga(K_D - K_D^2)\ga,
\eee
since, due to \rprop\ref{K-K^2}, the operator $\ga(K_D - K_D^2)\ga$ is non-negative.
Clearly,
\bee\non
\Big\|\fr{h^k}{\ga^2}\Big\|_\infty\leq \max\limits_{1\leq i\leq d} \|h_i\|_\infty^{|k|-2}.
\eee
On the other hand, due to (\ref{trK-K^2}), we have
\bee\lbl{trKK21}
\tr \ga(K_D - K_D^2)\ga = \tr \ga^2(K_D - K_D^2)=
\sli_{l=1}^d\tr h_l^2(K_D-K_D^2)
\leq\sli_{l=1}^d\Var\SSS_{h_l}.
\eee
Thus,
\bee\lbl{tr3}
|\tr (h^k K_D - h^kK_D^2)|\leq \max\limits_{1\leq i\leq d}
\|h_i\|_\infty^{|k|-2}\sli_{l=1}^d\Var\SSS_{h_l}.
\eee
We now  estimate the Hilbert-Schmidt norm of the commutators from (\ref{|tr-tr|}). Due to Proposition \ref{lem:[]}, for any $b\in\Zp$ we
have
\begin{align}\non
\|[h^{b},K_D]\|_{HS}&\leq
|b|\max\limits_{1\leq i\leq d}\|h_i\|_\infty^{|b|-1}
\sli_{l=1}^d\|[h_l,K_D]\|_{HS} \\ \lbl{HS1}
&\leq C(|b|, d)\max\limits_{1\leq i\leq d}\|h_i\|_\infty^{|b|-1}
\Big(\sli_{l=1}^d\Var\SSS_{h_l}\Big)^{1/2},
\end{align}
where in the last inequality we have used the second relation from (\ref{trK-K^2}).
Now (\ref{|tr-tr|}) joined with (\ref{tr3}) and (\ref{HS1}) implies the desired estimate.

{\it Step 2.} Assume that $j\geq 3$. Denote
\bee\lbl{GGG}
G:=h^{a^1}K_Dh^{a^2}K_D\cdots h^{a^{j-3}}K_Dh^{a^{j-2}},
\eee
so that
$\tr h^{a^1}K\cdots h^{a^j}K_D=\tr h^{a^1}K_D\cdots h^{a^j}K_D
= \tr G K_Dh^{a^{j-1}}K_Dh^{a^j}K_D$
(in particular, for $j=3$ we have $G=h^1$).
It suffices to show that
\bee\lbl{estind}
|\tr G K_Dh^{a^{j-1}}K_Dh^{a^j}K_D-\tr GK_D h^{a^{j-1}+a^j}K_D|\leq C(|k|, d)\max\limits_{1\leq i\leq
d}\|h_i\|^{|k|-2}_\infty\sli_{l=1}^d\Var \SSS_{h_l}.
\eee
A direct computation gives
\begin{align}\non
\tr G K_Dh^{a^{j-1}}K_Dh^{a^j}K_D &= \tr GK_D[h^{a^{j-1}},K_D][h^{a^j},K_D]
+\tr GK_Dh^{a^{j-1}}K_D^2h^{a^j}  \\ \lbl{tr2.0}
&-\tr GK_D^2h^{a^{j-1}}K_Dh^{a^j}
+ \tr GK_D^2h^{a^{j-1}+a^j}K_D.
\end{align}
Write
\begin{align}\non
|\tr G K_Dh^{a^{j-1}}K_Dh^{a^j}K_D&-\tr GK_D h^{a^{j-1}+a^j}K_D|
\leq |\tr GK_D[h^{a^{j-1}},K_D][h^{a^j},K_D]|
\\ \lbl{tr2.1}
&+|\tr GK_Dh^{a^{j-1}}K_D^2h^{a^j}
-\tr GK_D^2h^{a^{j-1}}K_Dh^{a^j}|  \\ \non
&+|\tr GK_D^2h^{a^{j-1}+a^j}K_D-\tr GK_D h^{a^{j-1}+a^j}K_D|=: I_1+I_2+I_3.
\end{align}
We estimate the terms $I_1,I_2,I_3$ separately. We have
\bee\lbl{tr2.1'}
I_1\leq \|GK_D\|\|[h^{a^{j-1}}, K_D]\|_{HS}\|[h^{a^j},K_D]\|_{HS}.
\eee
Recalling that $0\leq K_D\leq \mbox{Id}$, we obtain
\bee\lbl{GKbound}
\|GK_D\|\leq \max\limits_{1\leq i \leq d}\|h_i\|_\infty^{|k|-|a^{j-1}|-|a^j|}.
\eee
Then the relation (\ref{HS1}) implies
\bee\non
I_1\leq C(|k|, d)\max\limits_{1\leq i \leq d}\|h_i\|^{|k|-2}
\sli_{l=1}^d \Var\SSS_{h_l}.
\eee
Next,
\begin{align} \non
I_2&\leq |\tr GK_Dh^{a^{j-1}}K_D^2h^{a^j}-\tr GK_Dh^{a^{j-1}}K_Dh^{a^j}| \\ \non
&+|\tr GK_Dh^{a^{j-1}}K_Dh^{a^j}-\tr GK_D^2h^{a^{j-1}}K_Dh^{a^j}|=:I_2^\prime+I_2^{\prime\prime}.
\end{align}
Due to (\ref{trKK21}) and (\ref{GKbound}),
\begin{align}\non
I_2^\prime&=\Big|\tr GK_D \frac{h^{a^{j-1}}}{\ga}\ga(K_D^2-K_D)\ga\fr{h^{a^j}}{\ga}\Big| \\ \non
&\leq \|GK_D\|\Big\|\fr{h^{a^{j-1}}}{\ga}\Big\|_\infty \Big\|\fr{h^{a^j}}{\ga}\Big\|_\infty
\tr \ga(K_D-K_D^2)\ga
\leq \max\limits_{1\leq i\leq d}\|h_i\|^{|k|-2}_\infty\sli_{l=1}^d\Var \SSS_{h_l}.
\end{align}
In a similar way we get the same estimate for the terms $I_2^{\prime\prime}$ and $I_3$.
Then (\ref{tr2.1}) implies (\ref{estind}).
\qed

\subsection{Cumulants under the sine-process}

In this section we assume that $\KK=\KK_{\rm{sine}}$ is the sine-kernel
given by (\ref{iKsine}).
Using its special structure we
rewrite the traces from (\ref{B^N_k})
in an appropriate way,
representing them
through the Fourier transforms $\hat h_i$.

Let $k\in\Zp$, $v=(v_1,\ldots,v_{|k|})$ and
$a^1,\ldots,a^j\in\Zp$, $j\geq 1$, satisfy $a^1+\ldots+a^j=k$.
Denote
\begin{align}
\non
\hat h^{a^1,\ldots,a^j}&(v):
=\hat h_1(v_1)\ldots \hat h_1(v_{a^1_1})
\hat h_2(v_{a^1_1+1})\ldots \hat h_2(v_{a^1_1+a^1_2})
\ldots
\hat h_d(v_{a^1_1+\ldots+ a^1_{d-1}+ 1})\ldots \hat h_d(v_{|a^1|})
\\\non
&\hat h_1(v_{|a^1|+1})\ldots \hat h_1(v_{|a^1| + a^2_1})
\ldots
\hat h_d(v_{|a^1|+ a^2_1+\ldots+ a^2_{d-1}+ 1})\ldots \hat h_d(v_{|a^1|+|a^2|})
\hat h_1(v_{|a^1|+|a^2|+1})\ldots.
\end{align}
We abbreviate the relation above as
\bee\lbl{l_i}
\hat h^{a^1,\ldots,a^j}(v) = \prod\limits_{i=1}^{|k|} \hat h_{l_{i}}(v_i),
\eee
where $l_i=r$, $1\leq r \leq d$, if
$$
i\in\cup_{s=1}^d
\big[|a^1|+\ldots+|a^{s-1}|+a^s_1+\ldots+a^s_{r-1}+1,\,
 |a^1|+\ldots+|a^{s-1}|+a^s_1+\ldots+a^s_{r}\big].
$$
Let for $j\geq 2$
\begin{align}\lbl{Ja}
J^{|a^1|,\ldots,|a^j|}(v)
:=&-\max\big(0,\sli_{i=1}^{|a^1|}v_i,\sli_{i=1}^{|a^1|+|a^2|}v_i,\ldots, \sli_{i=1}^{|a^1|+\ldots+|a^{j-1}|}v_i\big) \\ \non
&-\max\big(0,-\sli_{i=1}^{|a^1|}v_i,-\sli_{i=1}^{|a^1|+|a^2|}v_i,\ldots, -\sli_{i=1}^{|a^1|+\ldots+|a^{j-1}|}v_i\big),
\end{align}
and for $j=1$ set $J^{|a^1|,\ldots,|a^j|}:=0$.
\bpp\lbl{lem:cumF} Let $\KK=\KK_{\rm{sine}}$ and vectors
$k, a^1,\ldots,a^j\in\Zp$, $|k|\geq 2$, $j\geq 1$, satisfy $a^1+\ldots+ a^j=k$. Then
\bee\lbl{cumF}
\tr h^{a^1} K \cdots h^{a^j}K_D
=\fr{1}{(2\pi)^{|k|}}\ili_{v_1+\ldots+ v_{|k|}=0} \hat h^{a^1,\ldots,a^j}(v) \max\big(2 + J^{|a^1|,\ldots,|a^j|}(v), 0\big)  \, dS,
\eee
where $dS$ is an elementary volume of the hyperplane $v_1+\ldots+ v_{|k|}=0$,
normalized in such a way that
$dS(v_1,\ldots,v_{|k|})=dv_1\ldots dv_{|k|-1}$.
\epp
{\it Proof. }
In this proof
we always consider the kernel $\KK$
as a function of one variable
\bee\lbl{Ksine1}
\KK(x)=\fr{\sin x}{\pi x}.
\eee
Denote the trace from the left-hand side of (\ref{cumF}) by $\rm{Tr}$.

{\it Step 1.} Assume first $j=1$.
We have
$$
\mbox{Tr}
			=\tr h^k K_D
			=\ilif h^k(x) \KK(0)\, dx
			=\fr{1}{\pi}\ilif h^k(x)\, dx
			=\fr{1}{\pi}\FF(h^k)(0).
$$
Denote by $*$ the convolution
and set
$\hat h^{*k}:=\hat h_1^{*k_1}*\ldots*\hat h_d^{*k_d}$.
Changing the order of the convolutions, we obtain
$\hat h^{*k}=*_{i=1}^{|k|} \hat h_{l_i},$
where we recall that the indices $l_i$ are defined below (\ref{l_i}).
Then, using that
$\FF(fg)=(2\pi)^{-1}\hat f*\hat g$ for $\forall f,g\in L^2(\mR)$
we get
\begin{align}
\mbox{Tr}\non
			&=\fr{1}{\pi(2\pi)^{|k|-1}}\hat h^{*k}(0)
			\\ \non
			&=\fr{2}{(2\pi)^{|k|}} \ili_{\mR^{|k|-1}}
				\hat h_{l_1}(-y_1)\hat h_{l_2}(y_1-y_2)
				\ldots
				\hat h_{l_{|k|-1}}(y_{|k|-2}-y_{|k|-1})
				\hat h_{l_{|k|}}(y_{|k|-1})
				\,dy_1\ldots dy_{|k|-1}.
\end{align}
Next we change the variables,
$v_1:=-y_1$
and
for $2\leq i \leq |k|-1$
we set
$v_i:=y_{i-1}-y_{i}$.
Then, denoting
$v_{|k|}:=-v_1-\ldots v_{|k|-1}$
(so that $y_{|k|-1}=v_{|k|}$)
and passing from the integration over $\mR^{|k|-1}$
to the integration over the hyperplane
$v_1+\ldots+ v_{|k|}=0$ in $\mR^{|k|}$,
we arrive at (\ref{cumF}):
$$
\mbox{Tr}
			=\fr{2}{(2\pi)^{|k|}} \ili_{v_1+\ldots+ v_{|k|}=0}
					 \prod\limits_{i=1}^{|k|} \hat h_{l_i}(v_i)
						\,dS.
$$

{\it Step 2.} Let now $j\geq 2$.
In this step we show that
\bee\lbl{step1s}
\mbox{Tr}
            = \fr{1}{(2\pi)^{|k|}}\ili_{\mR^j}
                                                            \hat h^{*a^1} (y_1-y_2)\hat \KK (y_2)
                                                            \hat h^{*a^2} (y_2-y_3)\hat \KK (y_3)
                                                            \ldots
                                                            \hat h^{*a^j} (y_j-y_1) \hat \KK (y_1)
                                                            \,dy_1\ldots dy_j.
\eee
We have
$$
\mbox{Tr}
            = \ili_{\mR^j}
                                    h^{a^1}(x_1)\KK(x_1-x_2) h^{a^2}(x_2)\KK(x_2-x_3)
                                    \cdots
                                    h^{a^j}(x_j)\KK(x_j-x_1)\,dx_1\ldots dx_j.
$$
Note that
$\FF\big(\KK(\cdot-b)\big)(y)=\hat \KK(y)e^{-iyb},$
for $\forall b\in\mR$.
Then, using that
$\lan f,g \ran_{L^2}=(2\pi)^{-1}\lan \hat f,\hat g \ran_{L^2}$
and
$\FF(fg)=(2\pi)^{-1}\hat f*\hat g$ for $\forall f,g\in L^2(\mR)$,
and that the function $\hat \KK$ is real,
we find
\begin{align} \non
\ilif \KK(x_j-x_1)h^{a^1}(x_1)\KK(x_1-x_2)\, dx_1
        &=\fr{1}{2\pi}\ilif
                \ov{\FF\big(\KK(x_j-\cdot)\big)}(y)\FF\big(h^{a^1}(\cdot)\KK(\cdot-x_2)\big)(y)\, dy
                \\ \non
        &=\fr{1}{(2\pi)^{|a^1|+1}}\ili_{\mR^2}\hat\KK(y_1) e^{iy_1x_j} \hat h^{*a^1}(y_1-y_2)\hat \KK(y_2) e^{-i y_2 x_2} \,dy_1dy_2.
\end{align}
Thus,
\bee\non
\mbox{Tr}=
            \fr{1}{(2\pi)^{|a^1|+1}}\ili_{\mR^{j+1}}
                \hat h^{*a^1}(y_1-y_2)\hat \KK(y_2) e^{-i y_2 x_2}  h^{a^2}(x_2) \KK(x_2-x_3)\ldots h^{a^j}(x_j)
                 e^{iy_1x_j}\hat\KK(y_1)
                \,dy_1dy_2dx_2\ldots d x_j.
\eee
Since
$$
\ds{\ilif e^{-i y_2 x_2}  h^{a^2}(x_2) \KK(x_2-x_3) \,dx_2
    = \FF\big(h^{a^2}(\cdot)\KK(\cdot-x_3)\big)= \fr{1}{(2\pi)^{|a^2|}}\ilif \hat h^{*a^2}(y_2-y_3) \hat \KK(y_3) e^{-iy_3x_3} \,dy_3},
$$
we obtain
\begin{align}\non
\mbox{Tr}=
            \fr{1}{(2\pi)^{|a^1|+|a^2|+1}}\ili_{\mR^{j+1}}
                & \hat h^{*a^1}(y_1-y_2)\hat \KK(y_2)
                \hat h^{*a^2}(y_2-y_3)\hat \KK(y_3) e^{-i y_3 x_3}
                \\\non
                &h^{a^3}(x_3) \KK(x_3-x_4)
                \ldots h^{a^j}(x_j)e^{iy_1x_j}\hat\KK(y_1)
                \,dy_1dy_2dy_3dx_3\ldots d x_j.
\end{align}
Continuing the procedure, finally we arrive at the formula
\begin{align}\non
\mbox{Tr}
    = \fr{1}{(2\pi)^{|a^1|+\ldots+|a^{j-1}|+1}}\ili_{\mR^{j+1}}
        &\hat h^{*a^1}(y_1-y_2)\hat\KK(y_2)
            \cdots
            \hat h^{*a^{j-1}}(y_{j-1}-y_j)\hat \KK(y_j)
            \\ \non
        &e^{-i y_j x_j}
            h^{a^j}(x_j)e^{iy_1x_j}\hat\KK(y_1)
            \,dy_1\ldots dy_jdx_j.
\end{align}
Then, taking the integral over $x_j$, we get (\ref{step1s}).

\ssk
{\it Step 2.} Writing the convolutions from  (\ref{step1s}) explicitly,
we obtain
\begin{align}\non
\mbox{Tr}
=\fr{1}{(2\pi)^{|k|}}\ili_{\mR^{|k|}} &\prod\limits_{i=1}^{|a^1|} \hat h_{l_i}(y_i-y_{i+1})\hat \KK(y_{|a^1|+1})
\prod\limits_{i=|a^1|+1}^{|a^1|+|a^2|}  \hat h_{l_i}(y_i-y_{i+1})\hat \KK(y_{|a^1|+|a^2|+1})
\\\non
&\ldots
\prod\limits_{i=|a^1|+\ldots+|a^{j-1}|+1}^{|k|}  \hat h_{l_i}(y_i-y_{i+1})\hat \KK(y_1)
\,dy_1\ldots dy_{|k|},
\end{align}
where we set $y_{|k|+1}:=y_1$.
Introducing the variables $y:=y_1$ and $v_i:=y_i-y_{i+1}$, $1\leq i\leq |k|-1$,
and using the relation
$
y_n=y-\sli_{i=1}^{n-1} v_i,
$
we obtain
\begin{align}\non
\mbox{Tr}=&\fr{1}{(2\pi)^{|k|}}\ili_{\mR^{|k|-1}} \prod\limits_{i=1}^{|k|-1}\hat h_{l_i}(v_i)
\hat h_{l_{|k|}}(-v_1-\ldots-v_{|k|-1})
\\\non
&
\ili_{-\infty}^\infty\hat \KK(y)\hat \KK(y-\sli_{i=1}^{|a^1|}v_i) \cdots
\hat \KK(y-\sli_{i=1}^{|a^1|+\ldots + |a^{j-1}|} v_i)\, dy \,d v_1\ldots d v_{|k|-1}.
\end{align}
Denoting
$v_{|k|}=-v_{1}-\ldots -v_{|k|-1}$
and passing from the integration over $\mR^{|k|-1}$
to that over the hyperplane
$v_{1}+\ldots +v_{|k|}=0$, we find
\bee\lbl{trRR}
\mbox{Tr}=\fr{1}{(2\pi)^{|k|}}\ili_{v_1+\ldots+ v_{|k|}=0} \hat h^{a^1,\ldots,a^j}(v)
\ili_{-\infty}^\infty\hat \KK(y)\hat \KK(y-\sli_{i=1}^{|a^1|}v_i) \cdots
\hat \KK(y-\sli_{i=1}^{|a^1|+\ldots + |a^{j-1}|} v_i)\, dy \,dS.
\eee
{\it Step 3.}
Using that the Fourier transform of the sine-kernel (\ref{Ksine1}) has the form  $\hat\KK=\mI_{[-1,1]}$,
by a direct computation we find
\bee\lbl{Soshid}
\ili_{-\infty}^\infty\hat \KK(y)\hat \KK(y-\sli_{i=1}^{|a^1|}v_i) \cdots
\hat \KK(y-\sli_{i=1}^{|a^1|+\ldots + |a^{j-1}|} v_i)\, dy = \max\big(2 + J^{|a^1|,\ldots,|a^j|}(v), 0\big),
\eee
where the function $J^{|a^1|,\ldots,|a^j|}$ is defined in (\ref{Ja}).
Then (\ref{trRR}) implies (\ref{cumF}).
\qed
\brr
In the proof of \rprop \ref{lem:cumF}
we use the special structure of the sine-kernel only in {\it Step 3}.
\err

\begin{cor}\lbl{lem:VarS}
For any bounded measurable functions $h,h_1,h_2:\mR\mapsto\mR$ with compact support
under the sine-process we have
\begin{align}\lbl{varS}
\Var\SSS_{h} &= \fr{1}{4\pi^2}
\Big(
2\ili_{|s|\geq 2} |\hat h(s)|^2\,ds
+ \ili_{|s|<2} |s| |\hat h(s)|^2\,ds
\Big),\\\lbl{covarsinus}
\Cov(\SSS_{h_1},\SSS_{h_2})&=
\fr{1}{4\pi^2}
\Re\Big(
2\ili_{|s|\geq 2} \hat h_1(s)\ov{\hat h_2(s)}\,ds
+ \ili_{|s|<2} |s| \hat h_1(s)\ov{\hat h_2(s)}\,ds
\Big).
\end{align}
\end{cor}
In particular, Corollary \ref{lem:VarS} implies
\bee\lbl{var12}
\Var\SSS_{h} \leq \|h\|_{1/2}^2,
\eee
where we 
recall that the seminorm $\|\cdot\|_{1/2}$ is defined in (\ref{pairing}).

{\it Proof}. We first prove the formula \eqref{varS}. Recall that the variance $\Var\SSS_{h}$ is given by (\ref{Var}).
Since $\KK(x,x)=1/\pi$, we have
\bee\lbl{Svar1}
\tr h^2 K_D=\ilif h^2(x)\KK(x,x)\,dx=\fr{\|h\|^2_{L^2}}{\pi} = \fr{\|\hat h\|^2_{L^2}}{2\pi^2}.
\eee
On the other hand, Proposition \ref{lem:cumF} implies
$$
\tr (h K_D)^2 = \fr{1}{4\pi^2} \ili_{v_1+v_2=0}\hat h(v_1)\hat h(v_2) \max (2-|v_1|,0)\, dS.
$$
Then, using that $\hat h(-s)=\ov{\hat h(s)}$ since the function $h$ is real, we get
\bee\lbl{Svar2}
\tr (h K_D)^2
=\fr{1}{4\pi^2}\ili_{-\infty}^\infty |\hat h(s)|^2 \max(2-|s|,0) \, ds
= \fr{1}{4\pi^2}\ili_{|s|<2} |\hat h(s)|^2 (2-|s|) \, ds.
\eee
Inserting (\ref{Svar1}) and (\ref{Svar2}) into  (\ref{Var}), we find
$$
\Var\SSS_{h}
= \fr{1}{4\pi^2}
\Big(
2\|\hat h\|^2_{L^2} - \ili_{|s|<2} |\hat h(s)|^2 (2-|s|) \, ds
\Big)
=\fr{1}{4\pi^2}
\Big(
2 \ili_{|s|\geq 2} |\hat h(s)|^2  \, ds + \ili_{|s|<2} |s| |\hat h(s)|^2  \, ds
\Big).
$$
To get formula  \eqref{covarsinus}, we note that
$$
\Cov(\SSS_{h_1},\SSS_{h_2})=
\fr12\big(\Var \SSS_{h_1} + \Var \SSS_{h_2} - \Var(\SSS_{h_1}-\SSS_{h_2}) \big).
$$
Since $\SSS_{h_1}-\SSS_{h_2}=\SSS_{h_1-h_2}$, 
the formula \eqref{covarsinus} follows from the identity \eqref{varS}.
\qed

\section{Central Limit Theorems for linear statistics}
\lbl{sec:CLTgen}
In this section we prove multidimensional Central Limit Theorems
\nolinebreak \ref{th:CLT} and \ref{th:CLTj}.

\subsection{Linear statistics with growing variance: Theorem \ref{th:CLT}}
\lbl{sec:CLT}

Let $d\geq 1$ and $h_1^N,\ldots, h^N_d:\mR^m\mapsto\mR$, $N\in\mN$,
be a family of bounded Borel measurable functions with compact supports.
Consider the corresponding vector of linear statistics
$$
\SSS_{h^N}:=(\SSS_{h^N_1},\ldots,\SSS_{h^N_d})
$$
as a random vector under a determinantal process given by a Hermitian kernel $\KK^N$.
Denote by $\MON$, $\Var_N$ and $\Cov_N$ the corresponding expectation, variance and covariance.
In this section we prove the Central Limit Theorem
for the vector $\SSS_{h^N}$
under assumption that the variances $\Var_N\SSS_{h^N_j}$ grow to infinity as $N\ra\infty$.
\btt\label{th:CLT}
Assume that there exists a sequence $V_N\ra\infty$ as $N\ra\infty$, $V_N >0$,
such that the following
two conditions hold.
\begin{enumerate}
  \item For all $1\leq i,j\leq d$ there exist the limits
\bee\lbl{thcov}
\fr{\Cov_N (\SSS_{h^N_i},\SSS_{h^N_j})}{V_N}\os{N\ra\infty}\ra b_{ij} ,
\eee
for some numbers $b_{ij}$.
\item We have
\bee
\lbl{th|h|}
\max\limits_{1\leq j\leq d} \|h^N_j\|_\infty = o(\sqrt{V_N}) \ass N\ra\infty.
\eee
\end{enumerate}
Let $\xi^N\in\mR^d$ be a random vector with components
\bee\lbl{gxi}
\xi^N_j=\fr{\SSS_{h^N_j}-\MON\SSS_{h^N_j}}{\sqrt{V_N}}.
\eee
Then for the family of distributions $\DD(\xi^N)$ we have the weak convergence $\DD(\xi^N)\raw \DD(\xi)$ as $N\ra\infty$, where $\xi$ is a
centred Gaussian random vector with the covariance matrix $(b_{ij}).$
\ett
Theorem \ref{th:CLT}  generalizes results obtained in works \cite{CL,So00,SoAB,So01},
where the Central Limit Theorems for various linear statistics were established,
under the assumption that $\Var_N\SSS_{h^N_j}\ra\infty$ as $N\ra\infty$.
More precisely, in papers \cite{CL} and \cite{So00}
the Central Limit Theorem was proven in the one-dimensional setting (i.e. $d=1$)
for the linear statistics
corresponding to a family of functions $h^N$ of the form
\bee\lbl{hni}
h^N(x)=\mI_{A}(x/N),
\eee
where $A$ is a bounded Borel set, so that $\SSS_{h^N}=\#_A$.
In \cite{SoAB} the author considered the linear statistics of the same form
under the Airy and Bessel processes.
He showed that their variances $\Var\SSS_{h^N}$
have the logarithmic growth
and proved a multidimensional Central Limit Theorem (i.e. $d\geq 1$).
In \cite{So01} a one-dimensional Central Limit Theorem was established
for a general family of
bounded measurable functions $h^N$ with compact supports,
under the assumptions that
\bee\lbl{assSosh}
\|h^N\|_\infty = o\big((\Var_N \SSS_{h^N})^\eps\big)
\qmb{and} \qu \MON\SSS_{|h^N|}=O\big((\Var_N \SSS_{|h^N|})^\delta\big),
\eee
for any $\eps>0$ and some $\delta>0$.
This result can not be applied
for the linear statistics corresponding to the family of functions (\ref{hni})
under the sine, Airy and Bessel processes.
Indeed, the variance in these cases has the logarithmic growth
while the expectation grows as $N^n,\, n>0$,
so that (\ref{assSosh}) fails.
Since in Theorem \ref{th:CLT} we do not impose assumption (\ref{assSosh}),
it covers all the Central Limit Theorems above.

Note that in Theorem~\ref{th:CLT} the multidimensional case $d>1$ follows from the one-dimensional one $d=1$ by the linearity of statistics $\SSS_{h_j^N}$ in $h_j^N$. Proof of this fact literally repeats the proof of Proposition~\ref{lem:d-2} below. 
However, we establish Theorem~\ref{th:CLT} directly in the multidimensional setting because it does not change the proof.
 
{\it Proof of Theorem \ref{th:CLT}.}
The proof uses a method developed in \cite{CL} and \cite{So00},
and is based on application of Corollary \ref{cor:cum} and Lemma \ref{lem:cumest}.
Since the normal law is specified by its moments
it suffices to show that the moments of the random vector $\xi^N$
converge to the moments of  $\xi$
(see \cite{F}, page 269).
Denote by $(A^N_k)_{k\in\Zp}$ and  $(A_k)_{k\in\Zp}$
the cumulants of $\xi^N$ and $\xi$ respectively, so that
\bee\non
A_k=\left\{
\begin{array}{cl}
0 &\mbox{ if } |k|\neq 2, \\
b_{ij} &\mbox{ if } k=e_i+e_j,
\end{array}
\right.
\eee
where $(e_l)$ is the standard base of $\mZ^d$.
Since the moments can be expressed through the cumulants,
it suffices to prove that
\bee\lbl{convA}
A_k^N\ra A_k \ass N\ra\infty \qmb{for any }k\in\Zp.
\eee
In the case $|k|\leq 2$ the convergence (\ref{convA}) is clear.
Indeed, due to (\ref{cumgeneral}), we have
$$
A^N_{e_i}=0 \qmb{and}\qu
A^N_{e_i+e_j}=
\Cov_N\Big(\fr{\SSS_{h^N_i}-\MON\SSS_{h^N_i}}{\sqrt{V_N}},
\fr{\SSS_{h^N_j}-\MON\SSS_{h^N_j}}{\sqrt{V_N}}\Big)= \fr{\Cov_N(\SSS_{h^N_i},\SSS_{h^N_j})}{V_N},
$$
so that  (\ref{convA}) follows from assumption (\ref{thcov}).
It remains to study the case $|k|\geq 3$.
By definition (\ref{gxi}) of the vector $\xi^N$ we have
\bee\lbl{ABrelations}
A_k^N=\fr{B_k^N}{V_N^{|k|/2}},
\eee
where $B_k^N$ are cumulants of the random vector $\SSS_{h^N}$.
Due to Corollary \nolinebreak \ref{cor:cum} joined with Lemma \ref{lem:cumest}, we have
\bee\non
|B^N_k|\leq C\max\limits_{1\leq i\leq d} \|h^N_i\|^{|k|-2}_\infty \sli_{l=1}^d \Var_N\SSS_{h^N_l}.
\eee
Then, assumptions (\ref{thcov}) and (\ref{th|h|}) imply
$
B^N_k= o(V_N^{|k|/2}),
$
if $|k|\geq 3$.
Now the desired convergence (\ref{convA})
follows from (\ref{ABrelations}).
\qed

\subsection{Joint linear statistics of growing and bounded variances: Theorem \ref{th:CLTj}}
\lbl{sec:CLTj}

Consider a family of measurable bounded functions with compact supports
$f^N_1,\ldots,f^N_p,$ $g^N_1,\ldots,g^N_q:\mR\mapsto\mR$,
where $N\in\mN$ and $p,q\geq 0$.
In this section we prove a multidimensional Central Limit Theorem  \ref{th:CLTj}
for the vector of the linear statistics
\bee\lbl{jSSS}
(\SSS_{f^N_1},\ldots,\SSS_{f^N_p}, \SSS_{g^N_1},\ldots, \SSS_{g^N_q}),
\eee
under the sine-process.
We assume that the functions $f^N_i$
are as in \rtheo \ref{th:CLT}
while the functions $g^N_j$ are supposed to be sufficiently regular
and for large $N$ asymptotically behave as $\bm g_j^\infty(\cdot/N)$,
for some functions $\bm g_j^\infty$ independent from $N$.
This situation is not covered by Theorem \ref{th:CLT}
since under our hypotheses
the variances $\Var\SSS_{g^N_j}$
do not grow at all,
so that condition (\ref{thcov}) fails.

Before formulating our assumptions
let us note that all of them except $f.1$
are automatically satisfied if
$f^N_i(x)=f_i\big(x/N\big),$
$g^N_j(x)=g_j\big(x/N\big),$
where the functions $f_i,\; g_j$ are bounded measurable with compact supports
and $g_j$ belong to the Sobolev space $H^{1/2}(\mR)$.
For the proof of this fact see Example \nolinebreak \ref{ex:jCLT} in the next section.

We assume that
there exist sequences $V_N, R_N\ra\infty$ as $N\ra\infty$,
$V_N,R_N>0$,
such that  for all $1\leq i \leq p$, $1\leq j\leq q$,
the the following hypotheses hold.
Let
\bee\lbl{boldg}
\bm f^N_i(x):=f^N_i(R_Nx) \qmb{and}\qu \bm g^N_j(x):=g^N_j(R_Nx).
\eee
\begin{description}
\item[f.1] \emph{Under the sine-process there exist the limits}
\bee\lbl{jconvbf}
\fr{\Cov (\SSS_{f^N_i},\SSS_{f^N_j})}{V_N}\ra b^f_{ij} \ass N\ra\infty,
\eee
\emph{for some numbers $b^f_{ij}$ and any $1\leq i,j\leq p$.}

\item[f.2] \emph{We have $\max\limits_{1\leq i\leq p}\|f_i^N\|_\infty = o(\sqrt{V_N})$ as $N\ra\infty$.}

\item[f.3]
\emph{We have $\max\limits_{1\leq i\leq p}\|\bm f^N_i\|_{L_2}=o(\sqrt{V_N})$ as $N\ra\infty$}.
\end{description}
Since 
$\|\bm f^N_i\|_{L_2}=R_N^{-1/2}\|f_i^N\|_{L_2}$, 
assumption $f.3$ just means that the norm $\|f_i^N\|_{L_2}$ 
grows slower than $(R_NV_N)^{1/2}$.
\begin{description}
\item[g.1]
\emph{The functions $g_j^N$ belong to the Sobolev space
$H^{1/2}(\mR)$
and $\bm g_j^N\ra \bm g_j^\infty$ as $N\ra\infty$ in $H^{1/2}(\mR)$,
for some functions $\bm g_j^\infty$ and any $j$. }

\item[g.2]
\emph{The functions $g_j^N$ are bounded uniformly in $N$. }

\end{description}
Before stating the theorem let us note that, due to the estimate (\ref{var12})
and the following obvious proposition,
assumption \emph{g.2} implies in particular that
\bee\lbl{Var<C}
\mbox{the variances $\Var\SSS_{g_j^N}$ are bounded uniformly in $N$.}
\eee

\bpp\lbl{lem:g1/2}
For any function $k\in H^{1/2}(\mR)$ and any $\de\neq 0$ we have
$\|k\|_{1/2}=\|k_\delta\|_{1/2}$,
where $k_\de(x):=k(\de x).$
\epp
\emph{Proof.}
Since $\hat k_\de(x)=\de^{-1}\hat k(\de^{-1} x)$,
we get
$$
4\pi^2 \|k_\delta\|_{1/2}^2  = \de^{-2}\ilif|v||\hat k (\de^{-1}v)|^2\, dv
=\ilif|u||\hat k (u)|^2\, du=4\pi^2\|k\|_{1/2}^2,
$$
where we set $u=\de^{-1}v$.
\qed

\ssk

\btt\lbl{th:CLTj}
Let a family of measurable bounded compactly supported functions
$f^N_1,\ldots,f^N_p$, $g^N_1,\ldots,g^N_q$, $p,q\geq 0$,
satisfies assumptions
\emph{f.1}-\emph{f.3},
\emph{g.1}-\emph{g.2} above.
Consider the vector of linear statistics (\ref{jSSS})
as a random vector under the sine-process.
Let $\xi^N=(\xi_f^N,\xi_g^N)\in\mR^{p+q}$ be a random vector with components
\bee\lbl{lstatth}
\xi^N_{f_j}=\fr{\SSS_{f^N_j}-\MO\SSS_{f^N_j}}{\sqrt{V_N}}, \qquad
\xi^N_{g_i}=\SSS_{g^N_i}-\MO\SSS_{g^N_i}.
\eee
Then
$\DD(\xi^N)\raw \DD(\xi)$ as $N\ra\infty$,
where $\xi=(\xi_f,\xi_g)\in\mR^{p+q}$
is a centred Gaussian random vector with the covariance matrix
$$\begin{pmatrix}
(b^f_{ij})&0 \\
0&(b^g_{kl})
\end{pmatrix},$$
where
$b^g_{kl}=\lan \bm g_k^\infty,\bm g_l^\infty \ran_{1/2}.$
In particular, the components $\xi_f$ and $\xi_g$ of the vector $\xi$ are independent.
\ett
\rtheo\ref{th:CLTj} applied to the functions (\ref{ex:jCLTfunc}) implies \rtheo\ref{th:CLTweak}
stated in \rsec\ref{sec:iCLT}.
If $q=0$ then \rtheo \ref{th:CLTj} is covered by \rtheo \ref{th:CLT},
while in the case $p=0,q=1$ it is proven
by Spohn \cite{Sp} and Soshnikov \cite{So00b,So01} in slightly less generality;
see the discussion in \rsec\ref{sec:iCLT}.
The case $p=q=1$ was not considered before and is the main novelty of the theorem,
while the general situation $p,q\geq 0$
follows from the three just mentioned cases 
(see Proposition~\ref{lem:d-2} below). 

Proof of \rtheo \ref{th:CLTj}  employs a method developed
in \cite{So00b} mixed with that related to the method used in the proof of \rtheo \ref{th:CLT}.

Note that the required regularity $H^{1/2}$ of the functions $g_i^N$
is optimal:
if we replace $H^{1/2}$ by $H^{1/2-\eps}$
then assertion of the theorem
will not be true any more.
Indeed, the indicator function $\mI_{[0,N]}$
belongs to the space $H^{1/2-\eps}$, for all $\eps>0$.
But the linear statistics
$\SSS_{\mI_{[0,N]}}=\#_{[0,N]}$ has (logarithmically) growing variance,
so that the indicator $\mI_{[0,N]}$ belongs to the class of functions $f_i^N$ but not $g_j^N$.

\subsection{Examples}
In this section we present two examples where assumptions \emph{f.2}-\emph{g.2}
\footnote{Here and below by \emph{f.2}-\emph{g.2} we mean 
\emph{f.2},\emph{f.3},\emph{g.1},\emph{g.2}.}
are satisfied.
We will use them in \rsec\ref{sec:FCLT},
when proving our main results, Theorems \nolinebreak \ref{th:FCLTS} and \ref{th:ergintFCLTS}.
\begin{example}\lbl{ex:jCLT}
Let
\bee\lbl{ex:jCLTfunc}
f^N_i(x)=f_i\Big(\fr{x}{N}\Big), \qu g^N_j(x)=g_j\Big(\fr{x}{N}\Big), \quad
1\leq i \leq p,\; 1\leq j \leq q,
\eee
where the functions $f_i,\; g_j$ are bounded measurable with compact supports
and $g_j$ belong to the Sobolev space $H^{1/2}(\mR)$.
Then assumptions \emph{f.2}-\emph{g.2} are fulfilled
with $R_N=N$, arbitrary sequence $V_N$
and
$
\bm g^\infty_j=g_j.
$
\end{example}
{\it Proof.} Assumptions \emph{f.2} and \emph{g.2} are obviously satisfied.
Fulfilment of assumptions \emph{f.3} and \emph{g.1}  immediately follows from
the fact that, due to (\ref{ex:jCLTfunc}),
we have $\bm f^N_i=f_i$ and $\bm g^N_j=g_j$.
\qed

\begin{example}
\lbl{ex:conv}
Assume that functions
$f^N_i$, $g^N_j$
satisfy assumptions \emph{f.2}-\emph{g.2}.
Take a bounded measurable function $\ph$ with compact support
such that $\ili_{-\infty}^\infty \ph(x)\,dx = 1$.
Then the functions
$$f_{\ph,i}^N:=\ph*f^N_i, \qu g_{\ph,j}^N:=\ph*g^N_j$$
also satisfy \emph{f.2}-\emph{g.2}
with the same sequences $V_N$, $R_N$ and
functions $\bm g^\infty_j$.
\end{example}
{\it Proof.}
Assumption \emph{f.2} follows from the identity
$$
\|f_{\ph,i}^N\|_\infty \leq  \|f_i^N\|_\infty \ili_{-\infty}^\infty |\ph(x)| \,dx  .
$$
Assumptions \emph{g.2} follows in the same way.
To get assumption \emph{f.3} we define the functions $\bm f^N_{\ph,i}$
as in (\ref{boldg}) and note that
$\hat {\bm f}^N_{\ph,i}(v)=\hat\ph(v/R_N)\hat {\bm f}^N_i(v)$.
Then
$$
\|{\bm f}^N_{\ph,i}\|_{L^2}= \fr{1}{\sqrt{2\pi}}\|\hat\ph\big(\cdot/R_N\big)\hat {\bm f}^N_i\|_{L^2}
\leq \fr{1}{\sqrt{2\pi}} \|\hat\ph\|_\infty\|\hat {\bm f}^N_i\|_{L^2}
= \|\hat\ph\|_\infty\|{\bm f}^N_i\|_{L^2}.
$$
Since $\ph\in L^1(\mR)$, we have
$\|\hat\ph\|_\infty<\infty$,
so that assumption \emph{f.3} follows.

The fact that the functions $g_{\ph,j}^N$ belong to
the space $H^{1/2}(\mR)$ is implied by the inequality
$$
\|g_{\ph,j}^N\|_{1/2}
\leq \|\hat\ph\|_\infty \|g_{j}^N\|_{1/2},
$$
which can be obtained similarly to the argument above.
To establish the convergence claimed in assumption \emph{g.1},
it suffices to show that
$\|\bm g_{\ph,j}^N-\bm g_j^N\|_{H^{1/2}}\ra 0$ as $N\ra\infty$.
Using that
$\hat {\bm g}^N_{\ph,j}(v) = \hat\ph(R_N^{-1} v)\hat {\bm g}^N_{j}(v)$
and
$\hat \ph (0)=\ilif\ph(x)\,dx=1$,
we obtain
\begin{align} \non
2\pi\|\bm g_{\ph,j}^N-\bm g_j^N\|_{H^{1/2}}^2
&=\ili_{-\infty}^\infty
\big(1+|v|\big)\big|\hat \ph \big(R_N^{-1}v\big) -  \hat \ph (0) \big|^2
\big|\hat {\bm g}^N_{j}(v)\big|^2\, dv
\\\non
&\leq
\max_{|v|\leq\sqrt{R_N}}
\big|\hat \ph \big(R_N^{-1}v\big) -  \hat \ph (0) \big|^2
\ili_{|v|\leq\sqrt{R_N}}
\big(1+|v|\big)\big|\hat {\bm g}^N_{j}(v)\big|^2\,dv
\\\non
&+
2\|\hat\ph\|^2_\infty
\ili_{|v|\geq\sqrt{R_N}}
\big(1+|v|\big)\big|\hat {\bm g}^N_{j}(v)\big|^2 \,dv.
\end{align}
Using  assumption \emph{g.1} for the functions $g_j^N$, the continuity of the function $\hat\ph$
and the relation $\|\hat\ph\|_\infty<\infty$,
we see that
both of the summands above go to zero as $N\ra\infty$.
\qed
\subsection{Beginning of the proof of Theorem \ref{th:CLTj} }
The rest of \rsec\ref{sec:CLTgen} is devoted to the proof of \rtheo \ref{th:CLTj}.
From now on we will skip the upper index $N$ in the notation
$f_i^N$, $g_j^N$, $\bm f_i^N$, $\bm g_j^N$.
We start by noting that, due to the linearity of statistics $\SSS_{f_i}$, $\SSS_{g_j}$ 
in $f_i,g_j$, it suffices to establish the theorem in the case $p,q\leq 1$, so when there is at most one function $f_i$ and one function $g_j$.
\bpp\lbl{lem:d-2}
Assume that Theorem 4.3 is established in the case $p=q=1$. 
Then it holds for any $p,q\geq 0$.
\epp
 \rprop\ref{lem:d-2} is proven in \rsec\ref{sec:prop}. Below we assume $p=q=1$.
 This simplification is not crucial: 
 with minor modifications, the proof we present below suits as well for the case of arbitrary $p,q$. 
 However, under this assumption the notation become simpler.
 
 Further on we skip the lower index $1$, so for the functions $f_1,g_1,\bm f_1, \bm g_1$ and 
 the numbers $b_{11}^f$, $b_{11}^g$
 we write
 $f,g,\bm f,\bm g, b^f, b^g$.
 To prove the theorem,
 it suffices to show that the cumulants 
 $(A_k^N)_{k\in\mZ^{2}_+}$  of the random vector $\xi^N$ satisfy
 \bee\lbl{jcumconv}
 A_k^N\ra A_k \qmb{as} \qu N\ra\infty,
 \eee
 where $k=(k_f,k_g)\in\mZ^{2}_+$
 and
 \bee\non
 A_k=\left\{
 \begin{array}{cl}
 	b^f&\qmb{if $k_f=2,\, k_g=0$}, \\
 	b^g&\qmb{if $k_f=0,\, k_g=2$}, \\
 	0&\qmb{otherwise}.
 \end{array}
 \right.
 \eee
By the definition (\ref{lstatth}) of the vector $\xi^N$,
for $|k|=1$ we have $A^N_k=0$ and for $|k|\geq 2$
\bee\lbl{jAB}
A^N_k=\fr{B^N_k}{(V_N)^{k_f/2}},
\eee
where $B^N_k$ are cumulants of the random vector $(\SSS_{f},\SSS_{g})$.
Further on we assume $|k|\geq 2$.
We single out four cases:
$k_g=0$;
$k_g\geq 1$ and $k_f\geq 3$;
$k_g\geq 1$ and $k_f=2$;
$k_g\geq 1$ and $k_f\leq 1$.
The last one turns out to be the most complicated,
so we study it separately in the next subsection.
The reason is that in this case
the denominator in (\ref{jAB}) grows too slowly
or does not grow at all,
so that estimates for the cumulants $B^N_k$
like those we use to study the other cases, do not suffice
in this situation
to prove the convergence $A^N_k\ra 0$ for $|k|\geq 3$.
Instead, we employ combinatorial techniques developed by Soshnikov in \cite{So00b}.
Note that we use the special form of the sine-kernel
only in this last case.

{\it Case 1: $k_g=0$.}
In this situation convergence (\ref{jcumconv})
is established in the proof of
Theorem \nolinebreak \ref{th:CLT}.
Indeed,
the cumulant $A^N_k$ in the present case
coincides with the cumulant $A^N_{k_f}$
of the random variable $\xi_{f}^N$.

\ssk
{\it Case 2: $k_g\geq 1$ and $k_f\geq 3$.}
 Set 
\bee\lbl{jh^N}
h=(h_1,h_2):=(f,g).
\eee
In view of  \rcor \ref{cor:cum},
the desired convergence immediately follows
from (\ref{jAB}) joined with the following proposition.
\bpp\lbl{lem:cumest'}
In the case $k_f \geq 3$ (while $k_g$ is arbitrary),
for any $a^1,\ldots,a^j\in \mZ^2_+$, $j\geq 1$,
satisfying $a^1+\ldots+a^j=k$,
we have
\bee\non
\tr h^{a^1}K\ldots h^{a^j}K_D - \tr h^k K_D = o(V_N^{k_f/2})
\ass N\ra\infty.
\eee
\epp
\rprop \ref{lem:cumest'} is obtained as a refinement
of \rlem \ref{lem:cumest},
adapted for the present situation.
Its proof is given in Section \nolinebreak \ref{sec:prop}.

{\it Case 3: $k_g\geq 1$ and $k_f=2$.}
Consider a partition $k=a^1+\ldots +a^j$ from
\rcor\ref{cor:cum}.
Let $a^i=(a^i_f,a^i_g)\in\mZ^{2}_+$,
so that $k_f=a^1_f+\ldots +a^j_f$ and $k_g=a^1_g+\ldots +a^j_g$.
Since $k_f=2$, there are only two possible situations:
\begin{description}
\item[S1] There is $1\leq l\leq j$
such that $a^l_f=k_f$ and for all $i\neq l$ we have $a^i_f=0$.

\item[S2] There are $1\leq l_1<l_2\leq j$
such that $a^{l_1}_f=a^{l_2}_f=1$,
while for all $i\neq l_1,l_2$ we have $a^i_f=0$.
\end{description}
\bpp\lbl{lem:jcum}
In the situation {\it S1} above we have
\bee\lbl{j1}
\tr h^{a^1}K\ldots h^{a^j}K_{D} -\tr h^k K_{D}
= o(V_N) \ass N\ra\infty.
\eee
In the situation {\it S2},
\bee\lbl{j2}
\tr h^{a^1}K\ldots h^{a^j}K_{D} -\tr fg^{k_g} K fK_{D}
= o(V_N).
\eee
\epp
Proof of Proposition \ref{lem:jcum} is given in Section \ref{sec:prop}.
Assume that a sequence $(\wid B^N_k)_{N\in\mN}$ satisfies
\bee\lbl{widB}
B_k^N - \wid B_k^N = o(V_N).
\eee
Then, in view of (\ref{jAB}) and equality $k_f=2$, we have
$$
\lim\limits_{N\ra\infty} A_k^N = \lim\limits_{N\ra\infty} \fr{\wid B^N_k}{V_N},
$$
in the sense that if one of the limits exists
then the other exists as well and the two are equal.
Due to Corollary \ref{cor:cum} joined with \rprop \ref{lem:jcum}, the choice
\bee\lbl{jcumB}
\wid B^N_k =
k! \big(
\tr f g^{k_g}K fK_{D} - \tr h^k K_D
\big)
\sli_{j=2}^{|k|}\fr{(-1)^{j+1}}{j}
\sli_{\sck{a^1,\ldots,a^j\in\mZ^2_+ \\ \mbox{\footnotesize{satisfying S2}} : \\ a^1+\ldots+a^j=k }}
\fr{1}{a^1!\ldots a^j!}
\eee
satisfies (\ref{widB}).
Then, to prove that $A_k^N\ra 0$ as $N\ra\infty,$
it suffices to show that the sum from the right-hand side of (\ref{jcumB}) vanishes, i.e.
$$
L_k:=\sli_{j=2}^{|k|}\fr{(-1)^{j+1}}{j}
\sli_{\sck{a^1,\ldots,a^j\in\mZ^2_+ \\ \mbox{\footnotesize{satisfying S2}} : \\ a^1+\ldots+a^j=k }}
\fr{1}{a^1!\ldots a^j!}=0.
$$
Let us subtract $L_k$ from the both sides of identity (\ref{gener}).
Using that
$|k|=k_g+k_f=k_g+2$,
we find
\begin{align}\non
L_k&=
-\sli_{j=1}^{k_g+2}
\fr{(-1)^{j+1}}{j}
\sli_{\sck{a^1,\ldots,a^j\in\mZ^2_+ \\ \mbox{\footnotesize{satisfying S1}} : \\ a^1+\ldots+a^j=k }}
\fr{1}{a^1!\ldots a^j!}
=
\sli_{j=1}^{k_g+1}
\fr{(-1)^{j}}{j}
\sli_{l=1}^j
\sli_{\sck{a^1,\ldots,a^j\in\mZ^2_+: \\ a^1+\ldots+a^j=k, \\ a_f^l=k_f}}
\fr{1}{a^1!\ldots a^j!}
\\ \lbl{Lk1}
&=
\sli_{j=1}^{k_g+1}
(-1)^{j}
\sli_{\sck{a^1,\ldots,a^j\in\mZ^2_+: \\ a^1+\ldots+a^j=k, \\ a_f^j=k_f}}
\fr{1}{a^1!\ldots a^j!}.
\end{align}
In the last sum from (\ref{Lk1})
the $f$-components $a^1_f,\ldots,a^j_f$ are defined uniquely,
$a^1_f=\ldots=a^{j-1}_f=0$ and $a^j_f=k_f$.
Then we can pass
from the summation over $a^1,\ldots,a^j\in\mZ^{2}_+$
to that over $a^1_g,\ldots, a^j_g$,
where $a^1_g,\ldots,a^{j-1}_g>0$
and $a_g^j\geq 0.$
Using that
$a^1!\ldots a^j!=k_f!a^1_g!\ldots a^j_g!$,
we obtain
\bee\lbl{LLLkkk}
L_k=
\sli_{j=1}^{k_g+1}
(-1)^{j}
\sli_{\sck{a^1_g,\ldots, a^{j-1}_g>0,\, a^j_g\geq 0: \\
\, a_g^1+\ldots+a_g^j=k_g }}
\fr{1}{k_f!a^1_g!\ldots a^j_g!}.
\eee
Next we separate the last sum from (\ref{LLLkkk}) into two parts,
over $a_g^1,\ldots,a_g^j$
such that $a^j_g=0$ and such that $a^j_g\neq 0$.
We find
\begin{align}\lbl{LkLk}
L_k=
\sli_{j=1}^{k_g+1}
(-1)^{j}
\Big(
\sli_{\sck{a^1_g,\ldots, a^{j-1}_g>0: \\ a^1_g+\ldots +a^{j-1}_g=k_g}}
\fr{1}{k_f!a^1_g!\ldots a^{j-1}_g!}
+\sli_{\sck{a^1_g,\ldots, a^j_g>0: \\ a^1_g+\ldots +a^{j}_g=k_g}}
\fr{1}{k_f!a^1_g!\ldots a^j_g!}
\Big).
\end{align}
Denote
$$
x_j:=\sli_{\sck{a^1_g,\ldots, a^j_g>0:\\  a^1_g+\ldots +a^{j}_g=k_g}}
\fr{1}{k_f!a^1_g!\ldots a^j_g!}.
$$
Since in the case $j=k_g+1$ the set
$\{a^1_g,\ldots, a^j_g>0: a^1_g+\ldots+a^j_g=k_g\}$ is empty,
the relation (\ref{LkLk}) takes the form
$$
L_k=-x_1+\sli_{j=2}^{k_g}(-1)^j(x_{j-1}+x_j)+(-1)^{k_g+1}x_{k_g}=0.
$$
This finishes the consideration of {\it Case 3}.

\subsection{Conclusion of the proof of Theorem \ref{th:CLTj}}

Here we consider the last case,
when $|k|\geq 2$,
\bee\lbl{case4}
k_g\geq 1 \qnd k_f\leq 1.
\eee
We will need the following
<<smoothing>> proposition which is 
established in
\rsec\ref{sec:prop}.
\begin{prop}
	\lbl{lem:smoothing}
	Assume that \rtheo\ref{th:CLTj}
	is proven when the assumption $g.1$
	is replaced by a
	stronger assumption
	\begin{description}
		\item[g.1$^{\bm\prime}$]
		The function $g$ belong to the
		Sobolev space $H^1(\mR)$
		and $\bm g\ra\bm g^\infty$ as $N\ra\infty$
		in $H^1(\mR)$, for some function $\bm g^\infty$.
	\end{description}
	Then it holds under the
	assumption  $g.1$
	as well.
\end{prop}
Further on we assume that the function $g$
satisfies condition $g.1'$.
Recall that the function $h$ is given by \eqref{jh^N}.
Set
$$
\bm h=(\bm h_1,\bm h_2):=(\bm f,\bm g). 
$$
Due to \rprop \ref{lem:cumF},
for any $a^1,\ldots,a^j\in\mZ^2_+$, $j\geq 1$,
satisfying $a^1+\ldots+a^j=k$,
the trace
$\tr h^{a^1}K\ldots h^{a^j}K_{D}$
has the form (\ref{cumF}).
Since $\hat h_i(s)=R_N\hat{\bm h}_i(R_Ns)$, the change of variables $u_l:=R_Nv_l$ transforms
(\ref{cumF}) to
\bee\lbl{1Bk}
\tr h^{a^1}K\ldots h^{a^j}K_{D} =
\fr{1}{(2\pi)^{|k|}}\ili_{u_1+\ldots+ u_{|k|}=0} F_N^{a^1,\ldots,a^j}(u)  \, dS,
\eee
where
\bee\lbl{jF}
F_N^{a^1,\ldots,a^j}(u):=
\hat {\bm h}^{a^1,\ldots,a^j}(u)
\max\big(0, 2R_N + J^{|a^1|,\ldots,|a^j|}(u)\big),
\eee
and the function $\hat {\bm h}^{a^1,\ldots,a^j}$ is defined as in (\ref{l_i}),
with $\hat h_{l_i}$ replaced by $\hat {\bm h}_{l_i}$.
In the present case, when $p=q=1$ and $k_f\leq 1$, the function 
$\hat {\bm h}^{a^1,\ldots,a^j}$ 
has a simplified form. 
To explain this, assume first 
$k_f=1$. 
Then there exists a unique $1\leq m\leq j$ such that $a_f^m=1$ while for $i\neq m$ we have $a_f^i=0$. 
Thus, $a^m=(1,a_g^m)$ 
and $a^i=(0,a_g^i)$ for $i\neq m$.
Set 
$$
A_f:=|a^1|+\ldots + |a^{m-1}|+1=a_g^1+\ldots+a_g^{m-1}+1.
$$
Then we have
$\hat {\bm h}_{l_{A_f}}=\hat {\bm f}$
and
$
\hat {\bm h}_{l_i}=\hat {\bm g}
$
for any $i\neq A_f$.
To cover the case $k_f=0$ we set 
\bee\lbl{hat_phi}
\hat{\bm \phi}:=
\left\{
\begin{array}{cl}
	\hat{\bm f} \quad&\mbox{if}\quad k_f=1, \\
	\hat{\bm g} \quad&\mbox{if}\quad k_f=0.
\end{array}
\right.
\eee
Then we obtain
\bee\lbl{hsimple}
\hat {\bm h}^{a^1,\ldots,a^j}(v)=
\hat{\bm \phi}(v_{A_f})\prod
\limits_{\sck{1\leq i \leq |k|,\\i\neq A_f }}
\hat{\bm g}(v_i),
\eee
where in the case $k_f=0$ we choose $A_f$ arbitrary. 
In particular, we see that if $k_f=1$ the function
$\hat {\bm h}^{a^1,\ldots,a^j}$
depends on the vectors $a^1,\ldots,a^j$ 
only through the number $A_f$ while in the case $k_f=0$ it is independent from $a^1,\ldots,a^j$.

Due to \rlem \ref{lem:cum} joined with (\ref{1Bk}), the cumulant $B^N_k$ takes the form
\bee\lbl{jB_k^N}
B^N_k
=\fr{k!}{(2\pi)^{|k|}}
\sli_{j=1}^{|k|}\fr{(-1)^{j+1}}{j}\sli_{\sck{a^1,\ldots,a^j\in\mZ^2_+: \\ a^1+\cdots+a^j=k}}
\fr{1}{a^1!\cdots a^j!}\ili_{u_1+\ldots+ u_{|k|}=0} F_N^{a^1,\ldots,a^j}(u) \, dS.
\eee
Let us split it into two components,
$$
B^N_k= B_{k,1}^N+B_{k,2}^N,
$$
where
\bee\lbl{jsplitB}
B^N_{k,1} =
\fr{k!}{(2\pi)^{|k|}}
\sli_{j=1}^{|k|}\fr{(-1)^{j+1}}{j}
\sli_{\sck{a^1,\ldots,a^j\in\mZ^2_+: \\ a^1+\cdots+a^j=k}}
\fr{1}{a^1!\cdots a^j!}
\ili_{\sck{u_1+\ldots+ u_{|k|}=0,\\ |u_1|+\ldots+|u_{|k|}|\leq R_N}}
F_N^{a^1,\ldots,a^j}(u)\, dS
\eee
and
\bee\non
B^N_{k,2} = \fr{k!}{(2\pi)^{|k|}}
\sli_{j=1}^{|k|}\fr{(-1)^{j+1}}{j}
\sli_{\sck{a^1,\ldots,a^j\in\mZ^2_+: \\ a^1+\cdots+a^j=k}}
\fr{1}{a^1!\cdots a^j!}
\ili_{\sck{u_1+\ldots+ u_{|k|}=0,\\ |u_1|+\ldots+|u_{|k|}|\geq R_N}}
F_N^{a^1,\ldots,a^j}(u)\, dS.
\eee
To finish the proof of the theorem it suffices to
check that under the assumption (\ref{case4})
assertions $B1$ and $B2$ below are satisfied
(note that, in fact, the proof of $B1$ does not use assumption (\ref{case4})).

{\bf B1.} We have
\bee\lbl{jB1=0}
B_{k,1}^N=0 \qmb{if}\qu |k|>2 \qu \forall N,
\eee
and if $|k|=2$,
\bee\lbl{jB1}
\fr{B^N_{k,1}}{V_N^{k_f/2}}\os{N\ra\infty}\ra
\left\{
\begin{array}{cl}
b^g &\mbox{ if } k_f=0, k_g=2, \\
0 &\mbox{ if } k_f=k_g=1.
\end{array}
\right.
\eee

{\bf B2.} We have
\bee\lbl{jA2}
\fr{B_{k,2}^N}{V_N^{k_f/2}}\os{N\ra\infty}\ra 0.
\eee
{\bf Proof of B1.}
For $|u_1|+\ldots+|u_k|\leq R_N$ we have
$\max\big(0, 2R_N + J^{|a^1|,\ldots,|a^j|}(u)\big)=2R_N + J^{|a^1|,\ldots,|a^j|}(u).$
So that,
\bee\lbl{INa}
F_N^{|a^1|,\ldots,|a^j|}(u)=\hat {\bm h}^{a^1,\ldots,a^j}(u)\big( 2R_N + J^{|a^1|,\ldots,|a^j|} (u)\big).
\eee
Then the integral from (\ref{jsplitB}) takes the form $I_1+I_2$,
where
$$
I_1:=2R_N\ili_{\sck{u_1+\ldots+ u_{|k|}=0,\\ |u_1|+\ldots+|u_{|k|}|\leq R_N}}
\hat {\bm h}^{a^1,\ldots,a^j}\, dS
\qnd
I_2:= \ili_{\sck{u_1+\ldots+ u_{|k|}=0,\\ |u_1|+\ldots+|u_{|k|}|\leq R_N}}
        \hat {\bm h}^{a^1,\ldots,a^j} J^{|a^1|,\ldots,|a^j|}\, dS.
$$
Changing the order in the product \eqref{hsimple},
we obtain
$$
I_1=2R_N\ili_{\sck{u_1+\ldots+ u_{|k|}=0,\\ |u_1|+\ldots+|u_{|k|}|\leq R_N}}
\hat {\bm \phi}(u_1)
\hat {\bm g}(u_{2})\cdots \hat {\bm g}(u_{|k|})
\, dS.
$$
Thus, the integral $I_1$ is independent from the choice of
the vectors $a^i$,
so that in the formula (\ref{jsplitB})
it can be put in front of the sums.
In view of  \rprop \ref{lem:gener} the sums vanish,
so that only the integral $I_2$ has an input to the term $B_{k,1}^N$:
\bee\lbl{BNk1}
B^N_{k,1}
=\fr{k!}{(2\pi)^{|k|}}
\sli_{j=1}^{|k|}\fr{(-1)^{j+1}}{j}
\sli_{\sck{a^1,\ldots,a^j\in\mZ^2_+: \\ a^1+\cdots+a^j=k}}
\fr{1}{a^1!\cdots a^j!}
\ili_{\sck{u_1+\ldots+ u_{|k|}=0,\\ |u_1|+\ldots+|u_{|k|}|\leq R_N}}
\hat {\bm h}^{a^1,\ldots,a^j}(u) J^{|a^1|,\ldots,|a^j|}(u) \, dS.
\eee
Denote by $\Sigma_{|k|}$ the symmetric group of degree $|k|$.
\bpp\lbl{lem:jpermut}
Let $l^1,\ldots,l^j\in\mN\sm\{0\}$, $j\geq 1$, satisfy $l^1+\ldots +l^j=|k|$. Then
\begin{align}\lbl{jpermut}
&\sli_{\sck{a^1,\ldots,a^j\in\mZ^2_+: \\ a^1+\cdots+a^j=k, \\|a^1|=l^1, \ldots, |a^j|=l^j}}
\fr{1}{a^1!\ldots a^j!}
\ili_{\sck{u_1+\ldots+ u_{|k|}=0,\\ |u_1|+\ldots+|u_{|k|}|\leq R_N}}
\hat {\bm h}^{a^1,\ldots,a^j}(u) J^{l^1,\ldots,l^j}(u) \, dS
\\\non
&=
\fr{1}{k!l^1!\ldots l^j!} \sli_{\sigma\in \Sigma_{|k|}}\;
\ili_{\sck{u_1+\ldots+ u_{|k|}=0,\\ |u_1|+\ldots+|u_{|k|}|\leq R_N}}
\hat {\bm \phi}(u_1)
\hat {\bm g}(u_{2})\cdots \hat {\bm g}(u_{|k|})
J^{l^1,\ldots,l^j}(u^\sigma)\,dS,
\end{align}
where $u^\sigma:=(u_{\sigma(1)},\ldots,u_{\sigma(|k|)})$.
\epp
Proof of Proposition \ref{lem:jpermut} is postponed to Section \ref{sec:prop}.
In view of the definition (\ref{Ja}) of the function
$J^{l^1,\ldots,l^j}$,
Proposition \ref{lem:jpermut} applied to (\ref{BNk1})
implies
\bee\lbl{jB1main}
B^N_{k,1}=\fr{1}{(2\pi)^{|k|}}\ili_{\sck{u_1+\ldots+ u_{|k|}=0,\\ |u_1|+\ldots+|u_{|k|}|\leq R_N}}
\hat {\bm \phi}(u_1)
\hat {\bm g}(u_{2})\cdots \hat {\bm g}(u_{|k|})
\big(G(u)+G(-u)\big)\,dS,
\eee
where
$$
G(u):=\sli_{j=1}^{|k|}\fr{(-1)^{j}}{j}
\sli_{\sck{l^1,\ldots,l^j >0: \\ l^1+\cdots+l^j=|k|}}
\sli_{\sigma\in \Sigma_{|k|}}
\fr{1}{l^1!\ldots l^j!}
\max\big(0,\sli_{i=1}^{l^1}u_{\sigma(i)},\sli_{i=1}^{l^1+l^2}u_{\sigma(i)},\ldots, \sli_{i=1}^{l^1+\ldots+l^{j-1}}u_{\sigma(i)}\big).
$$
The Main Combinatorial Lemma from \cite{So00b} (see page 1356 in  \cite{So00b}) states that for any real numbers
$u_1,\ldots, u_{|k|}$ satisfying $u_1+\ldots+ u_{|k|}=0$ we have
\footnote{It seems that in \cite{So00b} the factor 2 is omitted: it is written $G(u)$ instead of $2G(u)$.}
\bee\lbl{jSochG}
2G(u)=\left\{
\begin{array}{cl}
|u_1|=|u_2| &\qmb{if}\qu |k|=2, \\
0						&\qmb{if}\qu |k|>2.
\end{array}
\right.
\eee
Then, (\ref{jB1main}) implies (\ref{jB1=0}).
Now it remains only to study the term $B_{k,1}^N$ in the case $|k|=2$.

{\it Case $k_f=0$
and $k_g=2$}.
Since the function $\bm g$ is real,
we have $\hat {\bm g}(-s)\equiv \ov{\hat {\bm g}(s)}$.
Then, in view of (\ref{jB1main}) and (\ref{jSochG}), we get
\bee\lbl{jB1final}
B_{k,1}^N
=\fr{1}{(2\pi)^2}\ili_{\sck{u_1+u_{2}=0, \\ |u_1|+|u_2|\leq R_N}} \hat {\bm g}(u_1)\hat {\bm g}(u_2) |u_1| \, dS
= \fr{1}{4\pi^2} \ili_{-R_N/2}^{R_N/2} |s|\hat {\bm g}(s)\ov{\hat {\bm g}(s)}\, ds.
\eee
Due to assumption $g.1'$ (even $g.1$ suffices here),
the right-hand side of (\ref{jB1final}) converges to $b^g$,
so that we get (\ref{jB1}).

{\it Case $k_f=k_g=1$}. 
Relation (\ref{jB1main}) joined with (\ref{jSochG})  implies
$$
A_{k,1}^N:=
\fr{B^N_{k,1}}{V_N^{k_f/2}}
=
\fr{1}{4\pi^2\sqrt{V_N}} \ili_{-R_N/2}^{R_N/2}
|s|\hat {\bm f}(s)\ov{\hat {\bm g}(s)}\, ds.
$$
Using the Cauchy-Bunyakovsky-Schwarz inequality, we obtain
\bee\lbl{jA1}
|A_{k,1}^N|
\leq
\fr{1}{4\pi^2}
\Big(
\fr{1}{V_N} \ili_{-R_N/2}^{R_N/2} |\hat{\bm f}(s)|^2 \,ds
\Big)^{1/2}
 \Big(
\ili_{-R_N/2}^{R_N/2} |s|^2|\hat{\bm g}(s)|^2 \,ds
\Big)^{1/2}.
\eee
Due to assumption $f.3$, the first integral above
goes to zero as $N\ra\infty$.
Since, in view of assumption $g.1'$,
the second one is bounded uniformly in $N$,
the desired convergence follows.

{\bf Proof of B2.}
Since, by the definition,
$J^{|a^1|,\ldots, |a^j|}\leq 0$,
we have
$|F_N^{|a^1|,\ldots, |a^j|}|\leq 2R_N |\hat {\bm h}^{a^1,\ldots,a^j}|$,
see (\ref{jF}).
Thus, it suffices to prove that
\bee\lbl{3Fh}
V_N^{-k_f/2}
R_N\ili_{\sck{u_1+\ldots+ u_{|k|}=0,\\ |u_1|+\ldots+|u_{|k|}|\geq R_N}}
|\hat {\bm h}^{a^1,\ldots,a^j}(u)| \, dS \ra 0 \ass N\ra\infty,
\eee
for any $a^1,\ldots,a^j$.
Due to \eqref{hsimple},
changing the order in the product 
$\hat {\bm h}^{a^1,\ldots,a^j}$
we obtain
$$
|\hat {\bm h}^{a^1,\ldots,a^j}(u)|
	\leq	|\hat {\bm \phi}(u_1)\hat {\bm g}(u_2)\hat {\bm g}(u_3)\ldots\hat {\bm g}(u_{|k|})|.
$$
Excluding the variable $u_1$, we get
\begin{align}\non
R_N&\ili_{\sck{u_1+\ldots+ u_{|k|}=0,\\ |u_1|+\ldots+|u_{|k|}|\geq R_N}}
|\hat {\bm h}^{a^1,\ldots,a^j}(u)| \, dS
	\\\non
	&\leq
		R_N
		\ili_{|u_2|+\ldots+|u_{|k|}|\geq R_N/2}
		\big|\hat {\bm \phi}(-u_2-\ldots-u_{|k|})\hat {\bm g}(u_2)\ldots\hat {\bm g}(u_{|k|})\big|
		\,du_2\ldots du_{|k|}
\\\non
	&\leq
		2\ili_{|u_2|+\ldots+|u_{|k|}|\geq R_N/2} (|u_2|+\ldots+|u_{|k|}|)
		\big|\hat {\bm \phi}(-u_2-\ldots-u_{|k|})\hat {\bm g}(u_2)\ldots\hat {\bm g}(u_{|k|})\big|
		\,du_2\ldots du_{|k|}
\\\non
	&=2(|k|-1)\ili_{|u_2|+\ldots+|u_{|k|}|\geq R_N/2} |u_2|
		\big|\hat {\bm \phi}(-u_2-\ldots-u_{|k|})
	\hat {\bm g}(u_2)\ldots\hat {\bm g}(u_{|k|})\big|
		\,du_2\ldots du_{|k|}
\\\non
	&=:2(|k|-1) L^N.
\end{align}
Next we separate the cases
$|k|=2$ and $|k|>2$.

{\it Case $|k|=2$.}
Applying the  Cauchy-Bunyakovsky-Schwarz inequality, we obtain
\bee\non
L^N=\ili_{|u_2|\geq R_N/2}
		\big|u_2\hat{\bm \phi}(-u_2)\hat{\bm g}(u_2)\big| \, du_2
	\leq
		\|\hat{\bm \phi}\|_{L_2}\Big(\ili_{|u_2|\geq R_N/2}|u_2|^2|\hat{\bm g}(u_2)|^2\, du_2\Big)^{1/2}.
\eee
Assumptions $f.3$ and $g.1'$ (or $g.1$) imply that
\bee\lbl{phil2}
\|\hat{\bm \phi}\|_{L_2}\leq C V_N^{k_f/2}.
\eee
Then, using assumption $g.1'$, we find
$V^{-k_f/2}L^N\ra 0$ as $N\ra\infty$.
So that, we get (\ref{3Fh}).

{\it Case $|k|>2$.}
We have
\begin{align}\non
L^N
	&\leq\ili_{|u_3|+\ldots+|u_{|k|}|\geq R_N/4}
		\big|\hat {\bm g}(u_3)\ldots \hat {\bm g}(u_{|k|})\big|
		\ili_{|u_2|\leq R_N/4}
		\big|u_2	\hat {\bm \phi}(-u_2-\ldots-u_{|k|})\hat {\bm g}(u_2)\big|
		\,du_2\ldots du_{|k|}
	\\\non
	&+\ili_{\mR^{|k|-2}}
		\big|\hat {\bm g}(u_3)\ldots \hat {\bm g}(u_{|k|})\big|
		\ili_{|u_2|\geq R_N/4}
		\big|u_2	\hat {\bm \phi}(-u_2-\ldots-u_{|k|})\hat {\bm g}(u_2)\big|		\,du_2\ldots du_{|k|}
		\\\non
		&=:L^N_1+L^N_2.
\end{align}
Using the Cauchy-Bunyakovsky-Schwarz inequality, we find
\begin{align}\non
L^N_1
	&\leq
		C\|\hat {\bm \phi}\|_{L^2} \|\hat {\bm g}\|_{1}
		\ili_{|u_3|+\ldots+|u_{|k|}|\geq R_N/4}
		\big|\hat {\bm g}(u_3)\ldots \hat {\bm g}(u_{|k|})\big|\,du_2\ldots du_{|k|}
\qnd
\\\non
L_2^N
	&\leq
		\|\hat {\bm g}\|_{L^1}^{|k|-2}\|\hat {\bm \phi}\|_{L^2}
		\Big(\ili_{|u_2|\geq R_N/4} |u_2|^2 |\hat {\bm g}(u_2)|^2\,du_2\Big)^{1/2}.
\end{align}
In view of (\ref{phil2}) and assumption $g.1'$,
to see that
$V_N^{-k_f/2}L^N_2\ra 0$ as $N\ra\infty$
it suffices to show that the $L^1$-norm
$\|\hat {\bm g}\|_{L_1}$ is finite and bounded uniformly in $N$.
This follows from the estimate
\bee\lbl{harm}
\ilif |\hat {\bm g}(s)| \,ds
=\ilif\fr{|s|+1}{|s|+1}|\hat {\bm g}(s)|\, ds
\leq C\|\hat {\bm g}\|_{H^{1}}
\Big(\ilif \fr{ds}{(|s|+1)^{2}}\Big)^{1/2}
=C_1\|\hat {\bm g}\|_{H^{1}}.
\eee
To see that $V_N^{-k_f/2}L^N_1\ra 0$,
we need additionally prove that
the integral
$
\ili_{|s|>M} |\hat {\bm g}(s)| \,ds
$
converges to zero as $M\ra\infty$ uniformly in $N$.
This follows similarly.

\subsection{Proofs of auxiliary results}
\lbl{sec:prop}

In this section we establish Propositions \ref{lem:d-2}-\ref{lem:jpermut}
used in the proof of \rtheo \ref{th:CLTj}.

{\bf Proof of Proposition \ref{lem:d-2}.}
Let the functions $f_1,\ldots f_p$, $g_1,\ldots,g_q$ satisfy
assumptions $f.1$-$g.2$.
Set
$$\mB:=
\begin{pmatrix}
(b^f_{ij}) & 0 \\
0& (b^{g}_{kl})
\end{pmatrix},$$ 
so that $\mB$ is a $(p+q)\times(p+q)$-matrix.
To establish the proposition, it suffices to prove that for any $(t,s)\in\mR^{p+q}$ we have the convergence of the characteristic functions
\bee\lbl{xard-2}
\MO e^{i(\xi_f^N\cdot t+\xi_g^N\cdot s )}
\to
e^{-\frac12(t,s)\mB(t,s)^T}
\ass N\to\infty,
\eee
where $T$ stands for the transposition.
Note that 
$$
\xi^N_f\cdot t=\sli_{j=1}^p t_j\frac{\SSS_{f_j}-\MO\SSS_{f_j}}{\sqrt{V_N}}
=\frac{\SSS_{f\cdot t}-\MO\SSS_{f\cdot t}}{\sqrt{V_N}}=:\xi^N_{f\cdot t}.
$$
Similarly, we have $\xi^N_{g\cdot s}:=\xi_g^N\cdot s$.
Thus, it remains to study the "2-dimensional" characteristic function
$\MO e^{i(\xi^N_{f\cdot t}+\xi^N_{g\cdot s})}$.
We claim that the functions $f\cdot t$ and $g\cdot s$ satisfy assumptions of Theorem~\ref{th:CLTj}.
Indeed, assumptions $f.2$, $f.3$ and $g.2$ are obviously fulfilled. 
Assumption $g.1$ is also obvious, with  the function 
$\bm g^\infty\cdot s=\sli_{j=1}^N \bm g_j^\infty\cdot s_j$. 
Thus, it remains to check assumption $f.1$. We have
$$
\frac{\Var\SSS_{f\cdot t}}{V_N}=\frac{\sli_{i,j=1}^pt_i t_j\Cov{(\SSS_{f_i},\SSS_{f_j})}}{V_N}
\to \sli_{i,j=1}^p t_i t_j b_{ij}^f=t(b^f_{ij})t^T,
\ass N\to\infty,
$$
due to assumption $f.1$ for the functions $f_1,\ldots f_p$. 
So, the function $f\cdot t$ also satisfies $f.1$.
Here we emphasize that in Theorem~\ref{th:CLTj} we do not assume that $b_{ij}^{f,g}$ are different from zero (cf. Remark~\ref{rem:zerocovdelta}), 
so that the situation when $t(b^f_{ij})t^T=0$ 
does not pose problems.

Now, applying Theorem~\ref{th:CLTj} to the functions $f\cdot t$ and $g\cdot s$, we get the convergence
$
\DD(\xi_{f\cdot t}^N,\xi_{g\cdot s}^N)\raw \DD(\xi_{f\cdot t},\xi_{g\cdot s}),
$
where $(\xi_{f\cdot t},\xi_{g\cdot s})$ is a centred Gaussian random vector with the covariance matrix
$$\mB_2:=
\begin{pmatrix}
t(b^f_{ij})t^T & 0 \\
0& \|\bm g^\infty\cdot s\|^2_{1/2}
\end{pmatrix}.
$$
In particular, we have the convergence of the corresponding characteristic functions 
at the point $(1,1)$
$$
\MO e^{i(\xi^N_{f\cdot t}+\xi^N_{g\cdot s})}\to e^{-\frac12 (1,1)\mB_2(1,1)^T}.
$$
Now, to get \eqref{xard-2} it remains to note that $(1,1)\mB_2(1,1)^T=(t,s)\mB(t,s)^T$.
\qed

\ssk

\ssk
{\bf Proof of Proposition \ref{lem:cumest'}.}
We follow the scheme used in the proof of \rlem \ref{lem:cumest}.
Assume first that $j=2$ (the case $j=1$ is trivial).
Then we have estimates (\ref{|tr-tr|}) and (\ref{trKK2}).
Assumptions \emph{f.2} and \emph{g.2}  state that
\bee\lbl{fginfty}
\|f\|_\infty = o(\sqrt{V_N})
\qnd
\|g\|_\infty \leq C.
\eee
Then
$$
\Big\|\fr{h^k}{\ga^2}\Big\|_\infty
\leq C\|f\|_\infty^{k_f-2}
= o(V_N^{\fr{k_f-2}{2}}),
$$
where $\ga$ is defined in (\ref{gaaa}) with $d=2$.
Then, in view of (\ref{trKK2}) and (\ref{trKK21}),
we have
\bee\lbl{2trtr}
|\tr(h^kK_D -h^kK_D^2)|\leq
o(V_N^{\fr{k_f-2}{2}})(\Var\SSS_{f}+\Var\SSS_g)
= o(V_N^{k_f/2}).
\eee
Here we have used that
$\Var\SSS_{f}\leq CV_N$
and
$\Var\SSS_{g}\leq C$,
accordingly to assumption \emph{f.1} and (\ref{Var<C}).

Take any $c=(c_f,c_g)\in\mZ^{2}_+.$
If $c_f=0$ then,
due to Proposition \ref{lem:[]} joined with (\ref{trK-K^2}),
we have
\begin{align}\lbl{bf=0}
\|[h^c,K_D]\|_{HS}\leq
C\|g\|^{|c|-1}_\infty
\|[g,K_D]\|_{HS}
\leq
C_1\|g\|^{|c|-1}_\infty
(\Var\SSS_{g})^{1/2}\leq C_2.
\end{align}
If $c_f=1$ then
\begin{align}\lbl{bf=1}
\|[h^c,K_D]\|_{HS}&\leq
C\|f\|_\infty\|g\|^{|c|-2}_\infty
(\Var\SSS_{g})^{1/2}
+
C\|g\|^{|c|-1}_\infty
(\Var\SSS_{f})^{1/2}
\\\non
&\leq
o(\sqrt{V_N}) + C_1\sqrt{V_N}\leq C_2 V_N^{c_f/2}.
\end{align}
If $c_f\geq 1$, then arguing similarly we find
\bee\lbl{bf>1}
\|[h^c,K_D]\|_{HS}
=
o(V_N^{c_f/2}).
\eee
Take $a^1,a^2\in\mZ^2_+$ satisfying $a^1+a^2=k$.
Since $k_f\geq 3$,
the situation $a^1_f,a^2_f\leq 1$ is impossible.
Then, without loss of generality we assume that $a^2_f>1$ and get
\bee\lbl{3trtr}
\|[h^{a^1},K_D]\|_{HS}\|[h^{a^2},K_D]\|_{HS}
\leq CV_N^{a^1_f/2} o(V_N^{a^2_f/2})
= o(V_N^{k_f/2}).
\eee
Now estimates (\ref{2trtr}) and (\ref{3trtr}) imply that
the right-hand side of (\ref{|tr-tr|}) is bounded by $o(V_N^{k_f/2})$,
so that we get the desired inequality.
The case $j\geq 3$ can be studied in a similar way,
following the scheme of the proof of \rlem \ref{lem:cumest}.

\ssk
{\bf Proof of Proposition \ref{lem:jcum}.}
To get the desired estimates we revise the proof of \rlem \ref{lem:cumest},
additionally using the regularity of the function $g$
and the fact that the operator $K$ corresponds to the sine-kernel,
so that $K$ is a projection: $K^2=K$.
The latter relation will be used in estimates
analogous to (\ref{|tr-tr|}) and (\ref{tr2.1}),
to kill there the second and the second and the third terms
of the right-hand side correspondingly.
The problem here is that $K^2_D\neq K_D$.
Thus, our first aim is to reduce estimates on the operator $K_D$ to estimates on $K$.

Let $m\geq 2$ and $b^1,\ldots, b^m\in\mZ^2_+$. Since the supports $\supp h_i$ are compact
and the sine-kernel $\KK$ has the form (\ref{iKsine}),
the operators $Kh^{b^i}$ and $h^{b^i}K$ are Hilbert-Schmidt.
This implies that the operator $h^{b^1}K\ldots h^{b^m}K$ is of the trace class
as a product of Hilbert-Schmidt operators.
Jointly with cyclicity of the trace this provides
\bee\lbl{drop}
\tr h^{b^1}K\ldots h^{b^m}K_D = \tr h^{b^1}K\ldots  h^{b^m}K \mI_D
= \tr \mI_D h^{b^1}K \ldots h^{b^m}K  = \tr h^{b^1}K \ldots h^{b^m}K.
\eee
Thus, it suffices to establish relations (\ref{j1})-(\ref{j2}), where the index $D$
is dropped everywhere except the term $\tr h^kK_{D}$
(we do not know if the operator $h^kK$ is of the trace class,
so we can not argue as in (\ref{drop}) to drop the index $D$ there).

Now let us consider the situations \emph{S1} and \emph{S2}
separately.
If $j=1$, we are automatically in the case {\it S1}
and the left-hand side of (\ref{j1}) vanishes.
Further on we assume that $j\geq 2$.

{\it Case S1}. Let first $j=2$.
Since $a^1+a^2=k$,
we have
\bee\lbl{commutprod}
[h^{a^1},K] [h^{a^2},K]
= h^{a^1}Kh^{a^2}K + Kh^{a^1}Kh^{a^2} - Kh^{k}K -h^{a^1}K^2 h^{a^2}.
\eee
Using that for any Hilbert-Schmidt operators $A,B$ we have $\tr AB=\tr BA$
and that $K=K^2$, we obtain
$$
\tr Kh^{k}K= \tr K \mI_D h^{k} K
=  \tr  h^{k} KK \mI_D =\tr  h^{k} K_D.
$$
Similarly, $\tr h^{a^1}K^2 h^{a^2}=\tr  h^{k} K_D$
and
$\tr Kh^{a^1}Kh^{a^2}= \tr h^{a^1}Kh^{a^2}K$.
Then, due to (\ref{commutprod}),
$$
\tr [h^{a^1},K] [h^{a^2},K]
= 2\tr h^{a^1}Kh^{a^2}K - 2 \tr  h^{k} K_D.
$$
Consequently,
\bee\lbl{jHS2}
\big|\tr h^{a^1}K h^{a^2}K -\tr h^k K_D \big| \leq
\fr12\|[h^{a^1},K]\|_{HS}\|[h^{a^2},K]\|_{HS}.
\eee
Arguing as in (\ref{HS}) we find
$$
\|[h_i,K]\|_{HS}^2
=2\tr h_i^2 K^2 \mI_D-2\tr (h_iK_D)^2
=2\tr h_i^2 K_D-2\tr (h_iK_D)^2
=2\Var \SSS_{h_i},
$$
for $i=1,2,$
where in the last equality we have used (\ref{Var}).
Now for definiteness we assume that $a^1_f=k_f$ and $a^2_f=0$.
Then, similarly to (\ref{bf=0}) and  (\ref{bf>1}), we get
$$
\|[h^{a^2},K]\|_{HS}\leq C \qmb{and}\qu
\|[h^{a^1},K]\|_{HS}= o\big(V_N^{k_f/2}\big)=o(V_N).
$$
Thus, (\ref{jHS2}) implies the desired inequality (\ref{j1}).

Assume now $j\geq 3$.
Define the operator $G$ as in (\ref{GGG}) with $K_D$ replaced by $K$.
Then, literally repeating (\ref{tr2.0})-(\ref{tr2.1'}) with $K_D$ replaced by $K$ and using the identity $K^2=K$, we get
\bee\lbl{jjgeq3}
\big|\tr GK h^{a^{j-1}}K h^{a^j}K -\tr GK h^{a^{j-1}+a^j} K \big| \leq
\|GK\|\|[h^{a^{j-1}},K]\|_{HS}\|[h^{a^j},K]\|_{HS}.
\eee
Without loss of generality we assume that
$a^1_f=k_f$ and $a^i_f=0$ for $i\geq 2$
(in particular, $a^{j-1}_f=a^j_f=0$).
Then, arguing as above, we see that the Hilbert-Schmidt norms from (\ref{jjgeq3}) are bounded uniformly in $N$.
On the other hand,
\bee\lbl{GKj}
\|GK\|\leq \|f\|_\infty^{k_f-a^{j-1}_f-a^j_f}\|g\|_\infty^{k_g-a^{j-1}_g-a^j_g}
= o(V_N^{k_f/2})=o(V_N),
\eee
due to (\ref{fginfty}).
Thus, the right-hand side of (\ref{jjgeq3}) is bounded by $o(V_N)$.
Now, by the induction axiom, we get the desired inequality (\ref{j1}).

{\it Case S2}.
Without loss of generality we assume that
$a_f^1=a_f^n=1$ for some $n>1$, while for $i\neq 1,n$ we have $a_f^i=0$.
Consider first the situation when $j\geq 3$. Then $j,j-1\neq 1$.
If additionally $j,j-1\neq n$, then the norms
$\|[h^{a^{j}},K]\|_{HS}$ and $\|[h^{a^{j-1}},K]\|_{HS}$
are bounded uniformly in $N$.
Moreover, (\ref{GKj}) is satisfied, so that the right-hand side of inequality (\ref{jjgeq3})
is bounded by $o(V_N)$.
If one of the numbers $j-1$ or $j$ is equal to $n$,
then, arguing as in (\ref{bf=1}),
we see that
the Hilbert-Schmidt norm of the corresponding commutator
is majorated by $C\sqrt{V_N}$.
On the other hand, the product from
(\ref{GKj}) is then bounded by $o(V_N^{\fr{k_f-1}{2}})$ in this case.
Thus, the right-hand side of  (\ref{jjgeq3})
is majorated by
$C\sqrt{V_N} o(V_N^{\fr{k_f-1}{2}})
=o(V_N).$
Summing up, for $j\geq 3$ we obtain
\bee\lbl{jGKH1}
\tr GK h^{a^{j-1}}K h^{a^j}K -\tr GK h^{a^{j-1}+a^j} K
=o(V_N).
\eee
Let now $j=2$.
Then
$\tr h^{a^1}K h^{a^2}K= \tr fg^{a_g^1}K fg^{a_g^2}K$,
for some $1\leq m_1,m_2\leq p$.
Using cyclicity of the trace, we obtain
\begin{align} \non
\big| \tr fg^{a_g^1}K fg^{a_g^2}K &- \tr fg^{k_g}K fK\big|
\leq \big| \tr fg^{a_g^1}K fg^{a_g^2}K
 - \tr fg^{a_g^1}K g^{a_g^2}K fK\big|
\\ \lbl{inducj}
&+ \big| \tr fK fg^{a_g^1}K  g^{a_g^2}K  -
\tr fK fg^{k_g}K \big|
\leq o(V_N),
\end{align}
in view of (\ref{jGKH1}). Now the desired estimate (\ref{j2}) follows by induction
from (\ref{jGKH1}) and (\ref{inducj}).
\qed

{\bf Proof of Proposition \ref{lem:smoothing}.}
Let the functions $f, g$ satisfy
assumptions $f.1$-$g.2$.
Consider a smooth function $w:\mR\mapsto\mR$
$$
w(x)=\left\{
\begin{array}{cl}
Ce^{\fr{1}{x^2-1}} \quad &\mbox{if}\qu |x|<1,
\\
0 \quad &\mbox{if}\qu |x|\geq 1,
\end{array}
\right.
$$
where the constant $C$ is chosen in such a way that
$\ilif w(x)\,dx =1$.
Set $w_\eps(x):=\eps^{-1}w(\eps^{-1} x)$, where $0<\eps< 1$,
and let
$$
\bm g_{\eps}:=w_\eps*\bm g.
$$
{\it Step 1}. In this step we show that for any $\eps>0$ the function $g_{\eps}$
defined through the function $\bm g_{\eps}$
as in (\ref{boldg})
satisfies assumptions $g.1',\,g.2$
with $\bm g_{\eps}^\infty=w_\eps*\bm g^\infty$.
Fulfilment of $g.2$
follows from the inequality
$$
\|g_{\eps}\|_\infty
=\|\bm g_{\eps}\|_\infty
\leq \|\bm g\|_\infty \ilif w_\eps(x)\,dx=
\|\bm g\|_\infty=
\|g\|_\infty.
$$
Since the function $w_\eps$ is smooth,
the function $g_{\eps}$ also is,
so in particular $g_{\eps}$ belongs to the space $H^1(\mR)$.
Then, to get assumption $g.1'$ it suffices to show that
$\bm g_{\eps}\ra \bm g_{\eps}^\infty$ as $N\ra\infty$
in $H^1(\mR)$.
We have
\bee\lbl{gepsginf}
\|\bm g_{\eps}-\bm g_{\eps}^\infty\|_{H^1}^2 =
\fr{1}{2\pi}\ilif
\big( 1+|u|^2\big)
\big|\hat w_\eps\big|^2
\big|\hat{\bm g}-\hat{\bm g}^\infty\big|^2 \,du
\leq \| \hat\al_\eps \|_\infty \|\bm g-\bm g^\infty\|_{H^{1/2}}^2,
\eee
where $\hat\al_\eps:=(1+|u|^2)(1+|u|)^{-1} |\hat w_\eps|^2$.
Since the function
$\hat w_\eps$ is of the Schwartz class,
the norm $\|\hat\al_\eps\|_\infty$ is finite 
(although dependent from $\eps$).
Then, assumption $g.1$ for the function $g$ implies
that the right-hand side of (\ref{gepsginf}) goes to zero as $N\ra\infty$,
for each $\eps>0$.

{\it Step 2}. 
It remains to show that 
assertion of \rtheo\ref{th:CLTj} holds for the
functions $f,g$.
By assumption of the proposition,
it is satisfied for the
functions $f,g_{\eps}$.
So that,
for the random vector $(\xi^N_f, \xi^N_{g_\eps})$ defined as in
(\ref{lstatth})
and any $t,s\in\mR$  we have
\bee\lbl{convsmooth}
\MO e^{i(t\xi^N_f + s\xi^N_{g_\eps})}
\ra e^{-\fr{1}{2}(t,s)\mB_\eps(t,s)^T}
\ass N\ra\infty,
\eee
where
$\mB_{\eps}:=
\begin{pmatrix}
b^f & 0 \\
0& b^{g_\eps}
\end{pmatrix}$
and
$b^{g_\eps}: = \|\bm g_{\eps}^\infty\|_{1/2}^2$.
Let 
$\mB:=
\begin{pmatrix}
b^f & 0 \\
0& b^{g}
\end{pmatrix}$.
Then
$$
\Big| \MO e^{i(t\xi^N_f + s\xi^N_{g})} - e^{-\fr{1}{2}(t,s)\mB(t,s)^T} \Big|
\leq
I^N_{1,\eps} + I^N_{2,\eps} + I_{3,\eps},
$$
where
$$
I^N_{1,\eps} =
\Big| \MO e^{i(t\xi^N_f + s\xi^N_{g_\eps})} - e^{-\fr{1}{2}(t,s)\mB_{\eps}(t,s)^T} \Big|,
\qu
I^N_{2,\eps} =
\Big| \MO e^{i(t\xi^N_f+ s\xi^N_{g})}- \MO e^{i(t\xi^N_f + s\xi^N_{g_\eps})}\Big|
$$
and
$$
I_{3,\eps} =
\Big| e^{-\fr{1}{2}(t,s)\mB_{\eps}(t,s)^T}- e^{-\fr{1}{2}(t,s)\mB(t,s)^T} \Big|.
$$
In view of (\ref{convsmooth}), $I_{1,\eps}^N\ra 0$ as $N\ra\infty$ for any $\eps>0$ and any $t,s$.
Thus, to finish the proof of the proposition
it remains to show that
$I_{2,\eps}^N,\, I_{3,\eps}\ra 0$ as $\eps\ra 0$ uniformly in $N$,
for any $t,s$.
We have
$$
I^N_{2,\eps}
\leq \MO \big| e^{is\xi^N_{g}} - e^{is\xi^N_{g_\eps}} \big|
\leq \MO |s||\xi^N_{g} - \xi^N_{g_\eps}|
\leq |s|\sqrt{\Var (\xi^N_{g}-\xi^N_{g_{\eps}})}.
$$
Due to (\ref{var12}),
$$
\Var (\xi^N_{g}-\xi^N_{g_{\eps}})
=\Var \SSS_{g-g_{\eps}}
\leq \|g-g_{\eps}\|_{1/2}^2=\|\bm g- \bm g_{\eps}\|_{1/2}^2,
$$
where in the last identity we used \rprop\ref{lem:g1/2}.
For any $r>0$ we have
\begin{align*}
4\pi^2\|\bm g- \bm g_{\eps}\|_{1/2}^2
&	=\ilif |u||1-\hat w_\eps(u)|^2|\hat{\bm g} (u)|^2 \,du
\\
&\leq \sup\limits_{|x|\leq r} |1-\hat w_\eps(x)|^2
\ili_{-r}^r |u| |\hat{\bm g}|^2 \,du
+ (\|\hat w_\eps\|_\infty+1\big)^2\ili_{|u|\geq r} |u||\hat {\bm g}|^2 \,du.
\end{align*}
Due to assumption $g.1$ for the function $g$ and
the relation $\hat w_\eps(x)=\hat w(\eps x)\ra \hat w(0)=1$ as $\eps\ra 0$, which holds for any $x$, we see that the first term above goes to zero as $\eps\ra 0$, for any $r$ uniformly in $N$.
Using assumption $g.1$ again, we find that the second term goes to zero when $r\ra\infty$,
uniformly in $\eps$ and $N$.  Consequently,
\bee\lbl{bmgunN}
\|\bm g- \bm g_{\eps}\|_{1/2}^2 \ra 0 \ass \eps\ra 0
\qmb{uniformly in $N$},
\eee
so that $I^N_{2,\eps}$ also does.
To show that $I_{3,\eps}\ra 0$ as $\eps\ra 0$,
it suffices to prove that
$\|\bm g^\infty- \bm g^\infty_{\eps}\|_{1/2}\ra 0$ as $\eps\ra 0$.
This follows by taking the limit $N\ra\infty$ in (\ref{bmgunN}).
\qed

\ssk
{\bf Proof of Proposition \ref{lem:jpermut}.}
Set
$$
{\bm \eta}_1:=\hat {\bm \phi}
\qnd
{\bm \eta}_{2}=\ldots={\bm \eta}_{|k|}:=\hat{\bm g}.
$$
Then the sum from the right-hand side of (\ref{jpermut}) can be rewritten as
\bee\lbl{jrhspermut}
\sli_{\sigma\in \Sigma_{|k|}}\,
\ili_{\sck{u_1+\ldots u_{|k|}=0 \\ |u_1|+\ldots+|u_{|k|}|\leq R_N}}
{\bm \eta}_{\sigma(1)}(u_1){\bm \eta}_{\sigma(2)}(u_2)\ldots {\bm \eta}_{\sigma(|k|)}(u_{|k|})
J^{l^1,\ldots,l^j}(u)\,dS.
\eee
Fix any partition $a^1+\ldots+a^j=k$, where $|a^i|=l^i$ for all $i.$
The function $J^{|a^1|,\ldots,|a^j|}(u)$
depends on $u$ only through the
unordered sets
$\{u_1,\ldots, u_{|a^1|}\},\{u_{|a^1|+1},\ldots u_{|a^1|+|a^2|}\},$ $\ldots$.
Then the integral from the left-hand side of
(\ref{jpermut}),
corresponding to this partition,
coincides with
the integral from
(\ref{jrhspermut}),
corresponding to a permutation
$\sigma$,
iff
among the functions
${\bm \eta}_{\sigma(1)},\ldots,{\bm \eta}_{\sigma(|a^1|)}$
there are exactly
$a^1_f$ functions equal to $\hat {\bm f}$
and
$a^1_g$ functions equal to $\hat {\bm g}$;
among the functions
${\bm \eta}_{\sigma(|a^1|+1)},\ldots,{\bm \eta}_{\sigma(|a^1|+|a^2|)}$ there are
$a^2_f$ functions equal to $\hat {\bm f}$,
$a^2_g$ functions equal to $\hat {\bm g}$, and so on.
The number of such permutations can be found directly and is equal to
$$
\fr{k!l^1!\ldots l^j!}{a^1!\ldots a^j!}.
$$
Thus, the sum (\ref{jrhspermut}) can be rewritten as
$$
\sli_{\sck{a^1,\ldots,a^j\in\mZ^2_+: \\ a^1+\cdots+a^j=k, \\|a^1|=l^1, \ldots, |a^j|=l^j}}
\fr{k!l^1!\ldots l^j!}{a^1!\ldots a^j!}
\ili_{\sck{u_1+\ldots +u_{|k|}=0 \\ |u_1|+\ldots+|u_{|k|}|\leq R_N}}
\hat {\bm h}^{a^1,\ldots,a^j}(u) J^{l^1,\ldots,l^j}(u) \, dS.
$$
\qed

\section{Proofs of main results}
\lbl{sec:FCLT}
In this section we establish Propositions \ref{lem:i2}, \ref{lem:fd},   \ref{lem:covnotfixed} and Theorems \ref{th:FCLTS}, \ref{th:ergintFCLTS}.

\subsection{Proofs of Theorem \ref{th:FCLTS} and Propositions  \ref{lem:i2},\ref{lem:fd}}
\lbl{sec:FCLTS}
Here we prove Propositions \ref{lem:fd}, \ref{lem:i2}  and Theorem \ref{th:FCLTS}.

\ssk
{\bf Proof of Proposition \ref{lem:fd}.}

The number of particles $\#_{[0,t_iN]}$
coincides with the linear statistics $\SSS_{f_i^N}$,
where $f_i^N:=\mI_{[0,t_i N]}$.
Then the desired convergence would follow from the
Central Limit \rtheo \ref{th:CLT}
if we show that
\bee\lbl{iconvvar}
\fr{\Cov(\SSS_{f_i^N},\SSS_{f_j^N})}{\pi^{-2}\ln N}
\ra b_{ij} \ass N\ra\infty,
\eee
where $b_{ij}$ are given by \eqref{iCov}.
In the case $i=j$ convergence (\ref{iconvvar})
follows from (\ref{iVar}).
Assume that $i>j$.
Since $\SSS_{f_i^N}-\SSS_{f_j^N}=\#_{[t_jN,t_iN]}$,
due to (\ref{iVar}) we have
\bee\non
\Var(\SSS_{f_i^N}-\SSS_{f_j^N}) = \pi^{-2}\ln N + O(1).
\eee
Then (\ref{iconvvar}) follows from (\ref{iVar}) and the obvious relation
\bee\lbl{covlog1}
\Cov (\SSS_{f^N_i},\SSS_{f_j^N})
= \fr12 \big( \Var\SSS_{f^N_i} + \Var\SSS_{f^N_j} - \Var(\SSS_{f^N_i}-\SSS_{f^N_j})\big).
\eee
\qed

{\bf Proof of Theorem \ref{th:FCLTS}}.

{\it Step 1.}
In this step we show that
for any $0\leq t_1<\ldots<t_d\leq 1$,
\bee\lbl{fdconvFCLTS}
\DD(\eta^N,z^N_{t_1},\ldots,z^N_{t_d}) \raw \DD(\eta,z_{t_1},\ldots,z_{t_d}) \ass N\ra\infty.
\eee

Note that
\bee\lbl{etaNzN}
\eta^N=\fr{\SSS_{f^N}-\MO \SSS_{f^N}}{\pi^{-1}\sqrt{\ln N}}
\qmb{and}\qu
z_t^N=\SSS_{g_t^N}-\MO \SSS_{g_t^N},
\eee
where $f^N(x)=f(x/N)$, $g_t^N(x)=g_t(x/N)$ and
\bee\lbl{mainfg}
f(x):=\fr{1}{\tau}\ili_{0}^{\tau} \mI_{[0,s]}(x)\,ds,
\qu
g_t(x):=\ili_0^t \mI_{[0,s]}(x) \,ds - \fr{t}{\tau}\ili_0^{\tau} \mI_{[0,s]}(x) \,ds.
\eee
The following simple result
is established in the next section.
\bpp\lbl{lem:fvar}
We have
\begin{enumerate}
\item
$\Var\SSS_{f^N} = \ds{\fr{1}{2\pi^2}} \ln N + O(1).$
\item
$g_t\in H^1(\mR),$ for any $t\in [0,1]$.
\item
$\lan g_{t},g_{s}\ran_{1/2}=\,$right-hand side of (\ref{covarz}),
for any $t,s\in [0,1]$.

\end{enumerate}
\epp
We claim that the family of functions
$f^N,g^N_{t_1},\ldots,g^N_{t_d}$
satisfies assumptions of \rtheo \ref{th:CLTj}.
Indeed,
due to \rprop \ref{lem:fvar}(1),
assumption \emph{f.1} is fulfilled
with $V_N=\pi^{-2}\ln N$
and
$b_{11}^f=1/2$.
Assumptions \emph{f.2}-\emph{g.2} are fulfilled as well
with $R_N=N$ and $ \bm g^\infty_{t_i}=g_{t_i}$
since
we are in the situation of
Example \ref{ex:jCLT},
in view of \rprop \ref{lem:fvar}(2).
Then, in due to \rprop \ref{lem:fvar}(3),
\rtheo \ref{th:CLTj} implies the convergence (\ref{fdconvFCLTS}).

{\it Step 2}. In this step we show that
the family of measures $\{\DD(\eta^N,z^N), N\in\mN\}$
is tight in the space $\mR\times C([0,1],\mR)$.
To this end, it suffices to prove that
the family of measures $\{\DD(\eta^N), N\in\mN\}$
is tight in $\mR$
and the family $\{\DD(z^N), N\in\mN\}$
is tight in $C([0,1],\mR)$.
Indeed, then for any $\eps>0$ we will be able to find
compact sets $K_\eta\subset\mR$ and $K_z\subset C([0,1],\mR)$
such that
$$
\PR(\eta^N\in K_\eta)>1-\eps/2 \qmb{and}\qu \PR(z^N\in K_z)>1-\eps/2, \qmb{for all }N.
$$
Then we will have
\begin{align}\non
\PR\big((\eta^N,z^N)\in K_\eta\times K_z\big)
&= \PR(\eta^N\in K_\eta) - \PR(\eta^N\in K_\eta, z^N\notin K_z) \\ \non
&\geq \PR(\eta^N\in K_\eta) - \PR(z^N\notin K_z) > 1-\eps.
\end{align}

Tightness of the family of measures $\{\DD(\eta^N), N\in\mN\}$
follows from convergence (\ref{fdconvFCLTS})
since the weak convergence implies the tightness.
To show that the family $\{\DD(z^N), N\in\mN\}$ is tight,
we first formulate the following proposition.
\bpp\lbl{lem:regtight}
Consider a family of bounded measurable functions with compact supports
$h^N_t:\mR\mapsto\mR$,  $0\leq t \leq 1$, $N\in\mN$.
Assume that for each $t$ and $N$
the function $h^N_t$ belongs to the Sobolev space $H^{1/2}(\mR)$.
Assume also that there exist constants $C,\de>0$
such that for any $0\leq t,s \leq 1$ and $N\in\mN$ we have
\bee\lbl{contr}
\|h^N_0\|_{1/2}\leq C,
\qu
\|h^N_t-h^N_s\|^2_{1/2}\leq C(t-s)^{1+\de}.
\eee
Consider the random process
$$
\zeta^N_t:=\SSS_{h_t^N} -\MO\SSS_{h_t^N}, \qu 0\leq t\leq 1,
$$
under the sine-process.
Then there exists a continuous modification $\wid\zeta^N_t$
of the process $\zeta_t^N$
such that the family of measures
$\{\DD(\wid\zeta^N),N\in\mN\}$ is tight
in the space of continuous functions $C([0,1],\mR)$.
\epp
Proof of \rprop \ref{lem:regtight}
is given in the next section.
Now to get the desired tightness of the family $\{\DD(z^N)\}$,
it remains only to check that
assumption (\ref{contr}) is satisfied for the functions $g^N_t$.
Its first part is obvious since $g^N_0=0.$
Using that $\hat g_t^N(u)=N\hat g_t(Nu)$, we find
\bee\lbl{gg12}
\lan g^N_t,g^N_s\ran_{1/2}
=\fr{N^2}{4\pi^2}\ili_{-\infty}^\infty |u|\hat g_t(Nu)\ov{\hat g_s(Nu)} \,du
=\fr{1}{4\pi^2}\ili_{-\infty}^\infty |v|\hat g_t(v)\ov{\hat g_s(v)} \,dv
=\lan g_t,g_s\ran_{1/2}.
\eee
Then, recalling the notation
${\ds \theta(u)=\fr{u^2}{4\pi^2} \ln |u|, \, \theta(0)=0,}$
and using \rprop \ref{lem:fvar}(3), we obtain
$$
\fr12\|g^N_t-g^N_s\|_{1/2}^2
		=\fr{t-s}{\tau}
		 \Big(\theta(t)-\theta(s)
			+\theta(s-\tau)-\theta(t-\tau) - \fr{t-s}{\tau}\tht(\tau)\Big)
			-\tht(t-s)
=:\Theta(t,s) -\tht(t-s).
$$
Since the derivative $\theta^\prime(u)$ is bounded uniformly in $u\in (0,1)$,
we have
$\Theta(t,s) \leq C(t-s)^2$.
Since $|\tht(t-s)|\leq C(\delta)(t-s)^{1+\delta}$
for any $0<\delta<1$, the second part of assumption (\ref{contr}) is satisfied as well.

{\it Step 3.} In this step
we derive the required convergence (\ref{reqconv})
from the first two steps
by standard argument.
Since the family of measures
$\{\DD(\eta^{N},z^{N}), N\in\mN\}$
is tight,
by the Prokhorov Theorem
it is weakly compact.
Take a subsequence $N_k\ra\infty$
such that
$$
\DD(\eta^{N_k},z^{N_k})\raw \DD(\wid\eta,\wid z)
\ass k\ra \infty
\qmb{in}\qu
\mR\times C([0,1],\mR),
$$
where $\DD(\wid\eta,\wid z)$
is a limit point.
Due to (\ref{fdconvFCLTS}),
for any $0\leq t_1<\ldots<t_d\leq 1$ we have
$$
\DD(\wid\eta,\wid z_{t_1},\ldots, \wid z_{t_d})
=\DD(\eta,z_{t_1},\ldots, z_{t_d}).
$$
Since finite-dimensional distributions
specify a process,
all the limit points coincide
with $\DD(\eta,z)$,
so that we get the desired convergence.
Proof of the theorem is completed.
\qed

\ssk
{\bf Proof of \rprop\ref{lem:i2}}.

Consider first a cumulant $A_k^N$ with $k\geq 3$.
Denote by $(B_m^N)$ cumulants of the random variable $\SSS_{f^N}$,
where the function $f^N$ is defined above (\ref{mainfg}).
Due to Corollary \nolinebreak \ref{cor:cum}
joined with \rlem \ref{lem:cumest},
we have
$$|B^N_k|\leq C\|f^N\|^{k-2}_\infty \Var\SSS_{f^N}.$$
Since the norm $\|f^N\|_\infty$ is independent from $N$,
\rprop \ref{lem:fvar} implies
$|B_k^N|\leq C\ln N$.
In view of  (\ref{etaNzN}), we have
$\ds{ A^N_k=\fr{B^N_k}{(\pi^{-2}\ln N)^{k/2}} }.$
Then
$$
|A^N_k|\leq \fr{C\ln N}{(\ln N)^{k/2}}=\fr{C}{(\ln N)^{k/2-1}}.
$$
Since for $k\geq 3$ cumulants $A_k$ of the normal distribution vanish,
we get the desired estimate (\ref{cumtheo}).

For $k=2$ we have
$\ds{A^N_2=\Var\eta^N=\fr{\Var\SSS_{f^N}}{\pi^{-2}\ln N}}$ and $A_2=\Var\eta=1/2$.
Then the desired estimate follows from \rprop \ref{lem:fvar}.
Since $A^N_1=\MO\eta^N=0$ and $A_1=\MO\eta=0$, the proof of the proposition  is finished.
\qed

\subsection{Proofs of auxiliary propositions}

Here we establish Propositions \ref{lem:fvar} and \ref{lem:regtight}
used in the previous section.

{\bf Proof of \rprop \ref{lem:fvar}.}
{\it Item 1.} Since $f^N=\tau^{-1}\ili_0^{\tau} \mI_{[0,sN]} \,ds$, we have
$\SSS_{f^N}=\tau^{-1}\ili_0^{\tau} \SSS_{h^N_s}  \, ds$,
where $h^N_s=\mI_{[0,sN]}$.
Then, using the Fubini theorem, we get
\bee\lbl{varintts}
\Var\SSS_{f^N}
=\MO\Big(\fr{1}{\tau}\ili_0^{\tau} \SSS_{h^N_s} - \MO \SSS_{h^N_s} \, ds\Big)^2
=\fr{1}{\tau^2}\ili_0^{\tau}\ili_0^{\tau}\Cov(\SSS_{h^N_t},\SSS_{h^N_s})\, dsdt.
\eee
Let us represent the covariance
$\Cov(\SSS_{h^N_t},\SSS_{h^N_s})$
through the variances
$\Var\SSS_{h^N_t}$, $\Var\SSS_{h^N_s}$
as in (\ref{covlog1}).
Since
$\SSS_{h_s^N}=\#_{[0,sN]}$,
we have
$\SSS_{h_t^N}-\SSS_{h_s^N}=\#_{[sN,tN]}$,
if $t>s$.
Then, due to the logarithmic grows
of the variances (\ref{iVar}), we obtain
$$
\Cov(\SSS_{h^N_t},\SSS_{h^N_s})=\fr{1}{2\pi^2} \ln N + O(1),
$$
for $t\neq s$. It can be shown that the integral  $\ili_0^{\tau}\ili_0^{\tau} O(1) \, dsdt$ is
bounded uniformly in $N$.
Now (\ref{varintts}) implies the desired relation.
\qed

\ssk
{\it Item 2.}
Calculating the integrals from (\ref{mainfg}) explicitly, we see that
the functions $g_t$ are piecewise linear and continuous, so that
$g_t\in H^1(\mR)$.
Indeed,
for $0\leq t\leq \tau $ we have
\bee\lbl{gtx}
g_t(x)=
\left\{
\begin{array}{cl}
0 &\qmb{if } x \leq 0 \mbox{ or } x\geq \tau, \\
x(\tau^{-1}t-1) &\qmb{if } 0\leq x\leq t, \\
t(\tau^{-1}x-1) &\qmb{if } t \leq x\leq \tau.
\end{array}
\right.
\eee
For $\tau \leq t\leq 1$,
\bee\lbl{gtx1}
g_t(x)=
\left\{
\begin{array}{cl}
0 &\qmb{if } x \leq 0 \mbox{ or } x\geq t, \\
x(\tau^{-1}t-1) &\qmb{if } 0\leq x\leq \tau, \\
t-x &\qmb{if } \tau \leq x\leq t.
\end{array}
\right.
\eee

{\it Item 3.}
Since $g_0= g_\tau = 0$, in the case $t=0,\tau$ or $s=0,\tau$
the result is trivial.
Assume that $t,s\neq 0,\tau$.
By a direct computation we find
$$
\hat g_t(y)=\fr{h_t(y)}{y^2}
\qmb{where}\qu
h_t(y): =1- e^{-ity} -\fr{t}{\tau}(1 - e^{-i\tau y}).
$$
Then, using that $\ov{\hat g_t(y)}=\hat g_t(-y)$, we obtain
\bee\non
\lan g_t,g_s\ran_{1/2}
	=\fr{1}{2\pi^2}
	 \Ree\ili_0^\infty
			y \hat g_t(y)\ov{\hat g_s(y)}\,dy
	=\fr{1}{2\pi^2}
	 \Ree\ili_0^\infty \fr{h_t(y)\ov{h_s(y)}}{y^3}\, dy.
\eee
Integrating by parts two times we find
\bee\lbl{app1}
\ili_0^\infty \fr{h_t\ov h_s}{y^3}\, dy
		=-\fr{h_t\ov h_s}{2y^2}\Big|_0^\infty
		 -\fr{(h_t\ov h_s)^\prime}{2y}\Big|_0^\infty
		 +\ili_0^\infty\fr{(h_t \ov{h_s})^{\prime\prime}}{2y} \,dy,
\eee
where the prime stands for the derivative with respect to $y$.
We have
\bee\lbl{app2}
h_t^\prime(y)=ite^{-ity} - ite^{-i\tau y}
\qnd
h_t^{\prime\prime}(y)=t^2 e^{-ity} - \tau t e^{-i\tau y}.
\eee
Since $h_t(0)=h_t^\prime(0)=0$ for any $t$,
we have
$(h_t\ov h_s)(0)=(h_t\ov h_s)'(0)=(h_t\ov h_s)''(0)=0$, so that
the boundary terms from (\ref{app1})
vanish.
Then, using (\ref{app1}) and (\ref{app2}), by a direct computation we find
\bee\lbl{app3}
\Ree\ili_0^\infty \fr{h_t\ov h_s}{y^3}\, dy
	=\Ree\ili_0^\infty\fr{(h_t \ov h_s)^{\prime\prime}}{2y} \,dy
	=-\ili_0^\infty \fr{v(t,s)+v(s,t)}{2y},
\eee
where
$$
v(t,s)
	=\fr{(t-s)^2}{2}\cos\big((t-s)y\big)
	-t^2\big(1-\fr{s}{\tau}\big)\cos(ty)
	-\fr{s}{\tau}\big(\tau-t\big)^2\cos\big((t-\tau)y\big)
  +t(\tau-s)\cos (\tau y).
$$
\bpp\lbl{lem:cos}
Let $a_1,\ldots a_n,b_1,\ldots,b_n\in\mR\sm\{0\}$ and
$a_1+\ldots+a_n=0$. Then
\bee\lbl{intlog}
\ili_0^\infty\sli_{i=1}^n\fr{a_i \cos (b_iy)}{y} \,dy=-\sli_{i=1}^n a_i\ln|b_i|.
\eee
\epp
Observe that the last integral from (\ref{app3})
has the form (\ref{intlog}).
Then, applying \rprop \ref{lem:cos} we obtain the desired identity.

\ssk
{\it Proof of \rprop\ref{lem:cos}.}
Since $a_1+\ldots+a_n=0$, the integral under the question converges.
Take $\eps>0$ and write
$$
\ili_0^\infty\sli_{i=1}^n\fr{a_i \cos (b_iy)}{y} \,dy
	=\ili_0^\eps+\ili_\eps^\infty
		\sli_{i=1}^n\fr{a_i \cos (b_iy)}{y} \,dy
	=:I_0^\eps+I^\infty_\eps.
$$
Clearly, $I_0^\eps\ra 0$ as $\eps\ra 0$.
On the other hand,
$$
I^\infty_\eps=\sli_{i=1}^n a_i \ili_{\eps}^\infty \fr{\cos (b_i y)}{y} \,dy
=\sli_{i=1}^n a_i \ili_{|b_i|\eps}^\infty \fr{\cos y}{y} \,dy
=\sli_{i=1}^n a_i \ili_{|b_1|\eps}^\infty \fr{\cos y}{y} \,dy
+\sli_{i=2}^n a_i \ili_{|b_i|\eps}^{|b_1|\eps} \fr{\cos y}{y} \,dy.
$$
Since $a_1+\ldots+a_n=0$, this implies
$
\ds{I^\infty_\eps=\sli_{i=2}^n a_i \ili_{|b_i|\eps}^{|b_1|\eps} \fr{\cos y}{y} \,dy.}
$
Letting $\eps$ go to zero, we obtain
$$
I^\infty_\eps \volna \sli_{i=2}^n a_i \ili_{|b_i|\eps}^{|b_1|\eps} \fr{1}{y} \,dy
= -\sli_{i=2}^n a_i \ln\fr{|b_i|}{|b_1|} = -\sli_{i=1}^n a_i \ln{|b_i|}.
$$
\qed

\ssk
{\bf Proof of \rprop \ref{lem:regtight}.}
Due to the Kolmogorov-\v{C}entsov Theorem (see Theorem 2.8 in \cite{KaSh})
and Problem 2.4.11 from \cite{KaSh},
to prove the proposition
it suffices to show that
$$
(1)\;
\sup\limits_{N\in\mN}\MO|\zeta^N_0|^{2}<\infty
\qquad
(2)\;
\sup\limits_{N\in\mN}\MO(\zeta_t^N -\zeta_s^N)^2 \leq C(t-s)^{1+\de}
\mbox{ uniformly in } 0\leq s,t\leq 1.
$$
We have
\bee\lbl{Sg-g}
\MO(\zeta^N_t-\zeta^N_s)^2 = \Var\SSS_{h^N_t-h^N_s}.
\eee
Due to estimate (\ref{var12}) of Corollary \nolinebreak \ref{lem:VarS},
the right-hand side of (\ref{Sg-g}) is bounded by
$\|h^N_t-h^N_s\|^2_{1/2}$.
Then assumption (\ref{contr}) implies
item {\it (2)} above.
Assertion of item {\it (1)} follows in the same way,
$$
\MO |\zeta_0^N|^2 = \Var\SSS_{h_0^N} \leq \|h^N_0\|^2_{1/2}\leq C.
$$
\qed

\subsection{Proofs of Theorem \ref{th:ergintFCLTS} and Proposition \ref{lem:covnotfixed}}
\lbl{sec:ergintFCLTS}

{\bf Proof of Theorem \ref{th:ergintFCLTS}}

{\it Item 1.}
Denote $m:=\inf\supp\ph$ and $M:=\sup\supp\ph$.
It is easy to see that,
in view of (\ref{avnonzero}),
the function $\ph^N_t$ has the form
$$
\ph^N_t=\mI_{[M,m+Nt]} + r^N_t,
$$
where $|r^N_t|\leq C$ and the Lebesgue measure Leb$(\supp r^N_t)\leq C_1$,
with constants $C,C_1$ independent from $N$
(see figure \ref{pic:phi}).
Then
\bee \lbl{Varna2}
\Var\SSS_{\ph^N_t}=\Var\big(\SSS_{\mI_{[M,m+Nt]}} + \SSS_{r^N_t} \big)
=\Var\SSS_{\mI_{[M,m+Nt]}} + \Var\SSS_{r^N_t} + 2 \Cov (\SSS_{\mI_{[M,m+Nt]}} ,\SSS_{r^N_t}).
\eee
In view of (\ref{iVar}),
$\Var\SSS_{\mI_{[M,m+Nt]}} = \pi^{-2}\ln N + O(1).$
Clearly,
$\Var\SSS_{r^N_t}\leq C$, where $C$ is independent from $N$.
Then the desired relation follows from (\ref{Varna2}) joined with the
Cauchy-Bunyakovsky-Schwartz inequality
$$
\big|\Cov (\SSS_{\mI_{[M,m+Nt]}} ,\SSS_{r^N_t}) \big|
\leq \sqrt{\Var\SSS_{\mI_{[M,m+Nt]}} } \sqrt{\Var\SSS_{r^N_t}}.
$$

{\it Item 2.}
To get the  desired result
it suffices to note that
assumptions of \rtheo \ref{th:CLT}
are satisfied for
the family of functions $\ph_{t_1}^N,\ldots,\ph_{t_d}^N$,
with $V_N=\pi^{-2}\ln N$
and the covariance matrix $(b_{ij})$ from (\ref{iCov}).
Indeed, estimate (\ref{th|h|}) is obvious
since the functions $\ph_{t_i}^N$ are bounded uniformly in $N$.
Assumption (\ref{thcov}) follows from the logarithmic growth of the variance
by the argument similar to that used in the proof of \rprop\ref{lem:fd}.

\ssk
{\it Item 3.}
We follow the same strategy as when proving \rtheo \ref{th:FCLTS}.
We set
\bee\non
f_\ph^N:=\fr{1}{\tau}\ili_{0}^{\tau} \ph^N_s\,ds
\qnd
g_{\ph,t}^N:=\ili_0^t \ph^N_s \,ds - \fr{t}{\tau}\ili_0^{\tau} \ph^N_s \,ds.
\eee
Then we have
$$
\eta^N=\fr{\SSS_{f_\ph^N}-\MO \SSS_{f_\ph^N}}{\pi^{-1}\sqrt{\ln N}}
\qmb{and}\qu
z_t^N=\SSS_{g_{\ph,t}^N}-\MO \SSS_{g_{\ph,t}^N}.
$$
Take any $0\leq t_1<\ldots<t_d\leq 1$.
We claim that the functions
$f_\ph^N,g^N_{\ph,t_1},\ldots,g^N_{\ph,t_d}$
satisfy assumptions of \rtheo \ref{th:CLTj}
with
$V_N=\pi^{-2}\ln N$, $R_N=N$,
$b_{11}^f=1/2$ and
the functions $\bm g^\infty_{\ph,t_i}=g_{t_i}$,
where the $g_{t_i}$ are defined in (\ref{mainfg}).
Indeed, note that
$$
\ph_s^N=\ph*\mI_{[0,sN]}.
$$
Consider the functions $f^N$ and $g^N_t$, defined above (\ref{mainfg}).
We have
$$
f_\ph^N=\ph*\Big(\fr{1}{\tau}\ili_{0}^{\tau} \mI_{[0,sN]} \, ds\Big) = \ph*f^N.
$$
Similarly,
\bee\lbl{starrrrrr}
g^N_{\ph,t}= \ph*g_t^N.
\eee
As it was shown in the proof of \rtheo \ref{th:FCLTS},
the functions
$f^N,$ $g_{t_i}^N$
satisfy assumptions of \rtheo\ref{th:CLTj}
with $V_N$, $R_N$, $b_{11}^f$ as above and $\bm g^\infty_{t_i}=g_{t_i}$.
Then, due to Example \nolinebreak \ref{ex:conv},
the functions $f_\ph^N,g^N_{\ph,t_i}$ fulfil assumptions \emph{f.2}-\emph{g.2},
with the same $V_N$, $R_N$, $b_{11}^f$ and $\bm g^\infty_{\ph,t_i}=\bm g^\infty_{t_i}$.
To show that assumption \emph{f.1} is satisfied as well,
it suffices to prove that
$$
\Var \SSS_{f_\ph^N} = \fr{1}{2\pi^2}\ln N + O \big(\sqrt{\ln N}\big).
$$
In view of item {\it 1} of the theorem, this can be shown
by the argument used in the proof of \rprop\ref{lem:fvar}(1).
Now \rtheo\ref{th:CLTj} joined with \rprop\ref{lem:fvar}(3) implies the convergence
\bee\lbl{fdierg}
\DD(\eta^N,z^N_{t_1},\ldots,z^N_{t_d}) \raw \DD(\eta,z_{t_1},\ldots,z_{t_d}) \ass N\ra\infty,
\eee
where the random variable $\eta$ and the process $z_t$
are as in the formulation of \rtheo\ref{th:FCLTS}.

Next we show that the family of measures
$\{\DD(\eta^N,z^N), N\in\mN\}$
is tight. To this end, as in \rtheo\ref{th:FCLTS},
it suffices to prove that
the family of functions $g^N_{\ph,t}$ satisfies assumption (\ref{contr})
of \rprop \ref{lem:regtight}. The first estimate from (\ref{contr})
is obvious since $g^N_{\ph,0}=0$.
In view of the identity $\hat g_{\ph,t}^N= \hat\ph\hat g_t^N$
which follows from (\ref{starrrrrr}),
we have
$$
\|g_{\ph,t}^N-g_{\ph,s}^N\|_{1/2}\leq \|\hat\ph\|_\infty \|g^N_t-g^N_s\|_{1/2}.
$$
Since $\ph\in L^1(\mR)$, the norm $\|\hat\ph\|_\infty$ is finite.
Then it remains to establish the second estimate from (\ref{contr}) for
the functions $g^N_t$. But it was already done
in the proof of \rtheo\ref{th:FCLTS}.

Now, literally repeating arguments from {\it Step 3} of the proof of \rtheo\ref{th:FCLTS},
we see that the convergence of finite-dimensional distributions (\ref{fdierg})
together with the tightness of the family of measures
$\DD(\eta^N,z^N)$ implies the desired convergence
$\DD(\eta^N,z^N)\raw \DD(\eta,z)$.

\ssk
{\it Item 4.}
The proof literally repeats that of \rprop \ref{lem:i2}. The rate of convergence of the cumulants $A^N_2$ and $A_2$
in this case is different with that from \rprop \ref{lem:i2} because of the correction $O(\sq{\ln N})$ in item {\it 1} of the theorem (cf. (\ref{iVar})).
\qed

\msk

{\bf Proof of Proposition \ref{lem:covnotfixed}}

Set $f_t:=\mI_{[at,bt]}$ and $f^N_t(x):=f_t(x/N)$. 
Then we have 
$\#_{[aNt,bNt]}=\SSS_{f^N_t}$. 
Due to Corollary~\ref{lem:VarS},
$$
\Cov(\SSS_{f^N_t},\SSS_{f^N_s})=
\fr{1}{4\pi^2}
\Re\Big(
2\ili_{|u|\geq 2} \hat f^N_t(u)\ov{\hat f^N_s(u)}\,du
+ \ili_{|u|<2} |u| \hat f^N_t(u)\ov{\hat f^N_s(u)}\,du
\Big)=:\fr{1}{4\pi^2}(I_1+I_2).
$$
Note that $\ds{\hat f^N_t(u)=\fr{e^{-iuaNt}-e^{-iubNt}}{iu}}$.
Let us study the integral $I_1$. We find
$$
I_1=2\ili_{|u|\geq 2}\frac{\cos\big(N(at-as)u\big) + \cos\big(N(bt-bs)u\big) - \cos\big(N(at-bs)u\big) -\cos\big(N(bt-as)u\big)}{|u|^2}\,du.
$$
Due to the assumption $t\neq s$, $at\neq bs$ and $bt\neq as$, the arguments of the cosines do not vanish, so they oscillate fast for $N\gg 1$. Since the integral $I_1$ converges absolutely for every $N$, this implies $I_1\to 0$ as $N\to\infty$.
Let us turn to the integral $I_2$. We have
$$
I_2=2\ili_0^2\frac{\cos\big(N(at-as)u\big) + \cos\big(N(bt-bs)u\big) - \cos\big(N(at-bs)u\big) -\cos\big(N(bt-as)u\big)}{u}\,du.
$$
Changing the variable $v:=Nu$, we find
$$
I_2=2\ili_0^{2N}\frac{\cos\big((at-as)v\big) + \cos\big((bt-bs)v\big) - \cos\big((at-bs)v\big) -\cos\big((bt-as)v\big)}{v}\,dv.
$$
Due to Proposition~\ref{lem:cos}, we obtain
$$
I_2\to 2(\ln|at-bs| + \ln|bt-as| -\ln|at-as|-\ln|bt-bs|).
$$
\qed.

\section{Main order asymptotic for  determinantal processes with logarithmically growing variance}
\lbl{sec:FCLTw}


In this section we prove a generalized version of  \rprop  \ref{ith:FCLTw}
for an important class of determinantal processes.
The latter includes processes with logarithmically growing variance,
e.g. the sine, Bessel and Airy processes.

Let $h^N:[0,1]\times \mR^m\mapsto\mR$ be a family of Borel measurable bounded functions with compact supports.
Consider the linear statistics
$$
\SSS_{h_t^N}:=\sli_{x\in\XX} h^N(t,x)
$$
as a random variable under a determinantal process given by a Hermitian kernel $\KK^N$.
Denote by $\Var_N,\,\Cov_N$ and $\MON$ the corresponding variance, covariance and  expectation.
\bpp\label{th:FCLTw}
Assume that there exists a sequence
$V_N\ra\infty$ as $N\ra\infty$,
$V_N>0$,
such that the following
three conditions hold.
\begin{enumerate}
\item
There exists a constant $C$ such that for any $N$ and almost all $t\in [0,1]$ we have
\bee\lbl{1cond}
\fr{\Var_N\SSS_{h_t^N}}{V_N}\leq C.
\eee

\item
There exists $b\in\mR$ such that
for almost all
$(t,s)\in [0,1]^2$ we have
\bee\label{b}
\fr{\Cov_N (\SSS_{h_t^N},\SSS_{h_s^N})}{V_N}\os{N\ra\infty}\ra b.
\eee

\item
We have
\bee\non
\|h^N\|_\infty = o(\sqrt {V_N}).
\eee
\end{enumerate}
Denote
$$
\xi^N_t=\fr{\SSS_{h_t^N}-\MON\SSS_{h_t^N}}{\sqrt{V_N}}, \quad 0\leq t \leq 1.
$$
Then for any functions $\phi_1,\ldots,\phi_n\in L^1 [0,1]$, $n\geq 1$, we have
\bee\lbl{fxiN}
\DD\Big(\ili_{0}^{1} \phi_1(t)\xi^N_t \, dt,\ldots, \ili_{0}^{1} \phi_n(t)\xi^N_t \, dt\Big) \os{N\ra\infty}\raw
\DD\Big(\eta\ili_{0}^{1} \phi_1(t)\,dt,\ldots,\eta\ili_{0}^{1} \phi_n(t)\,dt \Big),
\eee
where $\eta\volna\NN(0,b)$.
\epp
The principal assumption of \rprop \ref{th:FCLTw} is (\ref{b}).
In particular, it is satisfied for determinantal processes
with the logarithmic growth of the variance.
Let us explain this on the following examples.
Consider the sine or the Bessel process and
the linear statistics corresponding to the function
$$
h^N(t,x)=\mI_{[0,Nt]}(x),
\qmb{so that}\qu
\SSS_{h^N_t}=\#_{[0,tN]}.
$$
It is known that for $0<a<b$ we have
$$
\Var\#_{[aN,bN]}\volna C\ln N \ass N\ra\infty,
$$
for the both processes  (see (\ref{iVar}) for the sine-process and \cite{SoAB} for the Bessel process).
Then, literally repeating arguments from the proof of the convergence (\ref{iconvvar}),
we obtain (\ref{b}) with $b=1/2$ and $V_N=C\ln N$.
The same holds for the Airy process,
if one puts
$h^N(t,x)=\mI_{[-Nt,0]}(x)$
(so that
$\SSS_{h_t^N}=\#_{[-tN,0]}$)
and
$a<b<0$.

It can be checked that in the examples above
the other assumptions of \rprop \ref{th:FCLTw}
are satisfied as well,
so that the convergence (\ref{fxiN}) takes place.
In particular,
taking $n=1$ and $\phi_1=\mI_{[0,t]}$,
we get  the following corollary,
which is a version of the main order asymptotic
from \rtheo\ref{th:FCLTS} for the Airy and Bessel processes.
Set
$$
\xi^N_{A,t}
	=\fr{\#_{[-tN,0]}-\MO\#_{[-tN,0]}}{\sqrt{\Var \#_{[-tN,0]}}}
\qnd
\xi^N_{B,t}
	=\fr{\#_{[0,tN]}-\MO\#_{[0,tN]}}{\sqrt{\Var \#_{[0,tN]}}}.
$$
\begin{cor}\lbl{cor:AB}
Under the Airy processes for any $0\leq t\leq 1$ we have
$\DD(\int_0^t \xi_{A,s}^N \,ds)\raw \DD(\eta t)$ as $N\ra\infty$,
where $\eta\volna\NN(0,1/2)$.
Under the Bessel process, we have $\DD(\int_0^t \xi_{B,s}^N \,ds)\raw \DD(\eta t)$.
\end{cor}

Similar result holds true for the ergodic integrals
under the shift operator (studied in \rsec\ref{sec:iergint}).
Let $\ph:\mR\mapsto\mR$ be a bounded measurable function with compact support
satisfying $\ili_0^1\ph(s)\,ds=1$.
Set
$$
\ph_{A,t}^N:=\ili_0^{tN} \ph(\cdot+u)\,du
\qnd
\ph_{B,t}^N:=\ili_0^{tN} \ph(\cdot-u)\,du.
$$
Denote
$$
\xi_{A,\ph,t}^N:=\fr{\SSS_{\ph^N_{A,t}} - \MO\SSS_{\ph^N_{A,t}}}{\sqrt{\Var\SSS_{\ph^N_{A,t}}}},
$$
and define  $\xi_{B,\ph,t}^N$ in the same way.
\begin{cor}
Assertion of Corollary \ref{cor:AB} holds if replace the processes
$\xi_{A,s}^N$ and $\xi_{B,s}^N$
by the processes $\xi_{A,\ph,s}^N$ and $\xi_{B,\ph,s}^N$.
\end{cor}


\ssk
{\it Proof of Proposition \ref{th:FCLTw}.}
Let $\phi_i^N(x):=\ili_{0}^{1}\phi_i(t)h^N(t,x)\,dt$.
Using the Fubini theorem, for any $1\leq i\leq n$ we obtain
$$
\ili_{0}^{1} \phi_i(t)\xi^N_t \, dt=\fr{\SSS_{\phi_i^N}-\MON\SSS_{\phi_i^N}}{\sqrt{V_N}}.
$$
We claim that the family of functions $\phi_i^N$ satisfies assumptions of \rtheo \ref{th:CLT}. Indeed, the only condition fulfilment of which is
not obvious is (\ref{thcov}). Let us check it. In view of the Fubini theorem, we have
$$
\fr{\Cov_N(\SSS_{\phi_i^N},\SSS_{\phi_j^N})}{V_N}
=\MON\Big(\ili_{0}^{1} \phi_i(t)\xi^N_t \, dt
\ili_{0}^{1} \phi_j(s)\xi^N_s \, ds\Big)
=\ili_{0}^{1}\ili_{0}^{1} \phi_i(t)\phi_j(s)\fr{\Cov_N(\SSS_{h^N_t},\SSS_{h^N_s})}{V_N} \, dtds.
$$
Then, using the dominated convergence theorem, (\ref{1cond}) and (\ref{b}) we get
\bee\lbl{newcov}
\fr{\Cov_N(\SSS_{\phi_i^N},\SSS_{\phi_j^N})}{V_N}
\ra
b \ili_{0}^{1}\ili_{0}^{1} \phi_i(t)\phi_j(s)\, dtds =:c_{ij}
\ass N\ra\infty,
\eee
so that assumption (\ref{thcov}) is fulfilled with $b_{ij}=c_{ij}$.
Now it remains to apply \rtheo \ref{th:CLT}.
Indeed, since the limiting vector $\xi$ obtained there
is Gaussian with the covariance matrix $(c_{ij})$,
it coincides in distribution with the random vector from the right-hand side of (\ref{fxiN}).
\qed
\ssk

\bsk

{\bf Acknowledgements.}
We are deeply grateful to   Vadim Gorin, Alexei Klimenko, Gaultier Lambert,  Leonid Petrov,  Alexander Sodin  and Mikhail Zhitlukhin for useful discussions.
We are deeply indebted to the referees for the careful reading of the 
paper and the very helpful comments and suggestions.
Both authors are supported by the grant MD 5991.2016.1 of the President of the Russian Federation.
A. Bufetov's research has received funding from the European Research Council
(ERC) under the European Union's Horizon 2020 research and innovation
programme (grant agreement No 647133 (ICHAOS)). It has also been
funded by the Russian Academic Excellence Project `5-100'
and
by the Gabriel Lam\'e Chair at the Chebyshev
Laboratory of the SPbSU,
a joint initiative of the French Embassy in the Russian Federation
and the Saint-Petersburg State University.

\end{document}